\documentclass[10pt,a4paper]{amsart}

\input{ivan_style.sty}

\begin{document}
	\setstcolor{red}
	\title[Regular polytopes, sphere packings and Apollonian sections]{Regular polytopes, sphere packings and Apollonian sections} 
	
	\thanks{Partially supported by the CNRS and the Austrian Science Fund (FWF), projects F-5503 and P-34763.}
	
	\author[Iv\'an Rasskin]{Iv\'an Rasskin}
	\address{LIS, Aix-Marseille Université, CNRS, Marseille, France}
	\email{ivan.rasskin@lis-lab.fr}
	
	
	
	\begin{abstract} 
		In this paper, we explore the geometry and the arithmetic of a family of polytopal sphere packings induced by regular polytopes. We prove that every integral polytope is crystallographic and we show that there are 11 crystallographic regular polytopes in any dimension. After introducing the notion  of Apollonian section, we study which Platonic crystallographic packings emerge as cross sections of the Apollonian arrangements of the regular 4-polytopes. Additionally, we compute the Möbius spectrum of every regular polytope. 
	\end{abstract}
	
	\maketitle
	
		\section{Introduction}

Apollonian circle packings and their generalizations are currently active research areas in geometric number theory \cite{kapovich2023superintegral,bogachev2024kleinian}.  In dimension 2, certain variants of integral Apollonian packings have been explored by substituting the building block with different circle packing modeled on polyhedra \cite{guettler,stange2015bianchi,Zhang+2018+71+110,aporingpacks,chait2020taxonomy,RR21_1}. Although every polyhedron can be employed to construct a packing, not all of them possess an integral structure like the Apollonian one, which allows for packings where the bends (the inverses of the radii) of all the circles are integers. The fundamental question of determining which polyhedra are \textit{integral} in this sense remains wide open \cite{KontorovichNakamura,chait2020taxonomy}.

\medskip

\begin{figure}[H]
	\centering
	\includegraphics[align=c,width=.247\textwidth]{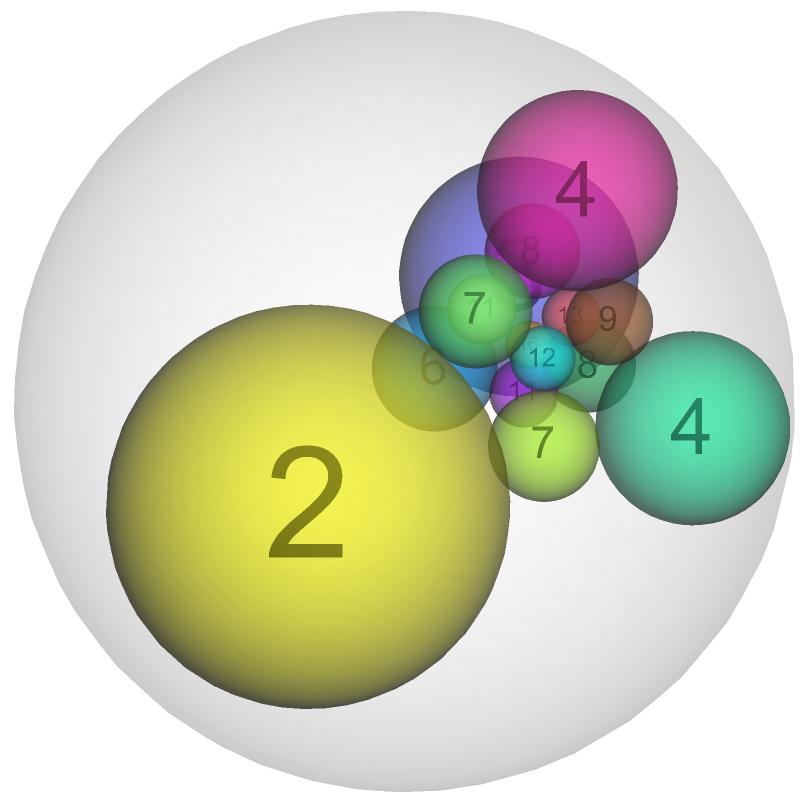}\includegraphics[align=c,width=.247\textwidth]{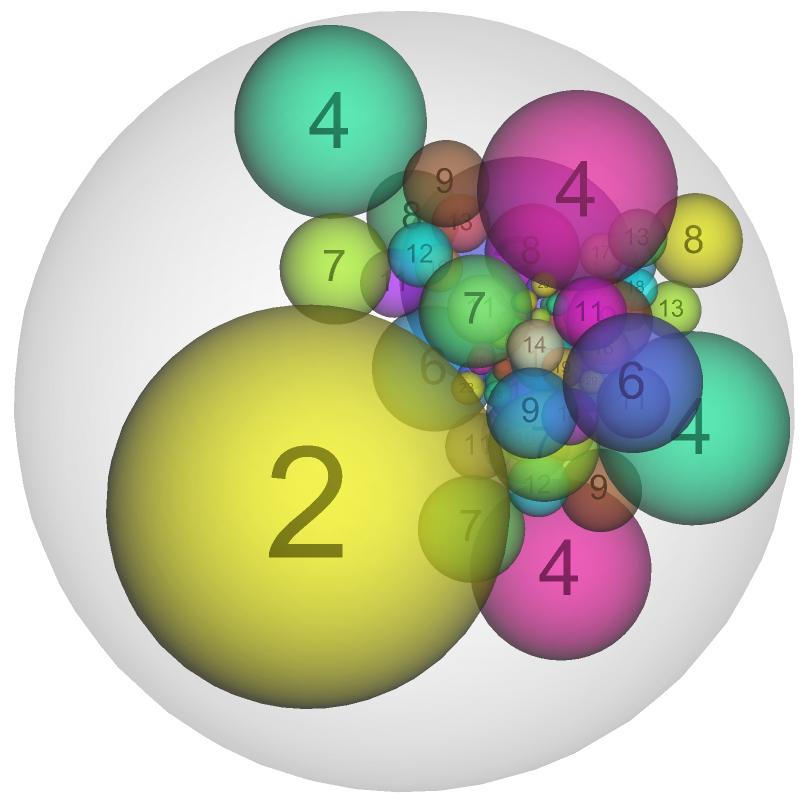}\includegraphics[align=c,width=.247\textwidth]{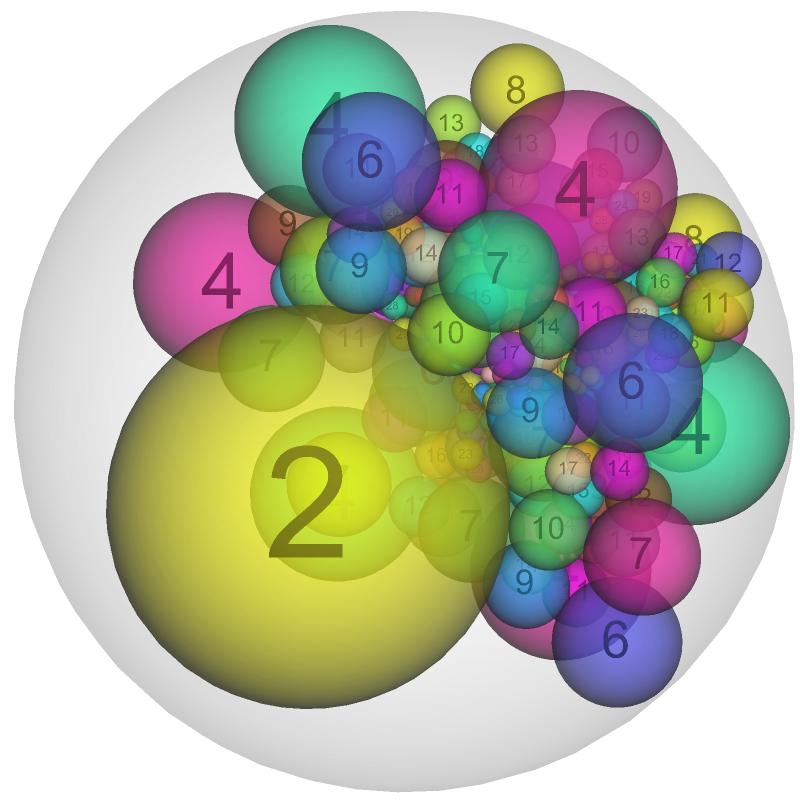}\includegraphics[align=c,width=.247\textwidth]{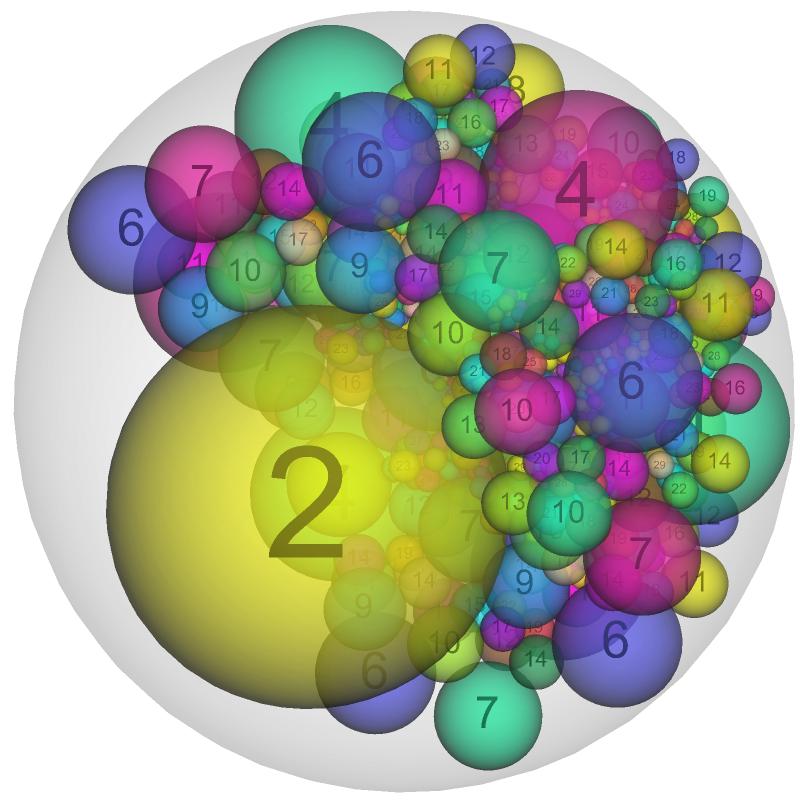}
	\vspace{-.2cm}
	\caption{
An integral hypercubic crystallographic packing after 0, 1, 2 and 3 iterations.	The numbers are the bends of the spheres.
	}
	\label{fig:hcube}
\end{figure}

Similarly, in dimension $3$, a family of crystallographic/Apollonian-like sphere packings arises by iteratively reflecting an initial sphere packing modeled on a $4$-polytope, as shown in Figure \ref{fig:hcube}. Integral crystallographic packings modeled on the $4$-simplex \cite{soddy1936, kontorovich2019soddy} and the orthoplex \cite{nakamura2014localglobal,Dias2014TheLP,SHEYDVASSER201941,RR2024links} have been extensively studied. However, unlike polyhedra, not every $4$-polytope is \textit{crystallographic}, meaning it may not serve as a suitable model for a crystallographic packing. 

\begin{figure}[H]
	\centering
	\begin{tabular}{cc}
		\includegraphics[align=c,width=0.45\textwidth]{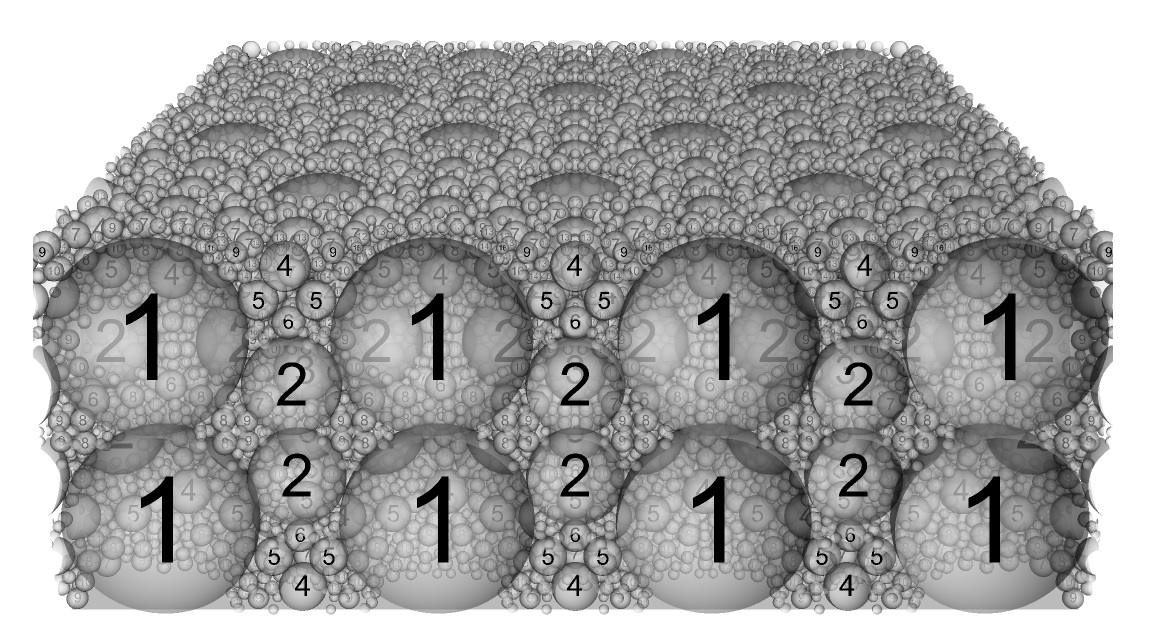}			&	\includegraphics[clip,trim=0 80 0 80,align=c,width=0.4\textwidth]{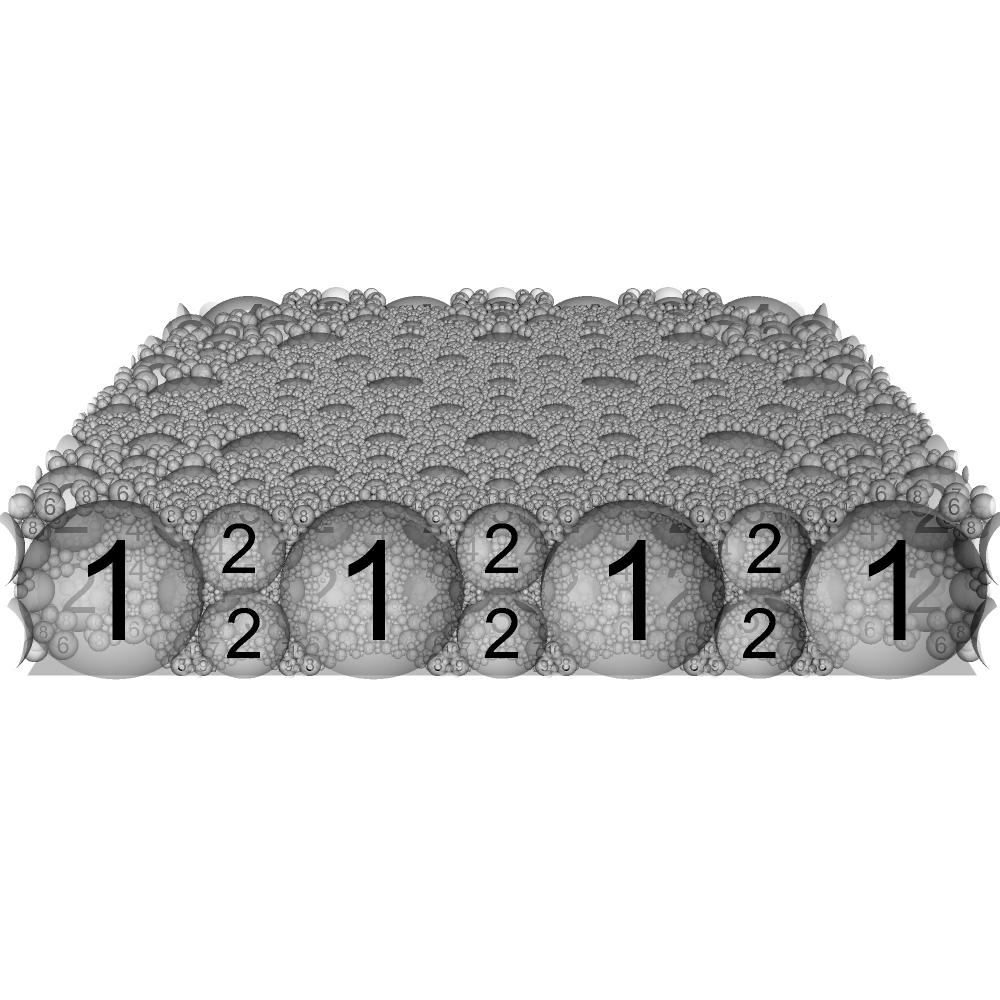}
	\end{tabular}	
	\caption{(Left) The integral hypercubic crystallographic packing $\mathscr P_{\{4,3,3\}}(0,0,1,2)$. (Right) The integral 24-cell crystallographic packing $\mathscr P_{\{3,4,3\}}(0,0,1,2)$.}
	\label{fig:conj24}
\end{figure}

One of the main questions regarding integral crystallographic packings is the \textit{local-to-global conjecture}, first introduced by Graham, Lagarias, Mallows, Wilks and Yan in \cite{apoGI} for Apollonian packings. The conjecture asserts that for any integral Apollonian circle packing, any sufficiently large number that avoids certain local modulo obstructions will be the bend of a circle in the packing. The local-to-global conjecture was disproven by Kertzer, Haag, Stange and Rickards in \cite{kertzer2024local}. This question has also been studied in the context of other integral crystallographic packings. For instance, Kontorovich demonstrated in \cite{kontorovich2019soddy} that the conjecture holds true for the family of integral simplicial crystallographic packings in dimension 3. Among the regular crystallographic sphere packings studied in this paper, two of them, depicted in Figure \ref{fig:conj24}, exhibit an interesting local-to-global behaviour, as their set of bends appears to have no local obstructions.
\begin{conj} The set of bends of the hypercubic $\mathscr P_{\{4,3,3\}}(0,0,1,2)$ and 24-cell $\mathscr P_{\{3,4,3\}}(0,0,1,2)$ crystallographic packings is $\mathbb N$.
\end{conj}
\vspace{-.5cm}

In this paper, we explore the crystallography and integrality of regular polytopes in any dimension. We present a necessary condition for the entries of the Gramian of the dual of an edge-scribable polytope to be integral (Lemma \ref{lem:integral}). This condition offers a straightforward method for proving the nonnintegrality of certain edge-scribable polytopes. 

The study of the cross-sections is a classic method for identifying patterns in crystallographic packings \cite{boyd,baragar2018higher}. We will introduce an algebraic approach based on the concept of \textit{Apollonian section}, which will be useful for identifing the Platonic cross-sections of the regular crystallographic packings in dimension 3. Some of these sections have been used as a geometric framework for deriving results in geometric knot theory \cite{RR2024links}.

There is an extensive family of different spectral invariants that have been prove to be useful for studying different algebraic and geometric graph properties. In this vein, Ramírez Alfonsín and the author defined in  \cite{RR21_1} a spectral invariant of edge-scribable polytopes called \textit{Möbius spectrum}. It is currently unknown if this invariant is complete or if it is related to other invariants. If so, it could be helpful for enumerating edge-scribable polytopes. In particular, in dimension $3$, since all $3$-polytopes are edge-scribable and are completely characterized by their graphs, the Möbius spectrum can be naturally extended to these graphs, usually called \textit{polyhedral graphs}. It would be interesting to investigate if co-spectral polyhedral graphs, with respect to the classical spectrum of the adjacency matrix, have different Möbius spectra. We conclude this paper by computing the Möbius spectrum of all the regular polytopes.

\subsection{Main contributions} 
\begin{enumerate}
	\item\label{Result1} We introduced a necessary condition for integral polytopes (Lemma \ref{lem:integral}), enabling us to establish that every integral polytope is crystallographic (Theorem \ref{thm:integral}).
	\item\label{Result2} We enumerate the 11 crystallographic regular polytopes in any dimension (Theorem \ref{thm:classification}) and determine their integrality.
	\item\label{Result3} We give a unified Descartes Theorem for the regular $4$-polytopes in terms of their Schläfli symbol (Theorem \ref{thm:4regdesth}). This theorem allows us to derive linear representations of the full symmetry groups (Corollary \ref{cor:lineargroup}) and the integrality conditions for constructing integral packings (Corollary \ref{cor:integrality}).
	\item\label{Result4} We show that the Platonic crystallographic packings can be obtained as cross sections of the Apollonian arrangements of the regular $4$-polytopes (Theorem \ref{thm:aposecs}) and study the integral invariance of these sections (Theorem \ref{cor:aposecsint}). 
	\item\label{Result5} We compute the Möbius spectrum of every regular polytope in dimension $d\ge3$ in terms of the number of vertices and the canonical length (Theorem \ref{thm:mobspec}). 
\end{enumerate}

\subsection{Organization of the paper} Section \ref{sec:preliminaries} presents all the necessary preliminaries on the classes of sphere packings and polytopes studied in this paper. In Section \ref{sec:crystalint}, we state and prove all our contributions. Finally, in Appendix \ref{sec:appendix} we explicit linear representations of the full symmetry groups, show integral packings/arrangements and describe the Platonic Apollonian sections for each regular $4$-polytope.
\medskip

\paragraph{\textbf{Acknowledgements:}} 
This work is partially included in the author's PhD thesis \cite{rasskin_thesis}. The author would like to thank his advisor Jorge Ram\'irez Alfons\'in, as well as, Benjamin Khlan, Thomas Serafini, Cesar Ceballos, Joseph Doolittle, Yannic Vargas, Michael Henry and Katherine Stange for their valuable suggestions and many fruitful discussions. The author also extends gratitude to the referee for their valuable insights and constructive criticisms, which greatly improved the presentation of this work.

\paragraph{\textbf{Figures:}} All the figures in this paper were done by the author's own software coded with Mathematica 13.1 \cite{Mathematica}.

\section{Preliminaries on sphere packings and regular polytopes}\label{sec:preliminaries}
In this section, we shall review some notions and definitions needed in the rest of the paper. We refer the reader to \cite{RR20,RR21_1,KontorovichNakamura,chait2020taxonomy} for more details. 
\subsection{Inversive coordinates}
An \textit{oriented hypersphere}, or simply \textit{sphere}, of $\wrd:=\mathbb{R}^d\cup\{\infty\}$, is the image of a spherical cap of $\mathbb S^d$ under the stereographic projection. Every sphere $S$ is uniquely defined by its center $c\in\wrd$ and its bend $b\in\mathbb R$ (the reciprocal of the \textit{signed} radius), or if $S$ is a half-space, by its normal vector $\widehat n\in \mathbb S^{d-1}$ pointing to the interior and the signed distance $\delta\in\ru$ between its boundary and the origin. The \textit{inversive coordinates} of $S$ are represented by the $(d+2)$-dimensional real vector
\begin{align}\label{eq:invcoord}
	\mathbf i(S)=
\begin{cases}
		(bc,\dfrac{\overline b-b}{2},\dfrac{\overline b+b}{2})^T&\text{ if }b\not=0,  \\
		\quad\\
		(\widehat n,\delta,\delta)^T&\text{otherwise}. \\
	\end{cases}
\end{align}
where $\overline{b}=b\|c\|^2-\frac{1}{b}$ is the \textit{co-bend} of $S$. The co-bend is the bend of $S$ after inversion through the unit sphere. The \textit{inversive product} of two spheres $S,S'$ of $\wrd$ is the real value
\begin{align}\label{eq:invprod}
	\langle S,S'\rangle= \mathbf i(S)^T\mathbf Q_{d+2}\mathbf i(S)
\end{align}
where $\mathbf Q_{d+2}$ is the diagonal matrix $\mathrm{diag}(1,\ldots,1,-1)$ of size $d+2$. The inversive product encodes the relative position of $S$ and $S'$ according to the following criteria:
\begin{align}\label{eq:invprodvalues}
	\langle S,S'\rangle
\begin{cases}
	<-1&\text{if } S\cap S'=\emptyset, \\
	=-1&\text{if  $\partial S$ and $\partial S'$ are tangent and } \mathrm{int}(S)\cap \mathrm{int}(S')=\emptyset,  \\
	=1&\text{if  $\partial S$ and $\partial S'$ are tangent and }  S\subseteq S'\text{ or }S'\subseteq S,  \\
	>1&\text{if } \partial S\cap \partial S'=\emptyset\text{ and } S\subset S'\text{ or }S'\subset S.  \\
	\end{cases}
\end{align}
An arrangement of spheres $\mathcal A$ in $\wrd$, possibly infinite, is a \textit{packing} if their interiors are mutually disjoint. The \textit{Gramian} of a finite arrangement $\mathcal A=(S_1,\ldots,S_n)$ is the matrix $\mathrm{Gram}(\mathcal{A})=(\langle S_i,S_j\rangle)_{1\le i,j\le n}$.  The group of Möbius transformations of $\wrd$ preserves the inversive product and acts linearly on the inversive coordinates. It acts as an orthogonal subgroup of $\mathrm {SL}_{d+2}(\mathbb R)$ with respect to $\mathbf Q_{d+2}$. In particular, the inversion through a sphere $S$ transforms the inversive coordinates through left multiplication with the matrix
\begin{align}\label{eq:invmatrix}
	\mathbf S= \mathbf{I}_{d+2}-2\mathbf i(S)^T \mathbf i(S)\mathbf{Q}_{d+2}
\end{align}
where $\mathbf{I}_{d+2}$ is the identity matrix of size $d+2$.  Consider an arrangement of spheres $\mathcal A=(S_1,\ldots,S_{d+2})$ in $\wrd$ such that the matrix $\mathbf A=(\mathbf i(S_1),\ldots,\mathbf i(S_{d+2}))$ has full rank. Let  $\mathbf b=(b_1,\ldots,b_{d+2})^T$ be the bend vector of $\mathcal A$. Then, the inversion through $S$ transforms $\mathbf b$ via left multiplication with the so-called \textit{bend matrix} \cite{chait2020taxonomy}
\begin{align}\label{eq:bendmatrix}
\mathbf B=(\mathbf A^{-1} \mathbf S\mathbf A)^T.
\end{align} 

\subsection{Polytopal sphere packings}
The \textit{polar}  of a subset $X\subset\eud$  is the subset 
$X^*=\{u\in\eud\mid \langle u, v\rangle \leq 1\text{ for all }v\in X\}$. The \textit{stereographic sphere} of a point $v\in \mathbb R^{d+1}$ \textit{outside} $\sd$ (i.e. with $\|v\|>1$) is the sphere $S_v$ of $\wrd$ obtained by the stereographic projection of the spherical cap $\{-v\}^*\cap \mathbb S^d$. The inversive coordinates of $S_v$ can be computed from the Cartesian coordinates $\mathbf c(v)$ of $v$ using the equation
\begin{align}\label{eq:coordaffinv}
	\mathbf{i}(S_v)=(\|v\|^2-1)^{-1/2}(\mathbf c(v),1).
\end{align}

\begin{figure}[H]
	\includegraphics[width=.3\linewidth]{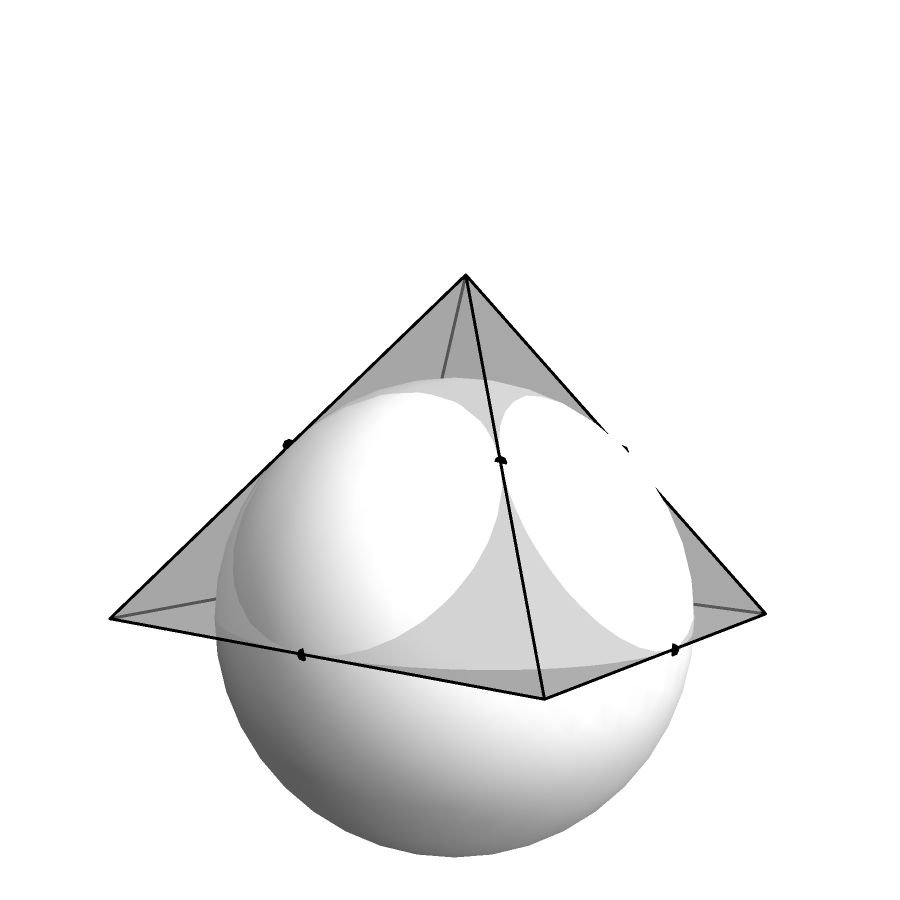}\includegraphics[width=.3\linewidth]{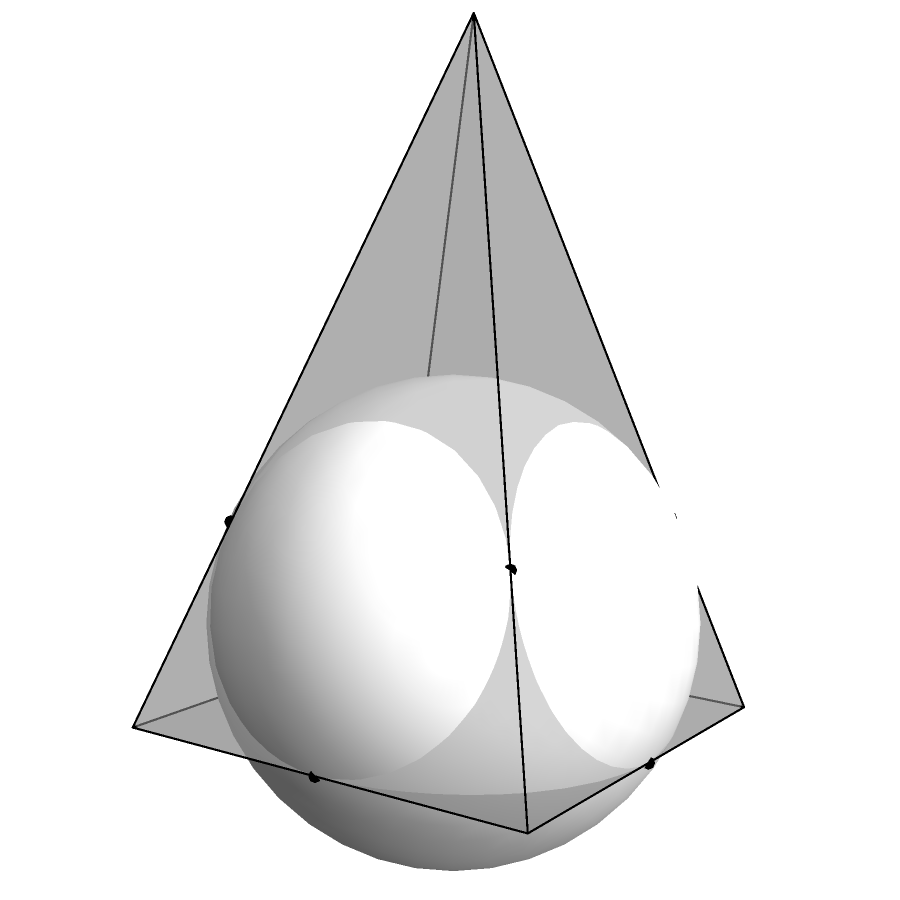}
	\vspace{-.3cm}
	\caption{An edge-scribed realization (left) and a canonical realization (right) of a $4$-pyramid.} 
	\label{fig:canonical}
\end{figure}

\begin{figure}[H]
	\centering
	\begin{tabular}{ccc}
		\includegraphics[trim=10 15 10 50,clip,height=4cm,align=c]{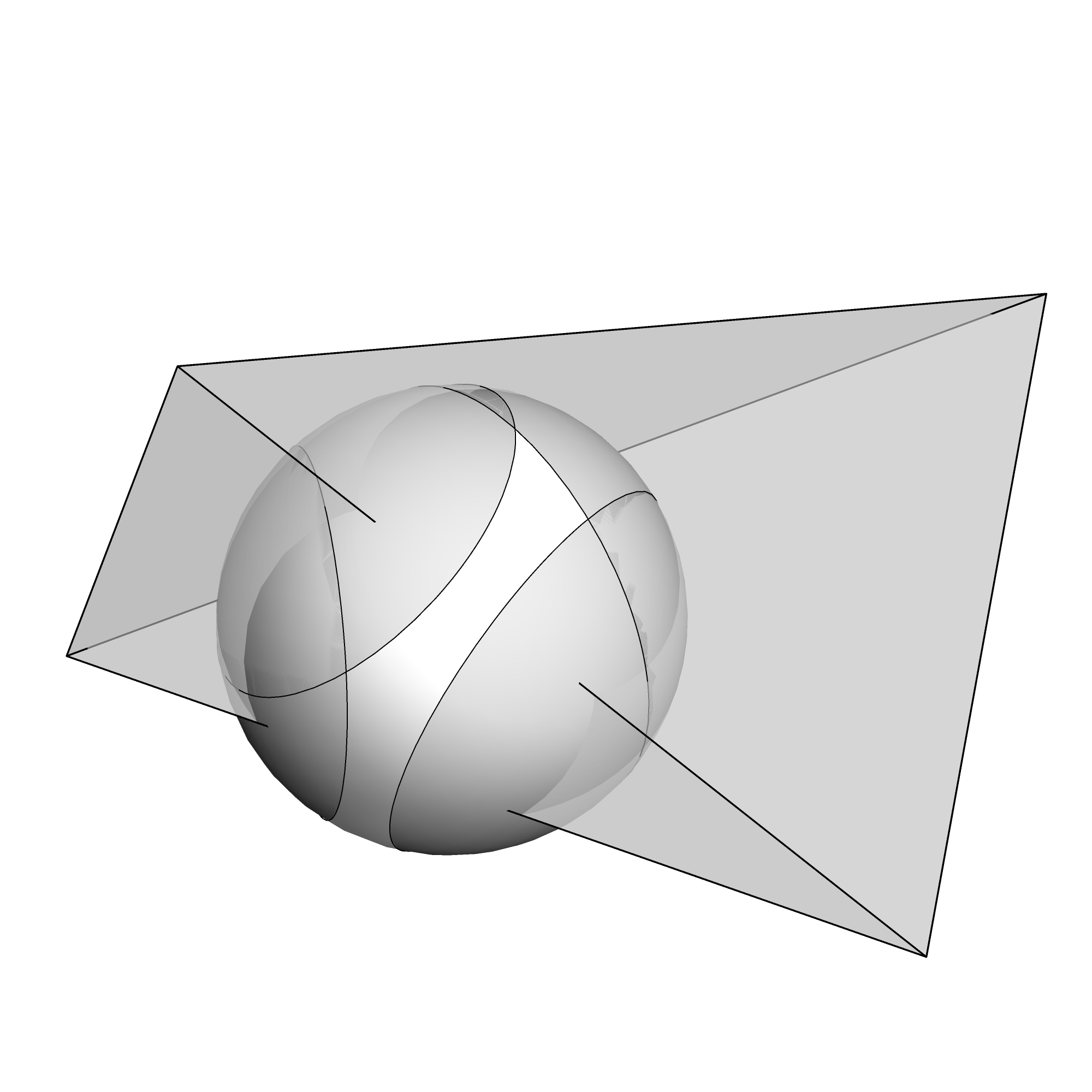}&\includegraphics[height=2.6cm,align=c]{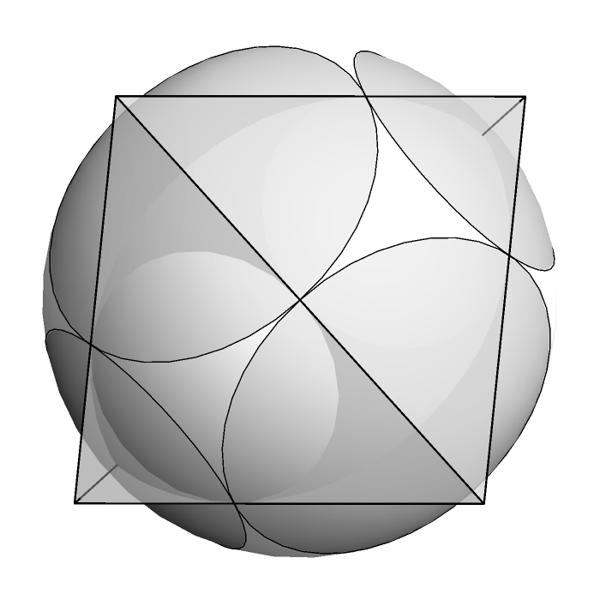}	 &\includegraphics[trim=0 0 20 30,clip,height=3cm,align=c,vshift=.1cm]{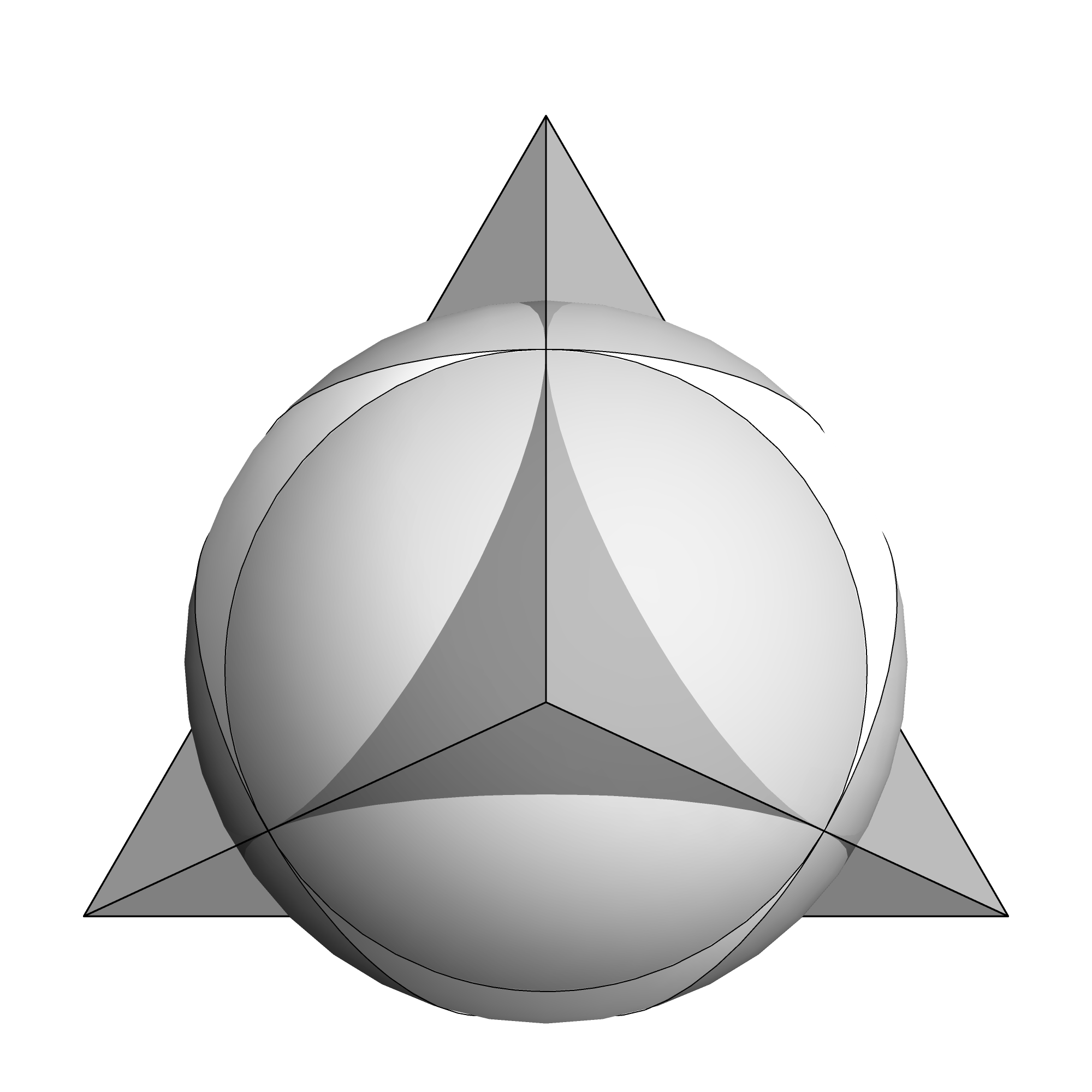} \\[-1.1cm]
		\reflectbox{\rotatebox[origin]{-150}{\begin{tikzpicture}[scale=.36] 
\clip (-8,-\prange) rectangle (5,8);
\draw[very thin, draw=black, fill=darkgray, fill opacity=\opac] (-2.062,0.2840) circle (1.689cm);
\draw[very thin, draw=black, fill=darkgray, fill opacity=\opac] (-2.754,3.696) circle (3.989cm);
\draw[very thin, draw=black, fill=darkgray, fill opacity=\opac] (0.4923,0.03083) circle (0.7515cm);
\draw[very thin, draw=black, fill=darkgray, fill opacity=\opac] (0.9696,-2.644) circle (2.733cm);

\end{tikzpicture}}}&\begin{tikzpicture}[scale=1] 
\begin{scope}[scale=2.1547]

\draw[very thin, draw=black, fill=darkgray, fill opacity=\opac] (0.4641,-0.2679) circle (0.4641cm);
\draw[very thin, draw=black, fill=darkgray, fill opacity=\opac] (-0.4641,-0.2679) circle (0.4641cm);
\draw[very thin, draw=black, fill=darkgray, fill opacity=\opac] (0,0.5359) circle (0.4641cm);
\draw[very thin, draw=black, fill=darkgray, fill opacity=\opac] (0.9282,0.5359) circle (0.4641cm);

\end{scope}

\end{tikzpicture}& \begin{tikzpicture}[scale=1] 
\begin{scope}[scale=2.1547]

\draw[very thin, draw=black, fill=darkgray, fill opacity=\opac] (0.4641,-0.2679) circle (0.4641cm);
\draw[very thin, draw=black, fill=darkgray, fill opacity=\opac] (-0.4641,-0.2679) circle (0.4641cm);
\draw[very thin, draw=black, fill=darkgray, fill opacity=\opac] (0,0.5359) circle (0.4641cm);
\draw[very thin, draw=black, fill=darkgray, fill opacity=\opac] (0,0) circle (0.07180cm);

\end{scope}

\end{tikzpicture}\\[-.5cm]
	\end{tabular}
	\caption{
		(Top) Three polyhedra with the spherical caps corresponding to their vertices.
		(Below) The arrangement projection of the three polyhedra. The last two are packings but only the third one is polytopal. 	}
	\label{fig:polytopal}
\end{figure}

Let $\P$ be a $d$-polytope with vertices outside the unit sphere. The \textit{arrangement projection} of $\P$ is defined as the arrangement $\mathcal A_\P$ of stereographic spheres of the vertices of $\P$ (see Figure \ref{fig:polytopal}). A $d$-polytope is termed \textit{edge-scribed} if its edges are tangent to the unit sphere. If, in addition, the barycenter of the contact points is the origin, it is referred to as \textit{canonical} (see Figure \ref{fig:canonical}). Canonical realizations are unique up to Euclidean isometries \cite{springborn2005unique}. A $d$-polytope is \textit{edge-scribable} if it admits an edge-scribed realization. In dimension $d\ge3$, all the edge-scribed realizations of an edge-scribable $d$-polytope $\P$ are equivalent up to Möbius transformations to a unique canonical realization $\P_0$ (see \cite{RR21_1} for more details).  

\begin{def1}\label{def:polytopalpacks}
	For every $d\ge2$, a sphere packing $\S_\P$ in $\wrd$ is \textit{polytopal} if there is an edge-scribable $(d+1)$-polytope $\P$ and a Möbius transformation $\mu$ such that $\S_\P=\mu( \mathcal A_{\P_0})$, where $\P_0$ is a canonical realization. 
\end{def1}
The combinatorial structure of a polytopal sphere packing $\S_\P$ is encoded by the corresponding edge-scribable polytope  $\P$. The vertices and the edges of $\P$ are in bijection to the spheres and the tangency relations of $\S_\P$. The facets of $\P$ correspond to the \textit{dual spheres} of $\S_\P$ which are the spheres forming the \textit{dual arrangement} $\S_\P^*:=\mu(\mathcal A_{\P_0^*})$. The \textit{Apollonian arrangement} of $\S_\P$ is defined as the orbit space $\mathscr{P}(\S_\P):= \langle \S_\P^*\rangle\cdot \S_\P$ where  $\langle \S_\P^*\rangle$ denotes the group generated by inversions through the dual spheres. All the Apollonian arrangements of $\P$ are equivalent up to Möbius transformations \cite{RR21_1}. 

\subsection{Crystallographic and integral polytopes}
In their work \cite{KontorovichNakamura}, Kontorovich and Nakamura introduced the concept of \textit{crystallographic packing}, which generalized families of highly symmetric sphere packings previously studied by Boyd \cite{boyd} and Maxwell \cite{maxwell}. These packings are dense infinite sphere packings derived from the limit set of a finitely generated hyperbolic reflection group. According to the Structure Theorem presented in \cite{KontorovichNakamura}, every crystallographic packing $\mathscr P$ can be decomposed as the action $\mathscr P=\langle \widetilde {\mathcal{C}}\rangle \cdot \mathcal C$, where $\mathcal C$ is a finite sphere packing called the \textit{cluster}, $\langle \widetilde {\mathcal{C}}\rangle$ is a geometrically finite subgroup of the group of Möbius transformations generated by the inversions through a finite arrangement of spheres $\widetilde {\mathcal{C}}$, called the \textit{co-cluster}, satisfying that every sphere of $\mathcal C$ is disjoint, tangent or orthogonal to every sphere of $\widetilde {\mathcal{C}}$ (see Figure \ref{fig:apollonianclassic}). Therefore, if an Apollonian arrangement $\mathscr{P}(\S_\P)= \langle \S_\P^*\rangle\cdot \S_\P$ of a polytopal sphere packing is also a packing, then it is crystallographic with $\S_\P$ as the cluster and $\S_{\P}^*$ as the co-cluster. In this case, due to the Möbius uniquenes of polytopal sphere packings, \textit{all} the Apollonian arrangements of $\P$ are also packings. Apollonian arrangements of $3$-polytopes are packings, but this is not true in general in higher dimensions \cite{RR21_1}. 
\begin{def1}\label{def:crystalpol}
	For every $d\ge3$, an edge-scribable $d$-polytope $\P$ is \textit{crystallographic} if any Apollonian arrangement $\mathscr{P}(\S_\P)=\langle \S_{\P}^*\rangle \cdot \S_{\P}$ is a packing.
\end{def1}
A crystallographic packing is said to be \textit{integral} if the set of bends is in $\mathbb Z$. We extend this definition to the Apollonian arrangements (not necessarily packings) of edge-scribable polytopes. 
\begin{def1}\label{def:intlpol}
	For every $d\ge3$, an edge-scribable $d$-polytope $\P$ is \textit{integral}\footnote{The definition of \textit{integral polytope} used in this paper differs from the one commonly employed in combinatorics, which involves polytopes with integer vertex coordinates, also known as \textit{lattice polytopes}.} if it admits an integral Apollonian arrangement  $\mathscr{P}(\S_\P)$.
\end{def1}

\subsection{Regular polytopes} A polytope is \textit{regular} if its symmetry group acts transitively on the set of flags.  Let us recall the list of regular polytopes in each dimension greater than 2. All $2$-polytopes admits a regular realization. The \textit{Platonic solids}—namely, the tetrahedron $\T^3$, the octahedron $\O^3$, the cube $\C^3$, the icosahedron $\I^3$ and the dodecahedron $\D^3$— are the five regular $3$-polytopes. The $4$-simplex $\T^4$, the orthoplex $\O^4$, the hypercube $\C^4$, the 600-cell $\I^4$ and the $120$-cell $\D^4$ are five regular $4$-polytopes which can be thought as a $4$-dimensional analogue of the Platonic solids. The remaining regular $4$-polytope, the 24-cell $\R^4$ (the notation is not standard), completes the list of regular $4$-polytopes. Classic realizations of the regular 3- and 4-polytopes are available in \cite{coxeter1973regular}. For every $d\ge3$, we shall denote by $\T^{d}$, $\O^{d}$ and $\C^{d}$ the $d$-dimensional analogue of the tetrahedron, octahedron and cube given by the following canonical realizations:
\begin{align}\label{eq:coordinatesreg}
	\T^d=\sqrt{\tfrac{2(d+1)}{d-1}}\mathrm{conv}(v_0,\ldots,v_d),&& \O^d=\sqrt 2 \,\mathrm{conv}(\pm e_1,\ldots,\pm e_d),&&	\C^d=\tfrac{1}{\sqrt{d-1}}\, \mathrm{conv}(\pm e_1\pm\ldots\pm e_d),
\end{align}
where $e_i$ is the canonical $i$-th vector of $\mathbb R^d$, $v_0=-\sum_{j=1}^{d}\sqrt{\frac{1}{j(j+1)}} e_j$, $v_i=\sqrt{\frac{i}{i+1}}e_i-\sum_{j=i+1}^{d}\sqrt{\frac{1}{j(j+1)}} e_j$ for every $i=1,\ldots,d-1$ and $v_d=\sqrt{\frac{d}{d+1}}e_d$. It is well-known that in dimension 5 or above, these three families are the only regular polytopes.

 \medskip

The polar of a regular $d$-polytope containing the origin is also a regular $d$-polytope. In particular,
\begin{align}\label{eq:coordinatespol}
	(\T^d)^*=\frac{1-d}2\T^d,&&(\O^d)^*
	=\sqrt{\frac{d-1}{2}}\C^d,&&(\C^d)^*
	=\sqrt{\frac{d-1}{2}}\O^d.
\end{align}

 \medskip

The \textit{Schläfli symbol} of a regular $d$-polytope $\P$ is the symbol $\{p_1,\dots,p_{d-1}\}$ which encodes the local structure and fully characterizes $\P$ \cite{schulte04}. It is defined recursively as follows: for each $(i+1)$-face $f$ of $\P$ with $i=1,\ldots,d-1$, $p_i$ represents the number of $i$-faces of $f$ that contain a given $(i-2)$-face of $f$ (with the convention that a $(-1)$-face is the empty set, which is contained in every face). The Schläfli symbol of every regular $d$-polytope with $d\ge3$ can be found in Table \ref{tab:mobspec}.

\subsection{The fundamental symmetries and basis of regular polytopes}
\begin{figure}[H]

			\begin{tikzpicture}[scale=1.1]
		\begin{scope}[xshift=-3cm]
					\node 	{\includegraphics[width=5cm]{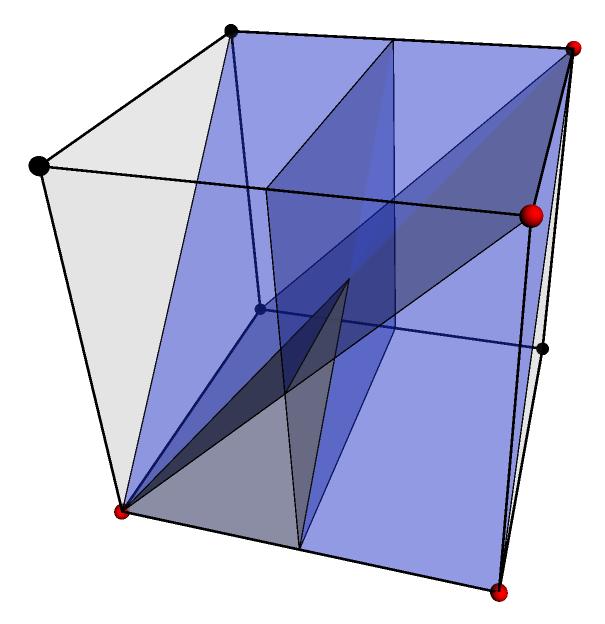}};
					\node at (-1.55,-1.75) {$v_1$};
					\node at (1.7,-2.3) {$v_2$};
					\node at (2.1,.7) {$v_3$};
					\node at (2.3,2.1) {$v_4$};
					\node at (0.4,1.4) {\color{white} $R_1$};
					\node at (1.3,1) {\color{white} $R_2$};
					\node at (1.2,-.5) {\color{white}  $R_3$};					
		\end{scope}		
				
		\begin{scope}[xshift=3cm,scale=.5]
\draw (-2.414,2.414) circle (2.414cm);
\draw (-.414,.414) circle (.414cm);
\draw (.414,-.414) circle (.414cm);
\draw (-.414,-.414) circle (.414cm);
\draw[red] (2.414,2.414) circle (2.414cm);
\draw[red] (2.414,-2.414) circle (2.414cm);
\draw[red] (-2.414,-2.414) circle (2.414cm);
\draw[red] (.414,.414) circle (.414cm);

\draw[thick,blue] (0,-5) -- (0,5) node at (.7,5.2) {$R_1$} ;
\draw[thick,blue] (-5,-5) -- (5,5) node at (5.6,4.8) {$R_2$} ;
\draw[thick,blue] (0,-1) circle (1.414cm) node at (1.9,-1) {$R_3$};			
		\end{scope}		
			\end{tikzpicture}
	\caption{(Left) A cube with a fundamental domain of its symmetry group (in dark gray), a fundamental basis ($v_1,v_2,v_3,v_4$) and the walls of the fundamental symmetries (in blue). 
	(Right) A cubic circle packing with the corresponding fundamental basis (in red) and the fundamental symmetries (in blue).
} 
	\label{fig:fundamentals}
\end{figure}

Let $\Phi=(f_0,\ldots,f_{d-1},f_{d}=\P)$ be a flag of a regular $d$-polytope $\P$. The simplex $\Delta_\Phi$, whose vertices are the barycenters of every  $f_i\in\Phi$, is a fundamental domain of the symmetry group $\mathfrak S(\P)$. This group is the finite Coxeter group generated by the reflections $r_1,\ldots,r_d$, which we call the \textit{fundamental symmetries} with respect to $\Phi$. Here, $r_i$ denotes the reflection through the hyperplane $R_i$ spanned by the facet of $\Delta_\Phi$ which is opposite to the barycenter of  $f_{i-1}\in\Phi$ (see Figure \ref{fig:fundamentals}). Since $\mathfrak S(\P)$ acts transitively on the set of flags, the group does not depend on the choice of the flag. We define the \textit{fundamental basis} of $\P$ with respect to $\Phi$, as the affine basis of vertices $(v_1,\ldots,v_{d+1})$ of $\P$, given by
\begin{align}
	v_1&:=f_0&&
	\text{and }&&	v_{k+1}:= r_1\cdots r_k(v_{k}) 
\end{align}
for every $k=1,\ldots,d.$
\subsection{The full symmetry group of the regular Apollonian arrangements} For each $d\ge2$, we will extend the definitions of fundamental symmetries and fundamental basis to any polytopal $d$-sphere packing $\S_\P$ derived from a regular  $(d+1)$-polytope $\P$, as illustrated in Figure \ref{fig:fundamentals}.  We also define a \textit{fundamental bend vector} of $\mathcal S_{\mathcal P}$ as the bend vector of a fundamental basis.

\medskip
Any Apollonian arrangement $\mathscr P(\mathcal S_\P)=\langle \S_\P^*\rangle\cdot \mathcal S_\P$ can be decomposed as $(\langle \S_\P^*\rangle\rtimes \mathfrak S(\S_\P))\cdot \{S_v\}$, where
$\mathfrak S(\mathcal S_\P)\simeq\mathfrak S(\P)$ is the symmetry group of $\S_\P$, and $S_v$ is the sphere corresponding to the vertex of a given flag $\Phi=(v=f_0,\ldots,f=f_d,\P)$. We call the group  $\Gamma(\S_\P):=\langle \S_\P^*\rangle\rtimes \mathfrak S(\S_\P)$ the \textit{full symmetry group} of the Apollonian arrangement $\mathscr P(\mathcal S_\P)$. This group can be seen as an analogue of the full symmetry group of the Apollonian-like packings defined by Baragar in \cite{baragar2018higher}. As $\mathfrak S(\P)$ is facet-transitive, $\Gamma(\S_\P)=\langle s_f\rangle\rtimes \mathfrak S(\S_\P)=\langle r_1,\ldots,r_{d+1},s_f\rangle$ where $r_1,\ldots,r_{d+1}$ are the fundamental symmetries of $\S_\P$ and $s_f$ is the inversion through the dual sphere $S_f$ corresponding to the facet $f\in\Phi$. In Figure \ref{fig:apollonianclassic}, we illustrate the classic Apollonian strip packing obtained as the orbit space of the Apollonian group $\langle \S_{\T^3}^*\rangle\cdot \mathcal S_{\T^3}$ and as the action of the full symmetry group $\Gamma(\S_{\T^3})\cdot \{S_v\}$.
\begin{figure}[H]
	\centering
	\begin{tabular}{cc}
		\includestandalone[align=c,scale=1.2]{tikzs/33apollonianmirrors}&\hspace{.5cm}	
		\includestandalone[align=c,scale=1.2]{tikzs/33standard}
	\end{tabular}
	\caption{The Apollonian strip packing  $\mathscr{P}_{\{3,3\}}$ obtained as crystallographic packing given by the action of the Apollonian group on a tetrahedral circle packing (left) and the full symmetry group $\Gamma_{\{3,3\}}$ on a single circle (right).
	}
	\label{fig:apollonianclassic}
\end{figure}

\subsection{The Platonic crystallographic packings}
\begin{figure}[H]
 	\centering
 	\begin{tabular}{cc}
 		\includegraphics[align=c,scale=1.2]{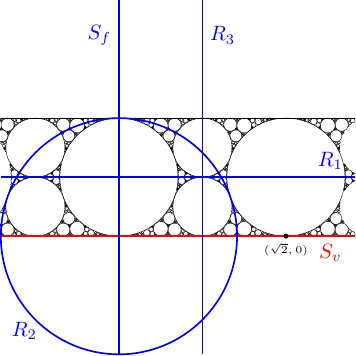}\hspace{.5cm}	&\includegraphics[align=c,scale=1.2]{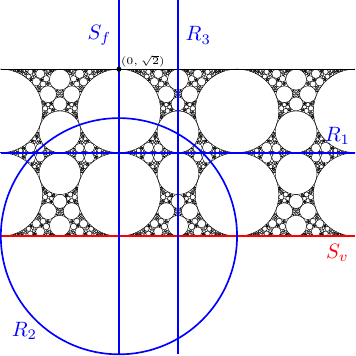}\\
 	\end{tabular}
 	\caption{The octahedral  $\mathscr{P}_{\{3,4\}}$ (left) and the cubic $\mathscr{P}_{\{4,3\}}$ (right) crystallographic packings.}
 	\label{fig:simpmirrors}
 \end{figure}

\begin{figure}[H]
 	\centering
 	\begin{tabular}{cc}
 		\includegraphics[align=c,scale=1.2]{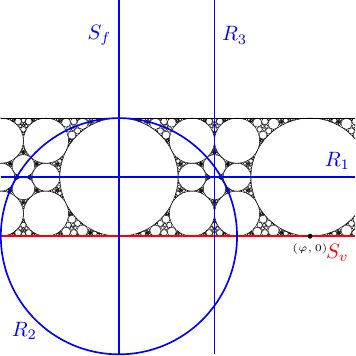}\hspace{.5cm}	&\includegraphics[align=c,scale=1.2]{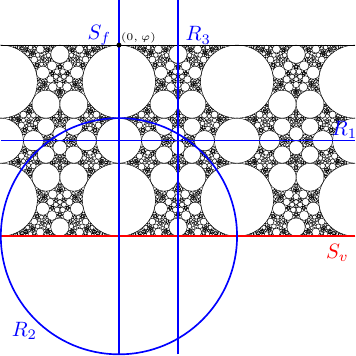}\\
 	\end{tabular}
 	\caption{The icosahedral  $\mathscr{P}_{\{3,5\}}$ (left) and the dodecahedral $\mathscr{P}_{\{5,3\}}$ (right) crystallographic packings.}
 	\label{fig:simpmirrors2}
 \end{figure}

The Apollonian packing depicted in Figure \ref{fig:apollonianclassic} is commonly known as the \textit{Apollonian strip packing} and serves as a canonical configuration for various purposes \cite{apoGI}. We extend this notion for every regular $(d+1)$-polytope $\P$ by stating that
a polytopal packing $\S_\P$ is \textit{strip} for a given flag $(v,\ldots,f,\P)$ if:
\begin{enumerate}
	\item The sphere $S_v\in \S_\P$ is the half-space $\{x_d\le0\}$.
	\item The dual sphere $S_f\in \S_\P^*$  is the half-space $\{x_1\le0\}$.
	\item The fundamental symmetry $r_2$ is the inversion through the unit sphere.
\end{enumerate}
For regular polytopes, strip packings are unique up to Euclidean isometries. We denote by $\mathscr P_{\{p_1,\ldots,p_d\}}$ and by $\Gamma_{\{p_1,\ldots,p_d\}}$ the Apollonian arrangement (up to Möbius transformations) and the full symmetry group of the regular polytope with Schläfli symbol ${\{p_1,\ldots,p_d\}}$. The Apollonian strip packing of Figure \ref{fig:apollonianclassic} corresponds to $\mathscr{P}_{\{3,3\}}$. In Figures \ref{fig:simpmirrors} and \ref{fig:simpmirrors2}, we illustrate the crystallographic strip packings $\mathscr{P}_{\{p,q\}}$  of the remaining Platonic solids, respectively. We refer to these five packings as the \textit{Platonic crystallographic packings.}

\section{Geometry and arithmetic of the regular polytopal sphere packings}\label{sec:crystalint}
An edge-scribable polytope $\P$ is crystallographic when the dihedral angles of $\P$, viewed as an hyperideal hyperbolic polytope, satisfy the \textit{crystallographic restriction} \cite{boyd}. This restriction dictates that the period of every rotation obtained as the product of two reflections through the facets is either $2,3,4,6,\infty$, imposing a condition on the dihedral angles. On the other hand, the dihedral angle $\alpha$ of two adjacents facets $f$ and $f'$ of $\P$ is equal to the \textit{intersection angle} of the corresponding dual spheres of $S_f,S_{f'}\in\S_\P^*$, as defined in \cite{RR20}. This angle can be computed from their inversive product by $\langle S_f,S_{f'}\rangle =\cos (\alpha)$. Therefore, the crystallographic restriction can be reformulated in terms of the inversive product of the dual spheres, as described in Lemma \ref{lem:crystal}.  
 \begin{lem}\label{lem:crystal} 	Let  $\P$ be  an edge-scribable $(d+1)$-polytope with $d\ge2$. Then  $\P$ is crystallographic if and only if, for any two dual spheres  $S_f$, $S_{f'}$ of any polytopal $d$-sphere packing $\S_\P$, we have
  $| \langle  S_f,S_{f'}\rangle |\in \lbrace
  \frac{\sqrt n}{2}\,|\,n\in\{0,1,2,3\}    \rbrace\cup[1,\infty)$.
 \end{lem} 
 
 In  \cite{chait2020taxonomy}, Chait-Roth, Cui and Stier studied the integrality of polyhedra with few vertices. To establish the integrality of edge-scribable polytopes, we will extend their methods.  We consider the group generated by the bend matrices of the dual inversions of any polytopal packing $\S_\P$. This group provides a linear representation $\langle \S_\P^*\rangle<\mathrm{SL}_{d+2}(\mathbb R)$ that acts on the set of bends of the Apollonian arrangement $\mathscr{P}(\S_\P)=\langle \S_\P^*\rangle\cdot\S_\P$. 
Notice that this representation does not depend on the packing but on the choice of the basis of $\S_\P$, which corresponds to the choice of an affine basis of vertices of $\P$.  
If all the bend matrices have integer entries, then  $\langle \S_\P^* \rangle<\mathrm{SL}_{d+2}(\mathbb Z)$ and the action of $\langle \S_\P^*\rangle$ on an initial integral packing $\S_\P$ produces an integral Apollonian arrangement. Alternatively, if the bend matrices have rational entries such that the set of denominators of the entries of all the matrices of $\langle \S_\P^*\rangle<\mathrm{SL}_{d+2}(\mathbb Q)$ is bounded, then any arrangement $\mathscr{P}(\S_\P)$ where the bends are rationals can be rescaled to obtain an integral arrangement.

\medskip

For proving nonnintegrality, similar ideas are employed. If there is a bend matrix with an irrational entry or there is a rational representation 
$\langle \S_\P^*\rangle<\mathrm{SL}_{d+2}(\mathbb Q)$ with unbounded denominators, then $\P$ is not integral. The authors of  \cite{chait2020taxonomy} demonstrated the unboundedness of the denominators of the entries for the powers of a matrix of  $\langle \S_\P^*\rangle<\mathrm{SL}_{d+2}(\mathbb Q)$ with infinite order. The following lemma provides a simpler method, akin to the characterization of integrality of Martin \cite{martin2024geometric} for Klenian arrengements induced by ideal class groups.

\begin{lem}\label{lem:integral}
	Let  $\P$ be  an edge-scribable $(d+1)$-polytope with $d\ge2$. If $\P$ is integral then, for any two dual spheres  $S_f$, $S_{f'}$ of any polytopal $d$-sphere packing $\S_\P$, we have $|\langle  S_f,S_{f'}\rangle|\in\lbrace \frac{\sqrt n}{2}\,|\,n\in\mathbb N\rbrace$.
\end{lem}

\begin{proof}
		We prove the lemma by contraposition. Let $\S_\P$ be a polytopal $d$-sphere packing with $d\ge2$ having two dual spheres $S_f$, $S_{f'}$  such that $|\langle S_f,S_{f'}\rangle|\not\in \lbrace\frac{\sqrt n}{2}\,|\,n\in\mathbb N\rbrace$. Let $\mathbf{M}\in\mathrm{SL}_{d+2}(\mathbb R)$ be the product of two bend matrices corresponding to the inversions through $S_f$, $S_{f'}$ with respect to some basis. If $\mathbf{M}$ has an irrational entry, then $\mathcal P$ is not  integral. Otherwise, $\mathbf{M}\in\mathrm {SL}_{d+2}(\mathbb Q)$ and in this case, we can combine the equations \eqref{eq:invprod}, \eqref{eq:invmatrix}, \eqref{eq:bendmatrix} with the hypothesis to obtain that
	\begin{align}\label{eq:tracerot}\mathrm {tr}(\mathbf{M})=d-2+4\langle S_f,S_{f'}\rangle^2\in\mathbb Q\setminus\mathbb Z.
	\end{align}

	Let us now suppose that the denominators of the entries of the powers of $\mathbf{M}$ are bounded. Hence, there is $r\in\mathbb N$ such that
	for every $n\in\mathbb N$, $r\mathbf M^n\in\mathrm{SL}_{d+2}(\mathbb Z)$. Then, the characteristic polynomial $\chi_{r\mathbf{M}^n}$ is a monic polynomial with integer coefficients. On the other hand, the eigenvalues of $r\mathbf{M}^n$ are $r\lambda_1^n\,\ldots,r\lambda_{d+2}^n$ where $\lambda_1,\ldots,\lambda_{d+2}$ are the eigenvalues of $\mathbf M$. Thus, for every $\lambda_i$ and for every $n\in\mathbb N$, we have that $\chi_{r\textbf{M}^n}(r\lambda_i^n)=0$, so $r\lambda_i^n$ is an algebraic integer. This implies that every $\lambda_i$ is also an algebraic integer. Hence, the coefficients of $\chi_M$, which include $\mathrm{tr}(\mathbf{M})$, are in $\mathbb Z$, a contradiction. Therefore, the denominators of the entries of the powers of $\mathbf M$ are unbounded, so $\P$ is not integral.	
\end{proof}

\begin{figure}[H]
	\centering
	\begin{tikzpicture}[scale=2.5] 
		\node at (-3.2,0) {	\includegraphics[align=c,width=.4\textwidth]{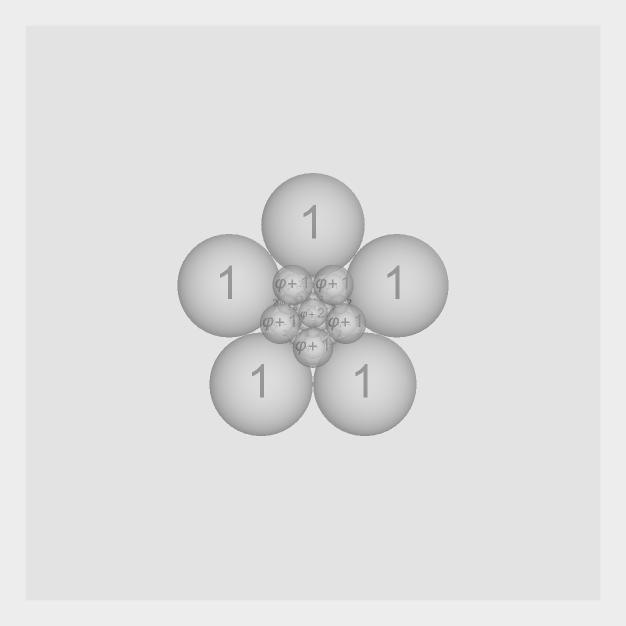}};
		\node at (0,0) {	\includegraphics[align=c,width=.4\textwidth]{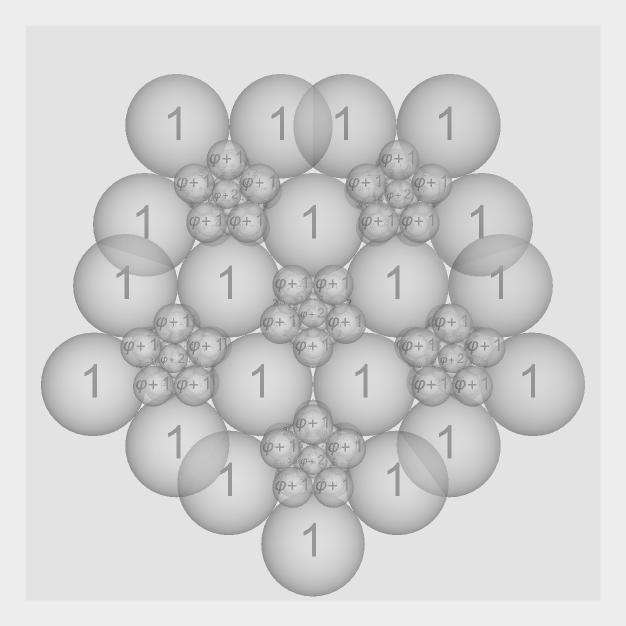}} ;
	\end{tikzpicture}
	\caption{
A $\mathbb Z[\varphi]$-integral polytopal strip packing of the 600-cell viewed from above and labelled with the bends (left) and the first iteration of its $\mathbb Z[\varphi]$-integral Apollonian arrangement having overlapping spheres (right).					
	}
	\label{fig:nonpack}
\end{figure}

Based on the work of Kontorovich and Nakamura in \cite{KontorovichNakamura},  Chait-Roth, Cui and Stier presented the list of uniform integral polyhedra (\cite{chait2020taxonomy}{Th. 26}). Lemma \ref{lem:integral} allows to easiliy identify an error in this list: the 6-prism is not integral because the Gramian of its dual contains the entry $-5/3\not\in\lbrace \frac{\sqrt n}{2}\,|\,n\in\mathbb N\rbrace$. Another straightforward consequence of the previous two lemmas is the following.
\begin{thm}\label{thm:integral}
Every integral polytope is crystallographic.
\end{thm}
Theorem \ref{thm:integral}  fails for integral polytopes defined in other number rings other than $\mathbb Z$. For instance, the 600-cell is integral in $\mathbb Z[\varphi]$, but is not crystallographic (see Figure \ref{fig:nonpack}). 

\subsection{The regular crystallographic packings}
Crystallographic polytopes exist only in dimension $3\le d\le 19$  \cite{bogachev2024kleinian}. In the case of regular polytopes, we have the following.
\begin{thm}\label{thm:classification}
	The only crystallographic regular polytopes are:
	\begin{enumerate}
		\item[$(d=3)$] $\T^3,\O^3,\C^3,\I^3,\D^3$.
		\item[$(d=4)$] $\T^4,\O^4,\C^4,\R^4,\D^4$.
		\item[$(d=6)$] $\O^6$.
	\end{enumerate}
	Moreover, all these are integral except $\I^3,\D^3,\D^4$ which are integral in $\mathbb Z[\varphi]$.
\end{thm}

\begin{proof}
	We start by regarding the crystallography. Every $3$-polytope is crystallographic. In dimension $4$, it can be easily checked that every regular $4$-polytope satisfies the conditions of Lemma \ref{lem:crystal} except the 600-cell. Let $\P$ be one of the three regular $d$-polytopes in the dimension $d\ge5$ and let $S_f,S_{f'}$ be two dual spheres of a polytopal $(d-1)$-sphere packing $\S_\P$ corresponding to two adjacent facets of $\P$. By combining equations \eqref{eq:coordaffinv}, \eqref{eq:coordinatesreg}, \eqref{eq:coordinatespol}, we obtain
\begin{align}\label{eq:prodadj}
	\langle S_f,S_{f'} \rangle=\left\{ \begin{array}{lcc}
		-\frac1{d-2}&&\text{ if }\P=\T^d\text{ or }\C^d,\\[.3cm]
		1-\frac2{d-2}&&\text{ if }\P=\O^d.
	\end{array}\right. 
\end{align}	
Therefore, the crystallographic restriction of Lemma \ref{lem:crystal} fails for every $\P\not=\O^{6}$. On the other hand, for any two dual spheres of $\S_{\O^6}$, we have that $\langle S_f,S_{f'} \rangle\in\{-2,-\frac32,-1,-\frac12,0,\frac12,1\}$,
so  $\O^{6}$ is crystallographic again by Lemma \ref{lem:crystal}.\\

We now discuss the integrality part.  The regular polytopes $\I^3,\D^3,\D^4$  are not integral in $\mathbb Z$ by Lemma \ref{lem:integral}, but in $\mathbb Z[\varphi]$ (see \cite{RR21_1}  and section \ref{sec:533} at the Appendix). The integrality of  $\T^3$ \cite{soddy1936}, $\O^3$ \cite{guettler}, $\C^3$ \cite{stange2015bianchi},  $\T^4$ \cite{soddy1936},  $\O^4$ \cite{nakamura2014localglobal,Dias2014TheLP} is well-known, while the cases of $\C^4$ and $\R^4$ are described in the sections \ref{sec:433} and \ref{sec:343}. We end by showing that $\O^6$ is integral. Let 
$\S_{\O^6}\subset\widehat{\mathbb R^5}$ be a strip $5$-sphere packing of $\O^6$. The bends of $\S_{\O^6}$ are integers. Let $(S_{v_1},\ldots,S_{v_7})\subset\S_{\O^6} $ corresponding to a fundamental basis with respect to a flag  $\Phi=(v,\ldots,f,\O^6)$. The corresponding bend matrices of the fundamental symmetries $r_1,\ldots,r_6$ with respect to $\Phi$  are
\begin{align*}
&	\mathbf{R}_1={\footnotesize\left(
		\begin{array}{ccccccc}
			& 1 &  &  &  &  &  \\
			1 &  &  &  &  &  &  \\
			&  & 1 &  &  &  &  \\
			&  &  & 1 &  &  &  \\
			&  &  &  & 1 &  &  \\
			&  &  &  &  & 1 &  \\
			1 & -1 &  &  &  &  & 1 \\
		\end{array}
		\right)}\quad\mathbf{R}_2={\footnotesize\left(
		\begin{array}{ccccccc}
			1 &  &  &  &  &  &  \\
			&  & 1 &  &  &  &  \\
			& 1 &  &  &  &  &  \\
			&  &  & 1 &  &  &  \\
			&  &  &  & 1 &  &  \\
			&  &  &  &  & 1 &  \\
			&  &  &  &  &  & 1 \\
		\end{array}
		\right)}\quad\mathbf{R}_3={\footnotesize\left(
		\begin{array}{ccccccc}
			1 &  &  &  &  &  &  \\
			& 1 &  &  &  &  &  \\
			&  &  & 1 &  &  &  \\
			&  & 1 &  &  &  &  \\
			&  &  &  & 1 &  &  \\
			&  &  &  &  & 1 &  \\
			&  &  &  &  &  & 1 \\
		\end{array}
		\right)}\quad\\
	&\mathbf{R}_4={\footnotesize\left(
		\begin{array}{ccccccc}
			1 &  &  &  &  &  &  \\
			& 1 &  &  &  &  &  \\
			&  & 1 &  &  &  &  \\
			&  &  &  & 1 &  &  \\
			&  &  & 1 &  &  &  \\
			&  &  &  &  & 1 &  \\
			&  &  &  &  &  & 1 \\
		\end{array}
		\right)}\quad\mathbf{R}_5={\footnotesize\left(
		\begin{array}{ccccccc}
			1 &  &  &  &  &  &  \\
			& 1 &  &  &  &  &  \\
			&  & 1 &  &  &  &  \\
			&  &  & 1 &  &  &  \\
			&  &  &  &  & 1 &  \\
			&  &  &  & 1 &  &  \\
			&  &  &  &  &  & 1 \\
		\end{array}
		\right)}\quad\mathbf{R}_6={\footnotesize\left(
		\begin{array}{ccccccc}
			1 &  &  &  &  &  &  \\
			& 1 &  &  &  &  &  \\
			&  & 1 &  &  &  &  \\
			&  &  & 1 &  &  &  \\
			&  &  &  & 1 &  &  \\
			1 &  &  &  &  & -1 & 1 \\
			&  &  &  &  &  & 1 \\
		\end{array}
		\right)}
	\end{align*}
	while the bend matrix of the inversion through the dual sphere $S_f$ is
	\begin{align*}
	\mathbf{S}_f={\footnotesize\left(
		\begin{array}{ccccccc}
			1&  &  &  &  &  &  \\
			& 1 &  &  &  &  &  \\
			&  & 1 &  &  &  &  \\
			&  &  & 1 &  &  &  \\
			&  &  &  & 1 &  &  \\
			&  &  &  &  & 1 &  \\
			-1 & 1 & 1 & 1 & 1 & 1 & -1 \\
		\end{array}
		\right)}
\end{align*}
Therefore, the full symmetry group $\Gamma_{\{3,3,3,3,4\}}=\langle\mathbf{R}_1,\ldots,\mathbf{R}_6,\mathbf{S}_f\rangle< \mathrm{SL}_{7}(\mathbb Z)$ so $\mathscr P(\S_{\O^6})$ is integral.
\end{proof}	

We denote the packings arising from one of the polytopes described in Theorem \ref{thm:classification} as the \textit{regular crystallographic packings}. These sphere packings also fall into the category of \textit{Boyd-Maxwell packings} introduced by Boyd in \cite{boyd}, which can be considered as a particular case of crystallographic packings where the co-cluster and the cluster are induced by a Coxeter simplex and its dual basis, respectively. Consequently, a Boyd-Maxwell packing is entirely determined by the hyperbolic Coxeter group generated by the reflections through the co-cluster, coinciding with the full symmetry group in the case of a regular crystallographic packing. In Table \ref{tab:coxgraphs}, we present the Coxeter graphs of the full symmetry group for the 11 regular crystallographic packings. The relations between the generators can be verified through direct computations on the matrices of the linear representations described in the Appendix \ref{sec:appendix}. The six graphs for $d=4,6$ can be found in various enumerations \cite{maxwell,chenlabbe,chen2016even}.

\begin{table}[H]
	\footnotesize
	\centering
	\begin{tabular}{cc}
		$\P$&$\Gamma_\P$\\
		\toprule
		$\T^3$ &	\begin{tikzpicture}[scale=.8]
			\draw[thick] (0,0) -- (3,0) node[pos=.833,above,inner sep=1pt] {$\infty$};
			
			\foreach \x in {0,1,2,3}
			{
				\node at (\x,0) [circle,fill=black,inner sep=0pt,minimum size=.15cm]  {};
			}
		\end{tikzpicture} \\
		$\O^3$ &	\begin{tikzpicture}[scale=.8]
			\draw[thick] (0,0) -- (3,0) 
			node[pos=.833,above,inner sep=1pt] {$\infty$}
			node[pos=.5,above,inner sep=1pt] {$4$}
			;
			
			\foreach \x in {0,1,2,3}
			{
				\node at (\x,0) [circle,fill=black,inner sep=0pt,minimum size=.15cm]  {};
			}
		\end{tikzpicture} \\
		$\C^3$ &	\begin{tikzpicture}[scale=.8]
			\draw[thick] (0,0) -- (3,0) 
			node[pos=.166,above,inner sep=1pt] {$4$}
			node[pos=.833,above,inner sep=1pt] {$\infty$};
			
			\foreach \x in {0,1,2,3}
			{
				\node at (\x,0) [circle,fill=black,inner sep=0pt,minimum size=.15cm]  {};
			}
		\end{tikzpicture}  \\
		$\I^3$ &	\begin{tikzpicture}[scale=.8]
			\draw[thick] (0,0) -- (3,0) 
			node[pos=.5,above,inner sep=1pt] {$5$}
			node[pos=.833,above,inner sep=1pt] {$\infty$};
			
			\foreach \x in {0,1,2,3}
			{
				\node at (\x,0) [circle,fill=black,inner sep=0pt,minimum size=.15cm]  {};
			}
		\end{tikzpicture} \\
		$\D^3$ &	\begin{tikzpicture}[scale=.8]
			\draw[thick] (0,0) -- (3,0) 
			node[pos=.166,above,inner sep=1pt] {$5$}
			node[pos=.833,above,inner sep=1pt] {$\infty$};
			
			\foreach \x in {0,1,2,3}
			{
				\node at (\x,0) [circle,fill=black,inner sep=0pt,minimum size=.15cm]  {};
			}
		\end{tikzpicture} \\
		
		$\T^4$ &	\begin{tikzpicture}[scale=.8]
			\draw[thick] (0,0) -- (4,0) 			
			node[pos=.125,above,inner sep=1pt] {}
			node[pos=.375,above,inner sep=1pt] {} 
			node[pos=.625,above,inner sep=1pt] {} 
			node[pos=.875,above,inner sep=1pt] {$6$};

			\foreach \x in {0,1,2,3,4}
			{
				\node at (\x,0) [circle,fill=black,inner sep=0pt,minimum size=.15cm]  {};
			}
		\end{tikzpicture}  \\
		
		$\O^4$
		& 	\begin{tikzpicture}[scale=.8]
			\draw[thick] (0,0) -- (4,0) 
			
			node[pos=.125,above,inner sep=1pt] {}
			node[pos=.375,above,inner sep=1pt] {} 
			node[pos=.625,above,inner sep=1pt] {$4$} 
			node[pos=.875,above,inner sep=1pt] {$4$};
			
			\foreach \x in {0,1,2,3,4}
			{
				\node at (\x,0) [circle,fill=black,inner sep=0pt,minimum size=.15cm]  {};
			}
		\end{tikzpicture} \\
		$\C^4$
		& 	\begin{tikzpicture}[scale=.8]
			\draw[thick] (0,0) -- (4,0) 
			node[pos=.125,above,inner sep=1pt] {$4$}
			node[pos=.375,above,inner sep=1pt] {} 
			node[pos=.625,above,inner sep=1pt] {} 
			node[pos=.875,above,inner sep=1pt] {$6$};

			\foreach \x in {0,1,2,3,4}
			{
				\node at (\x,0) [circle,fill=black,inner sep=0pt,minimum size=.15cm]  {};
			}
		\end{tikzpicture} \\
		
		$\R^4$
		& 	\begin{tikzpicture}[scale=.8]
			\draw[thick] (0,0) -- (4,0) 
			
			node[pos=.125,above,inner sep=1pt] {}
			node[pos=.375,above,inner sep=1pt] {$4$} 
			node[pos=.625,above,inner sep=1pt] {} 
			node[pos=.875,above,inner sep=1pt] {$6$};
			
			\foreach \x in {0,1,2,3,4}
			{
				\node at (\x,0) [circle,fill=black,inner sep=0pt,minimum size=.15cm]  {};
			}
		\end{tikzpicture}\\
		$\D^4$
		& 	\begin{tikzpicture}[scale=.8]
			\draw[thick] (0,0) -- (4,0) 
			node[pos=.125,above,inner sep=1pt] {$5$}
			node[pos=.375,above,inner sep=1pt] {} 
			node[pos=.625,above,inner sep=1pt] {} 
			node[pos=.875,above,inner sep=1pt] {$6$};

			\foreach \x in {0,1,2,3,4}
			{
				\node at (\x,0) [circle,fill=black,inner sep=0pt,minimum size=.15cm]  {};
			}
		\end{tikzpicture}\\
		$\O^6$
		& 	\begin{tikzpicture}[scale=.8]
			\draw[thick] (0,0) -- (6,0) 
			
			node[pos=.125,above,inner sep=1pt] {}
			node[pos=.375,above,inner sep=1pt] {} 
			node[pos=.75,above,inner sep=1pt] {$4$} 
			node[pos=.875,above,inner sep=1pt] {};
			
			\foreach \x in {0,1,2,3,4,5,6}
			{
				\node at (\x,0) [circle,fill=black,inner sep=0pt,minimum size=.15cm]  {};
			}
		\end{tikzpicture} 
	\end{tabular}
	\caption{The Coxeter graphs of the full symmetry groups of the 11 regular crystallographic packings.}
	\label{tab:coxgraphs}
\end{table}

\subsection{A Descartes theorem for the regular $4$-polytopes}\label{sec:platonic}

In \cite{RR21_1}, the authors presented the following Descartes quadratic equation satisfied by any fundamental bend vector $\mathbf{b}=(b_1,b_2,b_3,b_4)^T$ of a Platonic circle packing $\S_\P$ in terms of the Schläfli symbol of $\P$
\begin{align}\label{eq:descartes3d}Q_{\{p,q\}}(b_1,b_2,b_3,b_4)=\mathbf{b}^T\mathbf{Q}_{\{p,q\}} \mathbf{b}=0
\end{align}
where $\mathbf Q_{\{p,q\}}$ is the bisymmetric matrix
\begin{align}	\mathbf Q_{\{p,q\}}=
	\left(
	\begin{array}{ccccc}
		1 & a & b & -1  \\
		a & d & c & b \\
		b & c & d & a \\
		-1 & b& a & 1\\
	\end{array}
	\right)&&\text{ with }&&\begin{array}{rl}
		a&=-1-\omega _p-\omega _q\\
		b&=-1+\omega _p-\omega _q\\
		c&=-1-\omega _p^2+\omega _q^2\\		
		d&=\left(1+\omega _p+\omega _q\right)^2 \\
	\end{array}
\end{align}
and $\omega_{n}:=1+2\cos\frac {2\pi}n$. Notice that $\omega_3=0$, $\omega_4=1$, $\omega_5=\varphi$ and $Q_{\{3,3\}}$ is the classic Descartes quadratic form. The later was used to compute linear representations of the full symmetry groups and the integrality conditions of the Platonic solids. The equation \eqref{eq:descartes3d} was obtained from the following more general formula, 	proved by Ram\'irez Alfons\'in and the author in \cite{RR21_1}, defined for every polytopal sphere packing induced by a uniform polytope (which includes the regular family).

\begin{thm}[\cite{RR21_1}] \label{thm:poldescartes}
	Let $\mathcal S_\P$ a polytopal sphere packing where $\P$ is a uniform $(d+1)$-polytope with $d\ge2$. For any flag $(f_0,\ldots,f_d,f_{d+1}=\P)$ we have
	\begin{align}\label{eq:poldesth}
		(\kappa_{f_0}-\kappa_{f_1})^2+\ell_{f_2}^2(\kappa_{f_1}-\kappa_{f_2})^2+\sum_{i=2}^{d}\frac{1}{\ell_{f_{i+1}}^{-2}-\ell_{f_i}^{-2}}(\kappa_{f_i}-\kappa_{f_{i+1}})^2=\ell_\P^2\kappa_{\P}^2
	\end{align}
	where $	\kappa_{f_i}$ is the arithmetic mean of the bends of the spheres corresponding to all the vertices of $f_i$, and $\ell_i$ is the half edge-length of a canonical realization of $f_i$.
\end{thm}
The variables $	\kappa_{f_i}$ and $\ell_{f_i}$ are called the \textit{polytopal curvature} and the \textit{canonical length}. By proceeding similarly in the next dimension, we obtain the following unified Descartes Theorem for the regular $4$-polytopes.

\begin{thm} \label{thm:4regdesth} Let $\mathcal S_{\P}$ be a polytopal sphere packing where $\P$ is the regular $4$-polytope with Schläfli symbol $\{p,q,r\}$. For any fundamental bend vector $\mathbf b=(b_1,b_2,b_3,b_4,b_5)^T$ of $\mathcal S_{\P}$ we have
	\begin{align}\label{eq:platonicdescartes}
		Q_{\{p,q,r\}}(b_1,b_2,b_3,b_4,b_5)=\mathbf{b}^T\mathbf{Q}_{\{p,q,r\}} \mathbf{b}=0
	\end{align}	 
	where $\mathbf Q_{\{p,q,r\}}$ is the bisymmetric matrix
	\begin{align}
		\mathbf Q_{\{p,q,r\}}=
		\left(
		\begin{array}{ccccc}
			2 & a & b & c & -1+ \omega _r\\\
			a & d & e & h & c \\
			b & e & i & e & b \\
			c & h & e & d& a \\
			-1+ \omega _r & c & b & a & 2 \\
		\end{array}
		\right)
	\end{align}
	with
	\begin{align*}
		a&=-1-2 \omega _p-2 \omega _q-\omega _r\\
		b&=-1+\omega _p-\omega _q-\omega _r+\omega _p \omega _r-\omega _q \omega _r\\
		c&=-1+\omega _p+\omega _q-\omega _r-\omega _p \omega _r-\omega _q \omega _r\\
		d&=2+2 \left(1+\omega _p+\omega _q\right) \left(\omega _p+\omega _q+\omega _r\right)\\
		e&=-1-\omega _p-\omega _p^2+\omega _q^2+\omega _r^2-\omega _p \omega _r+\omega _q \omega _r-\omega _p^2 \omega _r+\omega _q^2 \omega _r+\omega _q \omega _r^2\\
		h&=-1+2 \omega _p-\omega _q-\omega _p^2-\omega _q^2+\omega _r^2-2 \omega _p \omega _q+2 \omega _p \omega _r+\omega _p^2 \omega _r+2 \omega _p \omega _q \omega _r+\omega _q^2 \omega _r+\omega _q
		\omega _r^2\\
		i&=2 \left(1+\omega _p+2 \omega _q+\omega _r+\omega _p^2+\omega _q^2+\omega _p \omega _q+2 \omega _q \omega _r+\omega _p^2 \omega _r+\omega _q^2 \omega _r-\omega _p \omega _r^2-\omega _p \omega _q
		\omega _r^2\right)
	\end{align*}
\end{thm}

\begin{proof}
	Let $\Phi=(f_0,f_1,f_2,f_3,f_4=\P)$ be the flag of $\P$ corresponding to a fundamental bend vector $\mathbf b=(b_1,b_2,b_3,b_4,b_5)$. The canonical lengths of $f_2$, $f_3$ and $f_4$ can be computed in terms of $\omega_p$, $\omega_q$ and $\omega_r$ by
		\begin{align}\label{eq:lvaluesfp}
		\ell_{f_2}=\sqrt{\frac{3-\omega_p}{1+\omega_p}},&&\ell_{f_3}=\sqrt{\frac{2-\omega_p-\omega_q}{1+\omega_p}},&&\ell_{f_4}=\sqrt{\frac{5-3(\omega_p+\omega_q+\omega_r)-\omega_q+\omega_p\omega_r}{(1+\omega_p)(3-\omega_r)}}.
	\end{align}
	
By replacing these values in the quadratic equation given by the Polytopal Descartes' Theorem and then solving and adding both solutions for each polytopal curvature $\kappa_{f_i}$, we obtain the following relations
	\begin{align}
		\label{eq:cycvertex}	\kappa_{f_0}=&b_1\\
		\label{eq:cycedge}	\kappa_{f_1}=&\frac{b_1+b_2}2 \\
		\label{eq:cycface}	\kappa_{f_2}=&\frac{b_1+b_2+b_3-\omega_p b_2}{3-\omega_p}\\
		\label{eq:cycfacet}	\kappa_{f_3}=&\frac{b_1+b_2+b_3+b_4-(\omega_p+\omega_q)(b_2+b_3)}{4-2(\omega_p+\omega_q)}\\
		\label{eq:cycpol}	\kappa_{f_4}=&\frac{b_1+b_2+b_3+b_4+b_5-(\omega_p+\omega_q+\omega_r)(b_2+b_3+b_4)-(\omega_q-\omega_p\omega_r)b_3}{5-3(\omega_p+\omega_q+\omega_r)-\omega_q+\omega_p\omega_r}
	\end{align}
	The above relations define a transition matrix $\mathbf T$ satisfying
	\begin{align}\label{eq:transition}
		\mathbf k=\mathbf T \mathbf b
	\end{align} 
	where $\mathbf k=(\kappa_{f_0},\kappa_{f_1},\kappa_{f_2},\kappa_{f_3},\kappa_{f_4})^T$. Let $\mathbf Q_\Phi$ the matrix of the quadratic form induced by Polytopal Descartes' Theorem for $\P$. Then, the quadratic equation \eqref{eq:poldesth} becomes
	\begin{align}
		\mathbf k^T\mathbf Q_\Phi\mathbf k=0\Leftrightarrow \mathbf b^T\mathbf T^T\mathbf Q_\Phi\mathbf T\mathbf b=0.
	\end{align} 
	It can be checked by direct computations that $\mathbf Q_{\{p,q,r\}}=2(1+\omega_p)(1+\omega_q)(1+\omega_r)(3-\omega_r)\mathbf T^T\mathbf Q_\Phi\mathbf T$.
\end{proof}

As it is done in \cite{RR21_1}, the transition matrix defined above can be used to compute the bend matrices of the fundamental symmetries $r_1,r_2,r_3,r_4$, and $s_f$ which generate the full symmetry group $\Gamma_{\{p,q,r\}}$. These bend matrices belong to $\mathrm{SL}_5(\mathbb{Z}[\omega_p,\omega_r])$ and are orthogonal with respect to $\mathbf Q_{\{p,q,r\}}$. We notice that for each of the six regular $4$-polytopes with Schläfli symbol $\{p,q,r\}$, we have that $\mathbb Z[\omega_p,\omega_q,\omega_r]=\mathbb Z[\omega_p,\omega_r]$ since $q=3,4$ so $\omega_q=0,1$.  

\begin{cor}[Linear representation of the full symmetry group] \label{cor:lineargroup} For each regular $4$-polytope $\{p,q,r\}$, $\Gamma_{\{p,q,r\}}$ is a discrete orthogonal subgroup of $\mathrm{SL}_5(\mathbb{Z}[\omega_p,\omega_r])$ with respect to $\mathbf Q_{\{p,q,r\}}$.	
\end{cor}
We denote by $\mathscr P_{\{p,q,r\}}(b_1,b_2,b_3,b_4)$ the Apollonian arrangement of the regular $4$-polytope with Schläfli symbol $\{p,q,r\}$, where the quadruple $(b_1,b_2,b_3,b_4)$ is made of four consecutive entries of a fundamental bend vector. Since $\mathscr P_{\{p,q,r\}}(b_1,b_2,b_3,b_4)$ is unique up to Euclidean isometries, the set of bends is fully determined by the quadruple. By solving the equation \eqref{eq:platonicdescartes} for $b_5$, we obtain the following.
	\begin{cor}[Integrality condition]\label{cor:integrality} Let $b_1,b_2,b_3,b_4$ be four consecutive entries of a fundamental bend vector of a polytopal sphere packing for a regular $4$-polytope $\{p,q,r\}$. If $b_1,b_2,b_3,b_4,\sqrt{\Delta_{\{p,q,r\}}}$ are in $\mathbb Z[\omega_p,\omega_r]$ where
			\begin{align}\label{eq:integralitycond}\Delta_{\{p,q,r\}}:=(\omega_r+1)(\omega_r-3)Q_{\{p,q\}}(b_1,b_2,b_3,b_4),
		\end{align}
		then $\mathscr P_{\{p,q,r\}}(b_1,b_2,b_3,b_4)$ is $\mathbb Z[\omega_p,\omega_r]$-integral.
	\end{cor}	
The reader can find linear representations of the full symmetry groups, the integrality conditions and examples of integral packings for the six regular 4-polytopes in the Appendix \ref{sec:appendix}

\subsection{Apollonian sections} \label{sec:aposections}
The main goal of this section is to study which Platonic crystallographic packings appear as cross sections of the 3-dimensional regular crystallographic packings. We define an \textit{Apollonian section} of an arrangement of spheres (not necessarily a packing) $\mathscr{P}:= \Gamma\cdot \S$ as a subset $\mathscr{S}:= G\cdot X$ where $G<\Gamma$ and $X\subset \S$. We say that $\mathscr S$ is \textit{geometric} if there is a sphere $\Sigma$, that we call the \textit{cutting sphere}, that is invariant under the action of $G$ and intersecting the interior of all the spheres of $\mathscr S$. Two Apollonian sections $\mathscr S=G\cdot X\subset\mathscr P $ and $\mathscr S'=G'\cdot X'\subset\mathscr P'$ are \textit{algebraically equivalent}, denoted by $\mathscr S\simeq\mathscr S'$,  if $G$ and $G'$ are isomorphic and there is equivariant bijection between $\mathscr S$ and $\mathscr S'$ with respect to the actions. If the bijection preserves the bends, we shall denote this equivalence by $\mathscr S\equiv\mathscr S'$. Let $\{p,q\}$ and $\{r,s,t\}$ be the Schläfli symbols of a regular $3$- and $4$-polytope respectively. We say that an Apollonian section $\mathscr S_{\{r,s,t\}}^{\{p,q\}}$ is \textit{Platonic} if it is geometric and it satisfies $\mathscr P_{\{p,q\}}\simeq\mathscr S_{\{r,s,t\}}^{\{p,q\}}\subset \mathscr P_{\{r,s,t\}}$. We notice that $\mathscr S_{\{r,s,t\}}^{\{p,q\}}$ is unique up to symmetries of $\mathscr P_{\{r,s,t\}}$. The notion of Apollonian section will allow us to prove the following.

\begin{thm}\label{thm:aposecs}
	There are the following relations between the Apollonian arrangements of the regular 3- and 4-polytopes
	\begin{equation}\label{eq:aposeclist}
		\begin{split}
			\mathscr P_{\{3,3\}}\prec\mathscr{P}_{\{3,3,3\}},\\
			\mathscr P_{\{3,3\}}, \mathscr P_{\{3,4\}}, \mathscr P_{\{4,3\}}\prec\mathscr{P}_{\{3,3,4\}},\\
			\mathscr P_{\{4,3\}}\prec\mathscr{P}_{\{4,3,3\}},\\
			\mathscr P_{\{3,4\}}, \mathscr P_{\{4,3\}}\prec\mathscr{P}_{\{3,4,3\}},\\
			\mathscr P_{\{3,3\}}, \mathscr P_{\{3,5\}}\prec\mathscr{P}_{\{3,3,5\}},\\
			\mathscr P_{\{5,3\}}\prec\mathscr{P}_{\{5,3,3\}},
		\end{split}	
	\end{equation}
	where $	\mathscr P_{\{p,q\}}\prec\mathscr{P}_{\{r,s,t\}}$ means that  $	\mathscr P_{\{p,q\}}$ can be obtained as a cross-section of  $\mathscr{P}_{\{r,s,t\}}$.
\end{thm}

\begin{proof}
For each relation $	\mathscr P_{\{p,q\}}\prec\mathscr{P}_{\{r,s,t\}}$ described in \eqref{eq:aposeclist}, we shall construct a Platonic Apollonian section $\mathscr P_{\{p,q\}}\simeq\mathscr S_{\{r,s,t\}}^{\{p,q\}}\subset \mathscr P_{\{r,s,t\}}$ by defining a homomorphism between the full symmetry groups $\phi_{\{r,s,t\}}^{\{p,q\}}:\Gamma_{\{p,q\}}\rightarrow \Gamma_{\{r,s,t\}}$ which can be easily deduced by comparing the action of the fundamental symmetries in the strip packings. All the homomorphisms and the cutting spheres $\Sigma_{\{r,s,t\}}^{\{p,q\}}$ are described in the Appendix \ref{sec:appendix}. These homomorphisms give us the following algebraic equivalences
	\begin{align}
	\mathscr P_{\{p,q\}}:=\Gamma_{\{p,q\}}\cdot\{C_v\}\simeq G_{\{r,s,t\}}^{\{p,q\}}\cdot\{S_v\}\subset \mathscr P_{\{r,s,t\}}
	\end{align}
where $G_{\{r,s,t\}}^{\{p,q\}}:=\phi_{\{r,s,t\}}^{\{p,q\}}(\Gamma_{\{p,q\}})<\Gamma_{\{r,s,t\}}$, and $C_v$ and $S_v$ are the initial circle and the initial sphere, respectively, in a strip packing of $\{p,q\}$ and $\{r,s,t\}$. Thus, we can set $\mathscr S^{\{p,q\}}_{\{r,s,t\}}:=G_{\{r,s,t\}}^{\{p,q\}}\cdot\{S_v\}$.
\end{proof}

\begin{conj}
	The list given in \eqref{eq:aposeclist} is complete.
\end{conj}

The bends of the circle packing obtained by the intersection of $\mathscr P_{\{r,s,t\}}$ with a cutting sphere $\Sigma$ depend on the intersection angle and the bend of $\Sigma$. Indeed, if $S$ and $\Sigma$ are two intersecting $d$-spheres for any $d\ge2$ with intersection angle $0<\alpha<\pi$, then, by combining the Law of cosines with Heron's Formula, we obtain that the bend of the  $(d-1)$-sphere $S\cap\Sigma$ is given by
\begin{align}\label{eq:intersection}
	b_{S\cap\Sigma}=\frac {\sqrt{b_S^2+b_\Sigma^2-2\cos(\alpha) b_S b_\Sigma}}{\sin(\alpha)}
\end{align}

\begin{figure}[H]
	\centering
	
	\begin{tikzpicture}
		
		\begin{scope}[xshift=-5.5cm,yshift=-.7cm]
			\node at (0,0) {\includegraphics[align=c,width=.34\textwidth]{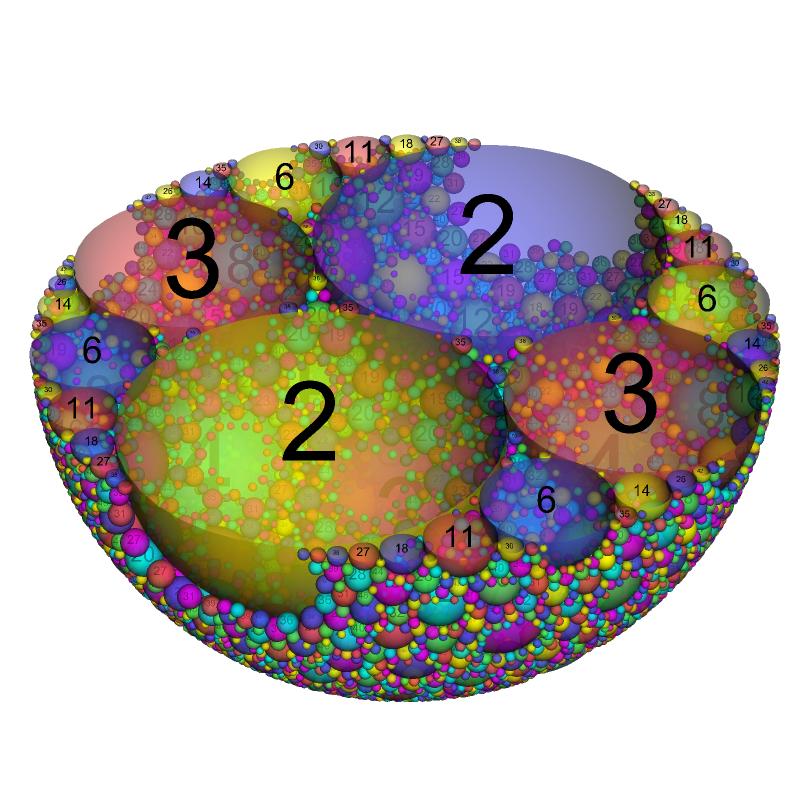}};
		\end{scope}
		\begin{scope}[xshift=-.1cm,yshift=-.1cm]
			\node at (0,0) {\includegraphics[align=c,width=.34\textwidth,clip,trim= 20 0 0 0]{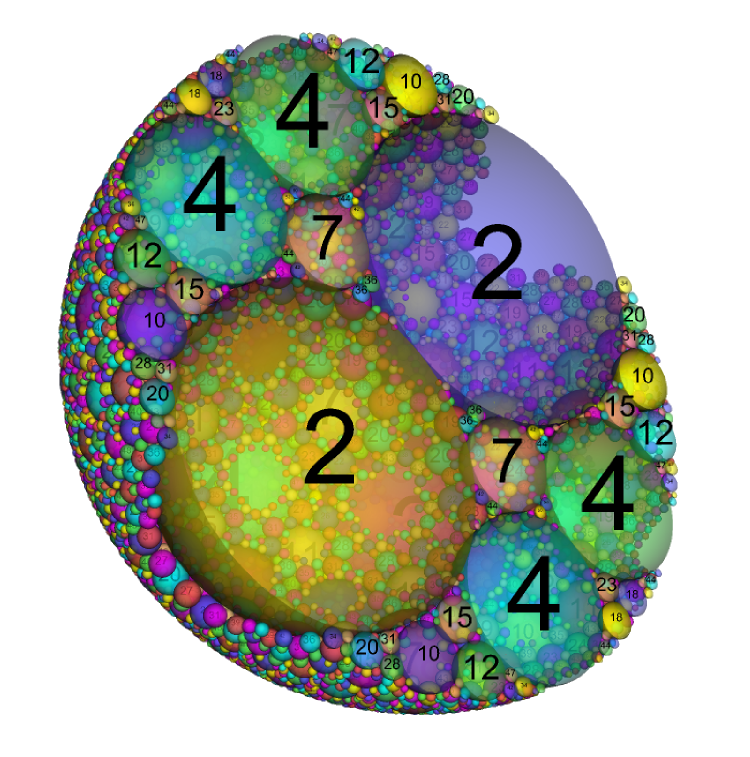}};
		\end{scope}
		\begin{scope}[xshift=5.5cm]
			\node at (0,0) {\includegraphics[align=c,width=.34\textwidth]{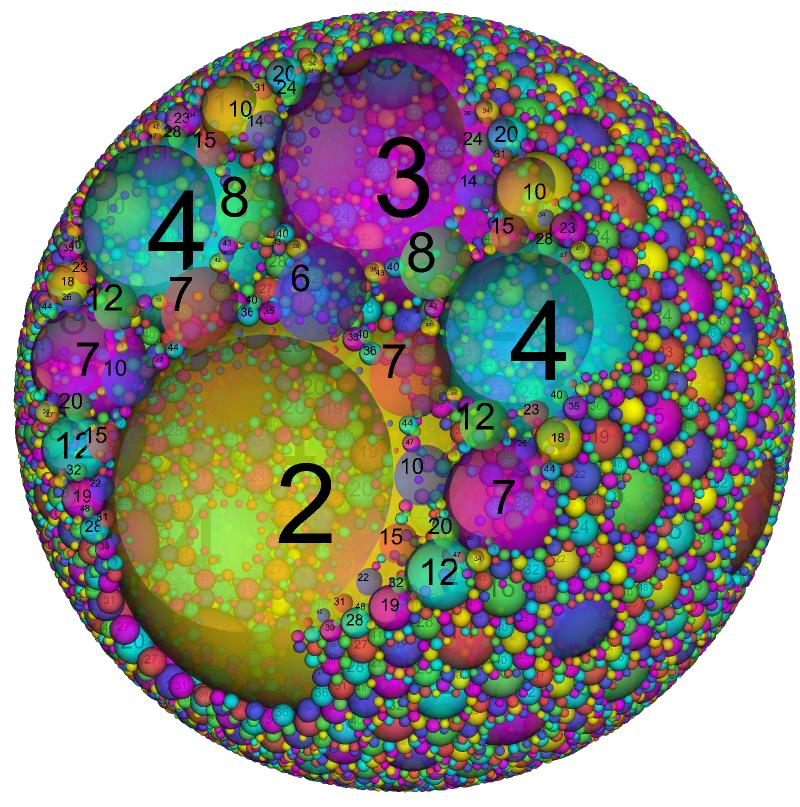}};
		\end{scope}
		
		\begin{scope}[yshift=-5.7cm]]
			
			\begin{scope}[xshift=-5.5cm]
				\node at (0,0) {\includegraphics[align=c,width=.34\textwidth]{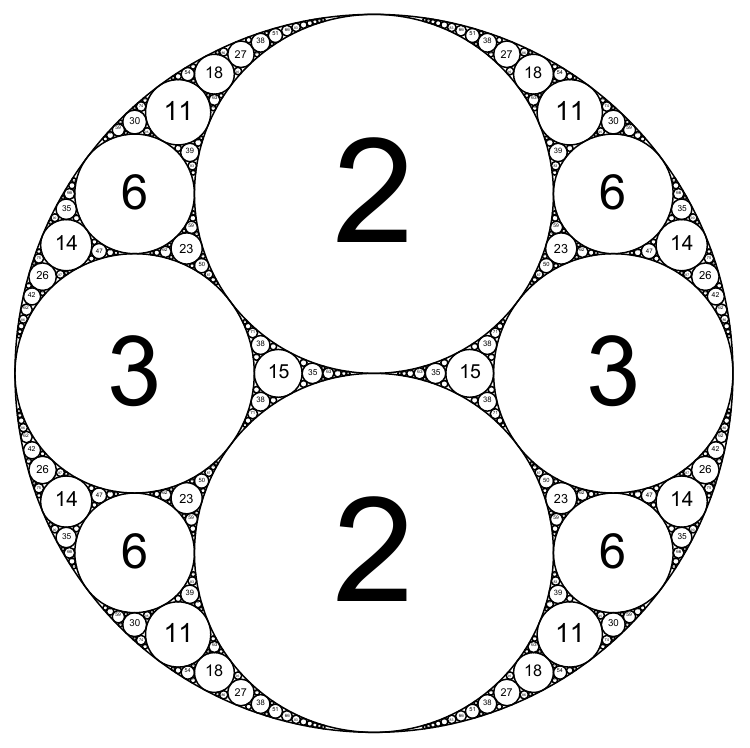}};
			\end{scope}
			\begin{scope}[xshift=0cm]
				\node at (0,0) {\includegraphics[align=c,width=.34\textwidth]{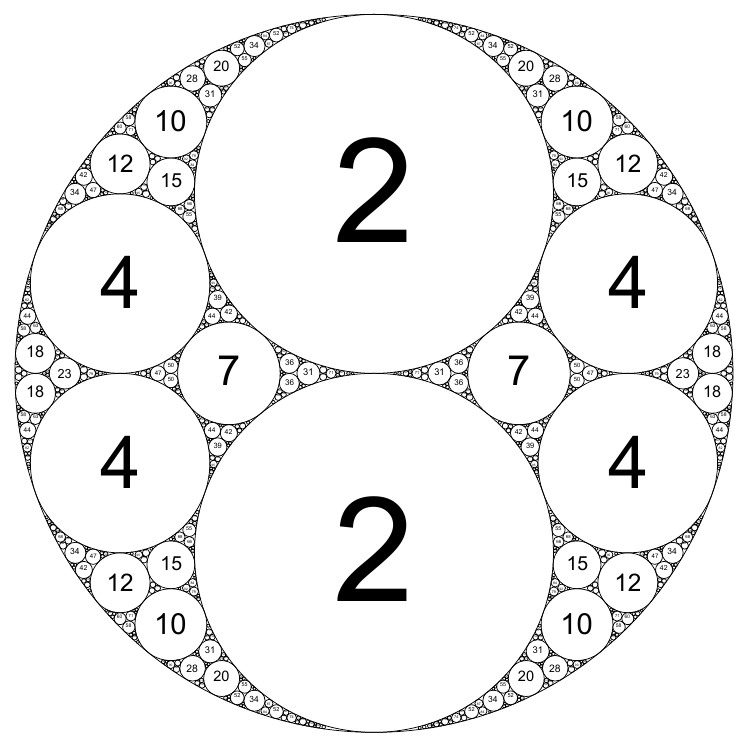}};
			\end{scope}
			\begin{scope}[xshift=5.5cm]
				\node at (0,0) {\includegraphics[align=c,width=.34\textwidth]{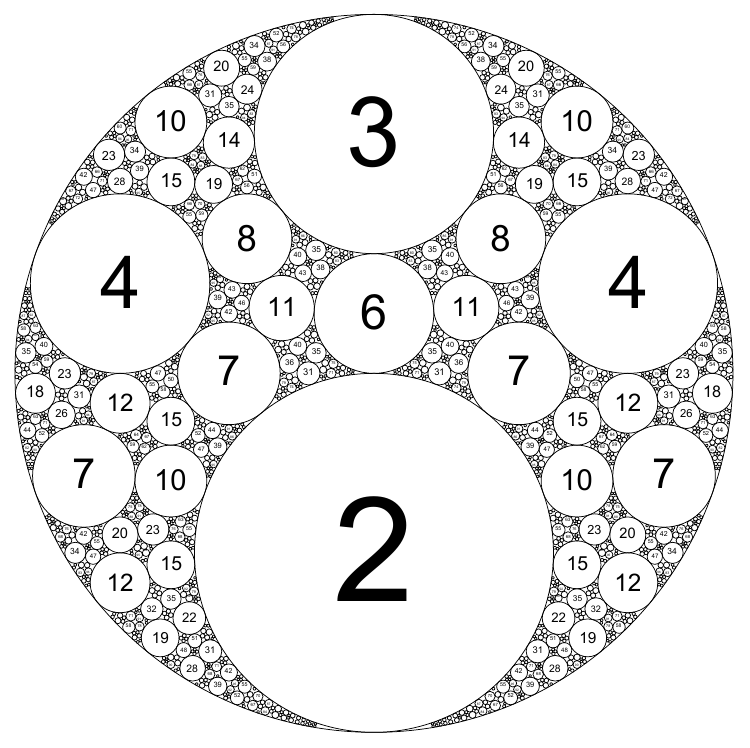}};;
			\end{scope}
			
		\end{scope}
		
	\end{tikzpicture}

	\caption{
		Top figures: The integral orthoplicial crystallographic packing $\mathscr P_{\{3,3,4\}}(-1,2,2,3)$ truncated with the cutting planes of an integral tetrahedral Apollonian section (left), an integral octahedral Apollonian section (center), and an integral cubic Apollonian section (right).
		Bottom figures: the corresponding integral crystallographic Platonic packings including the tetrahedral $\mathscr P_{\{3,3\}}(-1,2,2)$ (left), the octahedral $\mathscr P_{\{3,4\}}(-1,2,2)$, and the cubic $\mathscr P_{\{4,3\}}(-1,2,6)$.
	}
	\label{fig:arithemeticsections}
\end{figure}

In the context of the Platonic Apollonian sections of Theorem \ref{thm:aposecs}, the cutting sphere $\Sigma_{\{r,s,t\}}^{\{p,q\}}$ intersects all the spheres in $\mathscr S_{\{r,s,t\}}^{\{p,q\}}$ with the same angle. This follows from the invariance of $\Sigma_{\{r,s,t\}}^{\{p,q\}}$ under the group that generates $\mathscr S_{\{r,s,t\}}^{\{p,q\}}$. We remark that the sections mentioned above are all orthogonal, except for $\mathscr S_{\{3,3,4\}}^{\{4,3\}}$, where the intersection angle is $\pi/4$ (see section \ref{sec:334} in the Appendix). In particular, when the cutting sphere $\Sigma_{\{r,s,t\}}^{\{p,q\}}$ is a plane, then, by equation \eqref{eq:intersection}, we have that the bends of the spheres in $\mathscr S_{\{r,s,t\}}^{\{p,q\}}$ are equal to the bends of the circles in the Platonic packing $\mathscr P_{\{p,q\}}$ obtained by the intersection $\mathscr S_{\{r,s,t\}}^{\{p,q\}}\cap\Sigma_{\{r,s,t\}}^{\{p,q\}}$, except for $\mathscr S_{\{3,3,4\}}^{\{4,3\}}$, where the bends are rescaled by $1/\sqrt 2$. Therefore, when $\mathscr P_{\{r,s,t\}}$ is integral, the intersection with a cutting plane $\Sigma_{\{r,s,t\}}^{\{p,q\}}$ gives an integral packing $\mathscr P_{\{p,q\}}$ (up to rescaling by $1/\sqrt 2$ for $\mathscr S_{\{3,3,4\}}^{\{4,3\}}$). Figure \ref{fig:arithemeticsections} illustrates an integral tetrahedral, octahedral and cubic crystallographic packing obtained by the intersection of three different cutting planes with an integral orthoplicial crystallographic packing.

We shall demonstrate that every integral Platonic crystallographic packing can obtained in this manner. In \cite{RR21_1}, the integrality conditions for the three integral Platonic solids assert that if $b_1,b_2,b_3$ are three consecutive entries of a fundamental bend vector of a Platonic crystallographic packing $\mathscr{P}_{\{p,q\}}$ satisfying $b_1,b_2,b_3,\sqrt{\Delta_{\{p,q\}}}\in\mathbb Z$, where
 \begin{align}
	\Delta_{\{p,q\}}=(1+\omega_q)(b_1b_2+b_2b_3+b_3b_1-\omega_pb_2^2)
\end{align}
then $\mathscr{P}_{\{p,q\}}(b_1,b_2,b_3)$ is integral. The converse is also true. 
\begin{lem}\label{lem:converse}Let $\mathscr{P}_{\{p,q\}}(b_1,b_2,b_3)$ be an integral Platonic crystallographic packing. Then $\sqrt{\Delta_{\{p,q\}}}\in\mathbb Z$.
\end{lem}
\begin{proof}
	If $\mathscr{P}_{\{p,q\}}(b_1,b_2,b_3)$ is integral, then the fundamental bend vector $(b_1,b_2,b_3,b_4)$ is integral. By solving the Platonic Descartes equation \eqref{eq:descartes3d} for $b_4$, one obtains
	\begin{equation}
		b_4=b_1+(1-\omega_p+\omega_q)b_2+(1+\omega_p+\omega_q)b_3\pm2\sqrt{	\Delta_{\{p,q\}}}
	\end{equation}
	Therefore, in the integral case, $\{\omega_p,\omega_q\}=\{0,1\}$ so $m=2\sqrt{\Delta_{\{p,q\}}}$ is an integer with same parity as $m^2=4	(1+\omega_q)(b_1b_2+b_2b_3+b_3b_1-\omega_pb_2^2)$, which is an even integer. Hence, $m$ is even and $\sqrt{\Delta_{\{p,q\}}}$ is indeed an integer.
\end{proof}

\begin{thm}\label{cor:aposecsint}The following relations  holds:
		\begin{align}
			\mathscr P_{\{3,3\}}(b_1,b_2,b_3)&\equiv	\mathscr{S}_{\{3,3,3\}}^{\{3,3\}}\subset\mathscr P_{\{3,3,3\}}(b_1,b_2,b_3,b_1+b_2+b_3+2\sqrt{\Delta_{\{3,3\}}})\\
			\mathscr P_{\{3,3\}}(b_1,b_2,b_3)&\equiv	\mathscr{S}_{\{3,3,4\}}^{\{3,3\}}\subset\mathscr P_{\{3,3,4\}}(b_1,b_2,b_3,b_1+b_2+b_3+2\sqrt{\Delta_{\{3,3\}}})\\
			\mathscr P_{\{3,4\}}(b_1,b_2,b_3)&\equiv	\mathscr{S}_{\{3,4,3\}}^{\{3,4\}}\subset \mathscr P_{\{3,4,3\}}(b_1,b_2,b_3,b_1+2b_2+2b_3+2\sqrt{\Delta_{\{3,4\}}})\\
			\mathscr P_{\{4,3\}}(b_1,b_2,b_3)&\equiv	\mathscr{S}_{\{4,3,3\}}^{\{4,3\}}\subset \mathscr P_{\{4,3,3\}}(b_1,b_2,b_3,b_1+2b_3+\sqrt{\Delta_{\{4,3\}}})\\					
			\mathscr P_{\{3,4\}}(b_1,b_2,b_3)&\equiv	\mathscr{S}_{\{3,3,4\}}^{\{3,4\}}\subset \mathscr P_{\{3,3,4\}}(b_1,b_2,b_3,b_1+b_2+b_3+\sqrt{\Delta_{\{3,4\}}})\\				
			\mathscr P_{\{4,3\}}(b_1,b_2,b_3)&\equiv	\mathscr{S}_{\{3,3,4\}}^{\{4,3\}}\subset \mathscr P_{\{3,3,4\}}(b_1,b_2,b_3,b_1-b_2+b_3)
		\end{align}
		Moreover, for each relation described above, if the packing on the left-hand side is integral then the packing on the right-hand side is integral.
\end{thm}

\begin{proof}
	Let $\mathscr P_{\{3,3\}}(0,0,1)$ and $\mathscr P_{\{3,3,3\}}(0,0,1,1)$ be the crystallographic strip packing of the tetrahedron and the $4$-simplex, respectively.	The homomorphism $\phi_{\{3,3,3\}}^{\{3,3\}}:\Gamma_{\{3,3\}}\mapsto\Gamma_{\{3,3,3\}}$ described in the Appendix \eqref{eq:hom333} induces a bend-preserving equivariant bijection
	$$\mathscr P_{\{3,3\}}(0,0,1)\equiv\mathscr S_{\{3,3,3\}}^{\{3,3\}}\subset\mathscr P_{\{3,3,3\}}(0,0,1,1)$$
	which can be equally obtained by taking the intersection of the cutting sphere $\Sigma_{\{3,3,3\}}^{\{3,3\}}\cap\mathscr S_{\{3,3,3\}}^{\{3,3\}}\subset\mathscr P_{\{3,3,3\}}(0,0,1,1)$.	Let $\mathscr P_{\{3,3\}}(b_1,b_2,b_3)$ be any integral tetrahedral crystallographic packing. By solving \eqref{eq:descartes3d} on $b_4$, we have that  $\mathscr P_{\{3,3\}}(b_1,b_2,b_3)$ contains a tetrahedral circle packing with fundamental bend vector $(b_1, b_2, b_3, b_4=b_1+b_2+b_3+2\sqrt{\Delta_{\{3,3\}}})$. Due to the Möbius uniqueness of polytopal sphere packings, there is a Möbius transformation $\mu$ of $\widehat{\mathbb R^2}$ such that  $\mu:\mathscr P_{\{3,3\}}(0,0,1)\rightarrow\mathscr P_{\{3,3\}}(b_1,b_2,b_3)$. Up to reflection, there is a unique Möbius transformation $\widetilde\mu$ of $\widehat{\mathbb R^3}$ which acts as $\mu$ on the cutting sphere $\Sigma_{\{3,3,3\}}^{\{3,3\}}$. These two Möbius transformations induce the following mappings
	\begin{equation}
		\begin{tikzcd}
			\mathscr P_{\{3,3\}}(0,0,1) \arrow[d, "\mu"'] & \equiv& \mathscr S_{\{3,3,3\}}^{\{3,3\}}\arrow[d, "\widetilde\mu"'] & \subset & \mathscr P_{\{3,3,3\}}(0,0,1,1) \arrow[d, "\widetilde\mu"'] \\
			\mathscr P_{\{3,3\}}(b_1,b_2,b_3) & \equiv & \mathscr S_{\{3,3,3\}}^{\{3,3\}}& \subset & \mathscr P_{\{3,3,3\}}(b_1,b_2,b_3,b_4) 
		\end{tikzcd}
	\end{equation} 
	Indeed, since  $\Sigma_{\{3,3,3\}}^{\{3,3\}}$ is a plane cutting orthogonally the spheres of $\mathscr S_{\{3,3,3\}}^{\{3,3\}}$ and Möbius transformations preserve angles, then, by equation \eqref{eq:intersection}, we have that $b_1, b_2, b_3, b_4$ are also four consecutive entries of a fundamental bend vector of $\mathscr P_{\{3,3,3\}}(b_1,b_2,b_3,b_4)$. 	By Corollary \ref{cor:integrality}, we have that $\mathscr P_{\{3,3,3\}}(b_1,b_2,b_3,b_4)$ is also integral since $\sqrt{\Delta_{\{3,3,3\}}}=\sqrt{-3Q_{\{3,3\}}(b_1,b_2,b_3,b_4)}=0$, where $Q_{\{3,3\}}$ is the Descartes' quadratic form described in \eqref{eq:descartes3d}.\\ \medskip
	
	The same arguments apply to the relations involving the Platonic Apollonian sections of the form $\mathscr{S}_{\{p,q,r\}}^{\{p,q\}}$, specifically $\mathscr{S}_{\{3,3,4\}}^{\{3,3\}}$, $\mathscr{S}_{\{3,4,3\}}^{\{3,4\}}$, and $\mathscr{S}_{\{4,3,3\}}^{\{4,3\}}$. The strategy for handling the remaining two cases $\mathscr{S}_{\{3,3,4\}}^{\{3,4\}}$ and $\mathscr{S}_{\{3,3,4\}}^{\{4,3\}}$ is slightly different. The diffence arises because the forth bend $b_4$ of a fundamental bend vector of $\mathscr P_{\{3,3,4\}}(b_1,b_2,b_3,b_4)$ is not equal to the forth bend of a fundamental bend vector of $\mathscr P_{\{3,4\}}(b_1,b_2,b_3)$ or $\mathscr P_{\{4,3\}}(b_1,b_2,b_3)$.\\ \medskip
	
	 In the octahedral case, the value of $b_4$ corresponds to the polytopal curvature $\kappa_\P$ of an octahedral circle packing contained within $\mathscr P_{\{3,4\}}(b_1,b_2,b_3)$, which has a fundamental bend vector with three consecutive bends $b_1,b_2$, and $b_3$. The value of $\kappa_\P$ can be computed with equation \eqref{eq:poldesth}, which gives us 
	\begin{align}
		b_4=b_1+b_2+b_3+\sqrt{\Delta_{\{3,4\}}}=b_1+b_2+b_3+\sqrt{2(b_1b_2+b_2b_3+b_3b_1)}
	\end{align}
	Let us suppose that $\mathscr P_{\{3,4\}}(b_1,b_2,b_3)$ is integral. Checking the integrality of $\mathscr P_{\{3,3,4\}}(b_1,b_2,b_3,b_4)$ we obtain
	\begin{align*}
		\sqrt{\Delta_{\{3,3,4\}}}&=\sqrt{-Q_{\{3,3\}}(b_1,b_2,b_3,b_1+b_2+b_3+\sqrt{2(b_1b_2+b_2b_3+b_3b_1)})}\\
		&=\sqrt{2(b_1b_2+b_2b_3+b_3b_1)}\\
		&=	\sqrt{\Delta_{\{3,4\}}}\in\mathbb Z
	\end{align*}
	by Lemma \ref{lem:converse}. Hence, by Corollary \ref{eq:integralitycond}, $\mathscr{P}_{\{3,3,4\}}(b_1,b_2,b_3,b_4)$ is also integral.\\

	Finally, for the cubic Apollonian section $\mathscr{S}_{\{3,3,4\}}^{\{4,3\}}$, it is important to ocnsider the particularity that the cutting sphere does not intersect the spheres orthogonally. In this case, the inital bend-preserving equivariant bijection 
	is achieved by composing the intersection of $\mathscr{S}_{\{3,3,4\}}^{\{4,3\}}$ with the cutting sphere with a rescaling of factor $1/\sqrt 2$. This rescaling compensates for the intersecting angle factor, as detailed in \eqref{eq:intersection}. Now, four consecutive entries $b_1,b_2,b_3,b_4$ of a fundamental bend vector of $\mathscr P_{\{3,3,4\}}(b_1,b_2,b_3,b_4)$ correspond to the bends of four consecutive circles in a square face of  $\mathscr P_{\{4,3\}}(b_1,b_2,b_3)$. According to the work of Stange in \cite{stange2015bianchi}, these bends satisfy the following relationship
	\begin{align}
		b_1+b_3=b_2+b_4.
	\end{align}
		Then, if $\mathscr{P}_{\{3,4\}}(b_1,b_2,b_3)$ is integral
	\begin{align*}
		\sqrt{\Delta_{\{3,3,4\}}}&=\sqrt{-Q_{\{3,3\}}(b_1,b_2,b_3,b_1-b_2+b_3)	}\\
		&=\sqrt{b_1b_2+b_2b_3+b_3b_1-b_2^2}\\
		&=	\sqrt{\Delta_{\{4,3\}}}\in\mathbb Z
	\end{align*}
	again by Lemma \ref{lem:converse}.  By Corollary \ref{eq:integralitycond}, $\mathscr{P}_{\{3,3,4\}}(b_1,b_2,b_3,b_4)$ is also integral.
\end{proof}

\begin{figure}[H]
	\centering
	
	\begin{tikzpicture}
		\begin{scope}[xshift=-3.5cm]
			\node at (0,0) {\includegraphics[align=c,width=.45\textwidth,clip,trim= 0 30 30 60]{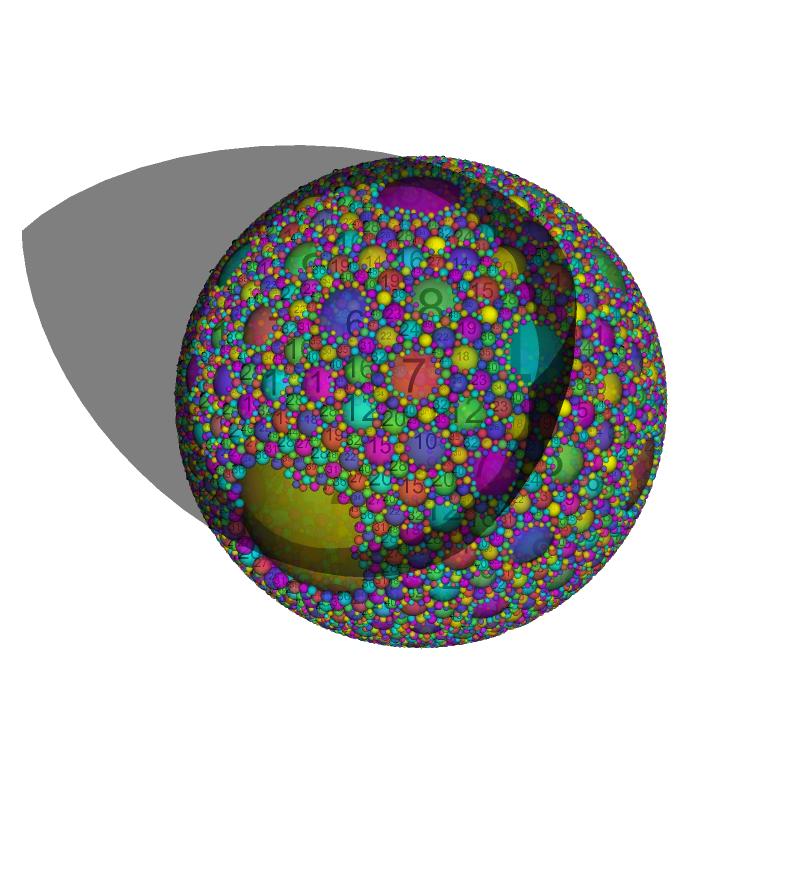}};				
		\end{scope}
		\begin{scope}[xshift=3.5cm]
			\node at (0,0) {\includegraphics[align=c,width=.45\textwidth,clip,trim= 0 30 30 60]{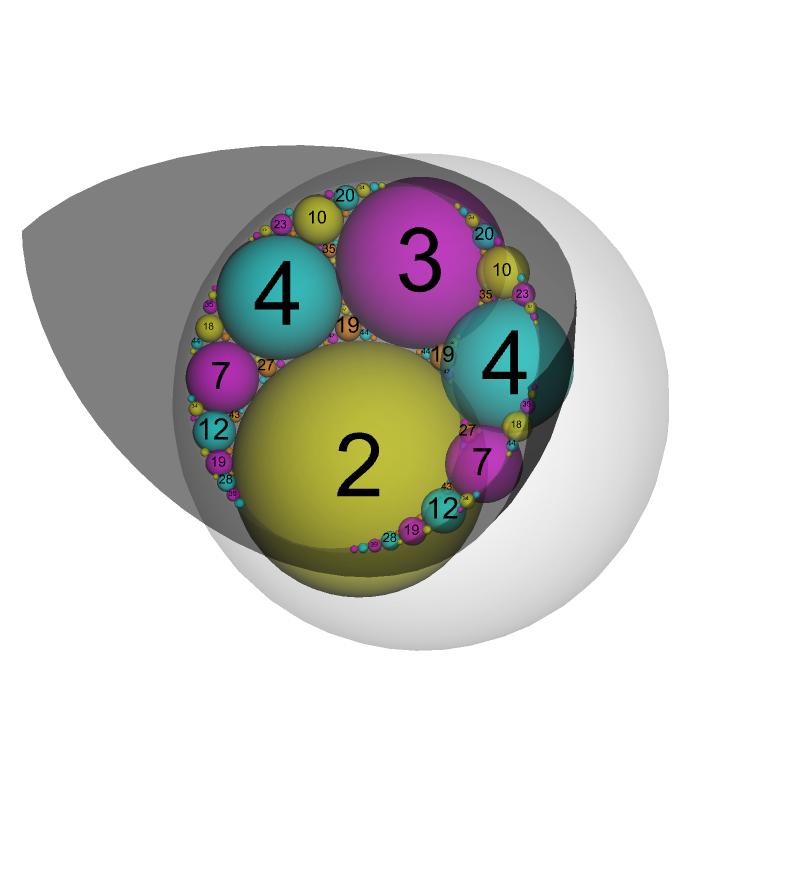}};
		\end{scope}

	\end{tikzpicture}
	\vspace{-1.75cm}

	\caption{
		(Left) The integral orthoplicial crystallographic packing $\mathscr P_{\{3,3,4\}}(-1,2,2,3)$ with a non-zero bend cutting sphere (in black) of a tetrehedral Apollonian section (right). The intersection yields the spherical integral Apollonian packing $\mathscr P_{\{3,3\}}(-1,2,3)$.
	}
	\label{fig:intersectingspheres}
\end{figure}

When the cutting sphere is not a plane, the circle packing obtained by the intersection can be thought of as a circle packing on the sphere, akin to the \textit{spherical} integral Apollonian packings introduced in \cite{lagarias2002beyond} and studied in \cite{eriksson2007apollonian}. Particularly for the tetrahedral sections, one obtains a \textit{spherical} analogue of Descartes' Theorem due to Mauldon \cite{mauldon1962sets}, which states that the bends of four mutually tangent circles on a sphere of bend $b_\Sigma$ are related by
\begin{align}
(b_1+b_2+b_3+b_4)^2	=2(b_1^2+b_2^2+b_3^2+b_4^2)+4b_\Sigma^2.
\end{align} 
In Figure \ref{fig:intersectingspheres}, we illustrate the spherical integral Apollonian packing with initial bends $(-1,2,3)$ obtained as the intersection of a cutting sphere of bend $b_\Sigma=1$ with a tetrahedral section of the integral orthoplicial crystallographic packing of initial bends $(-1,2,2,3)$. 
\vspace{-.25cm}

\subsection{The Möbius spectrum of the regular polytopes}\label{sec:mobspec}

In \cite{RR21_1}, the authors introduced an invariant of edge-scribable $d$-polytopes $\P$ with $d\ge3$ called the \textit{Möbius spectrum} $\mathfrak M(\P)$. This is defined as the multiset of the eigenvalues of the Gramian of any polytopal sphere packing $\S_\P$. A similar notion of the Möbius spectrum has been used in \cite{boyd,chenlabbe} for enumerating Boyd-Maxwell sphere packings. Due to the Möbius uniqueness of edge-scribable polytopes, $\mathfrak M(\P)$ does not depend on the packing.
 
 \medskip
 In this final section, we compute the Möbius spectrum of every regular polytope $\P$ in terms of the number of vertices and the canonical length. We shall use a lemma established in \cite{RR21_1} based on the Lorentzian model of the space of spheres. In this model, a sphere $S$ of $\widehat{\mathbb R^{d-1}}$ corresponds to a unique vector $\mathbf x_S$ of the Lorentz space $\mathbb R^{d,1}$ satisfying $\langle \mathbf x_S,\mathbf x_S\rangle=1$. The inversive coordinates and the inversive product of spheres correspond to the Cartesian coordinates and the inner product of the Lorentzian vectors, respectively. The Lorentzian vector $\mathbf x_v$ of a point $v\in\mathbb R^{d}$ outside $\mathbb S^{d-1}$ is the Lorentzian vector of the stereographic sphere of $v$. Additionally, for each $i=0,\ldots,d$ the \textit{Lorenzian barycenter} of a $i$-face $f$ of a $d$-polytope $\P\subset\mathbb R^{d}$ (including the case $f=\P$) whose vertices are outside $\mathbb S^{d-1}$, is $$\mathbf x_f:=\frac{1}{|V(f)|}\sum_{v\in V(f)}\mathbf x_{v}.$$
  
  \begin{lem}[\cite{RR21_1}]\label{lemma:lp}
	Let $\P$ be a uniform edge-scribed $d$-polytope with $d\ge 3$. For each face $f$ of $\P$
	\begin{align}
		\langle \mathbf x_f,\mathbf x_\P\rangle=-\ell_\P^{-2}.
	\end{align}
\end{lem}

\begin{thm}\label{thm:mobspec}
	For every $d\ge3$, the Möbius spectrum of every regular $d$-polytope $\P$ with $n> d$ vertices is
	\begin{align}
		\mathfrak M(\P)=(\frac nd(1+\ell_{\P}^{-2})_{(d)},-n\ell_{\P}^{-2},0_{(n-d-1)})
	\end{align}
	where $\ell_{\mathcal P}$ denotes the canonical length of $\mathcal P$.
\end{thm}
\begin{proof}
	Let $\S_\P$ be a polytopal $(d-1)$-sphere packing where $\P$ is a regular $d$-polytope with $d\ge3$ and $n$ vertices. 
	The spheres of $\S_\P$ correspond to a full-rank collection of $n$ unit vectors in $\mathbb R^{d,1}$. Let $\mathbf V$ be the matrix of the Cartesian coordinates of these vectors.  Since $\mathrm{Gram}(\S_\P)$ is a real symmetric matrix, its eigenvalues $\lambda_1,\ldots,\lambda_{n}$ are real. Moreover, the rank of $\mathrm{Gram}(\S_\P)$ is equal to the rank of  $\mathbf V$, so there are $\lambda_1,\ldots,\lambda_{d+1}$ non-zero eigenvalues and $\lambda_{d+2}=\ldots=\lambda_n=0$. By combining the Spectral Theorem for real symmetric matrices with the definition of $\mathrm{Gram}(\S_\P)$, we have the
\begin{equation}
	\mathbf {UDU}^{-1}=\mathbf{ V Q}_{d+1} \mathbf V^T
\end{equation} 
where $\mathbf U$ is an orthogonal matrix, $\mathbf D=\mathrm {diag}(\lambda_1,\ldots,\lambda_{n})$ and $\mathbf Q=\mathrm{diag} (1,\ldots,1,-1)$ of size $d+1$. Then, $\mathbf {D}=\mathrm{Gram}(\mathbf W)$ with $\mathbf W=\mathbf U^T\mathbf V$. This implies that $\mathbf {D}$ is the Gramian of a collection of $n$ vectors of $\mathbb R^{d,1}$ containing an orthogonal basis without light-like vectors of $\mathbb R^{d,1}$, so this basis must be made of $d$ space-like vectors and one time-like vector. Therefore, there are $d$ positive eigenvalues $\lambda_1,\ldots,\lambda_d$ and one negative eigenvalue $\lambda_{d+1}$.	The symmetry of $\P$ implies that $\lambda_1=\cdots=\lambda_d$. By adding all the rows of $\mathrm{Gram}(\mathcal S_\P)-\lambda_i I_n$ to the last one, we obtain a $n$-vector whose $i$-th entry is equal to
	\begin{align}
		\sum_{j=1}^{n}\langle \mathbf{x}_{v_i},\mathbf{x}_{v_j}\rangle-\lambda_i=\langle \mathbf{x}_{v_i},\sum_{j=1}^{n}\mathbf{x}_{v_j}\rangle-\lambda_i=\langle \mathbf{x}_{v_i},n\mathbf{x}_{\P}\rangle-\lambda_i=-n\ell_{\P}^{-2}-\lambda_i
	\end{align}
	by Lemma \ref{lemma:lp}. Therefore, $-n\ell_{\P}^{-2}$ is a negative root of the characteristic polynomial of $\mathrm{Gram}(\S_\P)$, so $\lambda_{d+1}=-n\ell_{\P}^{-2}$. By combining this with the equation $\mathrm{tr}(\mathrm{Gram}(\S_\P))=d\lambda_1+ \lambda_{d+1} =n$,	we obtain that $\lambda_1=\frac nd(1+\ell_{\P}^{-2})$.
\end{proof}
In Table \ref{tab:mobspec}, we present the canonical lengths (adapted  from  \cite{coxeter1973regular}) and the Möbius spectrum of every regular polytope in every dimension equal or greater than 3. 
\begin{table}[H]
	\small
	\centering
	\begin{tabular}{lcccrrl}
		Dim.&Regular polytope& Schläfli symbol&Canonical length&	\multicolumn{3}{c}{Möbius spectrum}\\
		\toprule
		\multirow{4}{*}{$d\ge3$} 	
		&$d$-Simplex $\T^d$&
		$\{3_{(d-1)}\}$ &$\sqrt{\frac{d+1}{d-1}}$&$2_{(d)},$&$-d+1$ \\
		\cmidrule{2-7}  
		&$d$-Cross polytope $\O^d$
		&$\{3_{(d-2)},4\}$&$1$ &$4_{(d)}$,&$-2d$,&$0_{(d-1)}$\\
		\cmidrule{2-7}   
		&$d$-Cube $\C^d$
		&$\{4,3_{(d-2)}\}$&$\sqrt{\frac{1}{d-1}}$&$2^{d}_{(d)}$,&$-2^{d} (d-1)$,&$0_{(2^{d}-d-1)}$\\
		\hline
		\multirow{3}{*}{$d=3$} 
		&Icosahedron $\I^3$
		&$\{3,5\}$&$\varphi^{-1}$&$4(1+\varphi^2)_{(3)}$,&$-12\varphi^2$,&$0_{(8)}$\\ 
		\cmidrule{2-7} 
		&Dodecahedron $\D^3$&$\{5,3\}$ &$\varphi^{-2}$&$20\varphi^2_{(3)}$,&$-20\varphi^4$,&$0_{(16)}$\\                         
		\hline
		\multirow{4}{*}{$d=4$}
		&24-cell $\R^4$&$\{3,4,3\}$ &$3^{-1/2}$& $24_{(4)}$,&$-72$,&$0_{(19)}$ \\
		\cmidrule{2-7}
		&600-cell $\I^4$&
		$\{3,3,5\}$&$5^{-1/4}\varphi^{-3/2}$&$120\varphi^2_{(4)}$,&$-120(1-4\varphi^2)$,&$0_{(115)}$\\
		\cmidrule{2-7}
		&120-cell $\D^4$&
		$\{5,3,3\}$&$3^{-1/2}\varphi^{-3}$&$1200\varphi^4_{(4)}$,&$-1800\varphi^6$,&$0_{(595)}$	
	\end{tabular}
	\caption{Notations, Schläfli symbol, canonical lengths and Möbius spectrum of the regular $d$-polytopes for $d\ge3$. The indices in parentheses indicate multiplicity.}
	\label{tab:mobspec}
\end{table}

	\printbibliography[	title={References}	] 

@article{apoGI,
author = {Graham, R.
and Lagarias, J. C.
and Mallows, C. L.
and Wilks, A. R.
and Yan, C. H.},
title={{A}pollonian Circle Packings: Geometry and Group Theory I. The {A}pollonian Group},
journal={Discrete {\&} Computational Geometry},
year={2005},
day={01},
volume={34},
number={4},
pages={547-585},
doi={10.1007/s00454-005-1196-9},
url={https://doi.org/10.1007/s00454-005-1196-9}
}

@article{aporingpacks,
  title={Apollonian ring packings},
  author={Bolt, Adrian and Butler, Steve and Hovland, Espen},
  journal={Connections in Discrete Mathematics: A Celebration of the Work of Ron Graham},
  pages={283},
  year={2018},
  publisher={Cambridge University Press}
}

@Article{baragar2018higher,
author={Baragar, A.},
title={Higher dimensional {A}pollonian packings, revisited},
journal={Geometriae Dedicata},
year={2018},
month={Aug},
day={01},
volume={195},
number={1},
pages={137-161},
issn={1572-9168},
doi={10.1007/s10711-017-0280-7},
url={https://doi.org/10.1007/s10711-017-0280-7}
}

@article{bogachev2024kleinian,
	title={Kleinian sphere packings, reflection groups, and arithmeticity},
	author={Bogachev, Nikolay and Kolpakov, Alexander and Kontorovich, Alex},
	journal={Mathematics of Computation},
	volume={93},
	number={345},
	pages={505--521},
	year={2024}
}

@article{boyd,
author = {D. W. Boyd},
title = {{A new class of infinite sphere packings.}},
volume = {50},
journal = {Pacific Journal of Mathematics},
number = {2},
publisher = {Pacific Journal of Mathematics, A Non-profit Corporation},
pages = {383 -- 398},
year = {1974},
doi = {pjm/1102913226},
URL = {https://doi.org/}
}

@article{chait2020taxonomy,
	title={A taxonomy of crystallographic sphere packings},
	author={Chait-Roth, D. and Cui, A. and Stier, Z.},
	journal={Journal of Number Theory},
	volume={207},
	pages={196--246},
	year={2020},
	publisher={Elsevier}
}

@article{chen2016even,
  title={Even More Infinite Ball Packings from {L}orentzian Root Systems},
  author={Chen, H.},
  journal={The Electronic Journal of Combinatorics},
  doi={10.37236/4989},
  year={2016},
  volume={23},
  number={3}
}

@Article{chenlabbe,
author={Chen, H.
and Labb{\'e}, J. P.},
title={{L}orentzian {C}oxeter systems and {B}oyd--{M}axwell ball packings},
journal={Geometriae Dedicata},
year={2015},
month={Feb},
day={01},
volume={174},
number={1},
pages={43-73},
issn={1572-9168},
doi={10.1007/s10711-014-0004-1},
url={https://doi.org/10.1007/s10711-014-0004-1}
}

@book{coxeter1973regular,
  title={Regular Polytopes},
  author={Coxeter, H. S. M.},
  isbn={9780486614809},
  lccn={73084364},
  series={Dover books on advanced mathematics},
  year={1973},
  publisher={Dover Publications}
}

@article{Dias2014TheLP,
  title={The {L}ocal-{G}lobal {P}rinciple for Integral Generalized {A}pollonian Sphere Packings},
  author={Dias, D.},
  journal={arXiv: Number Theory},
  year={2014}
}

@article{eriksson2007apollonian,
	title={Apollonian circle packings: number theory II. Spherical and hyperbolic packings},
	author={Eriksson, Nicholas and Lagarias, Jeffrey C},
	journal={The Ramanujan Journal},
	volume={14},
	number={3},
	pages={437--469},
	year={2007},
	publisher={Springer}
}

@article{guettler,
author = {Guettler, G. and Mallows, C. L.},
year = {2008},
month = {01},
pages = {},
title = {A generalization of {A}pollonian packing of circles},
volume = {1},
journal = {Journal of Combinatorics},
doi = {10.4310/JOC.2010.v1.n1.a1}
}

@article{kapovich2023superintegral,
	title={On superintegral kleinian sphere packings, bugs, and arithmetic groups},
	author={Kapovich, M. and Kontorovich, A.},
	journal = {Journal für die Reine und Angewandte Mathematik (Crelles Journal)},
	number={0},
	year={2023},
	publisher={De Gruyter}
}

@inproceedings{kertzer2024local,
	title={The local-global conjecture for Apollonian circle packings is false},
	author={Kertzer, Clyde and Haag, Summer and Stange, Katherine E and Rickards, James},
	booktitle={2024 Joint Mathematics Meetings (JMM 2024)},
	 year={2024},
	organization={AMS}
}

@article{kontorovich2019soddy,
	title = {The local-global principle for integral Soddy sphere packings},
	journal = {Journal of Modern Dynamics},
	volume = {15},
	number = {0},
	pages = {209-236},
	year = {2019},
	author = {A. Kontorovich},
}

@article {KontorovichNakamura,
	author = {Kontorovich, A. and Nakamura, K.},
	title = {Geometry and arithmetic of crystallographic sphere packings},
	volume = {116},
	number = {2},
	pages = {436--441},
	year = {2019},
	doi = {10.1073/pnas.1721104116},
	publisher = {National Academy of Sciences},
	issn = {0027-8424},
	journal = {Proceedings of the National Academy of Sciences}
}

@article{lagarias2002beyond,
	title={Beyond the Descartes circle theorem},
	author={Lagarias, Jeffrey C and Mallows, Colin L and Wilks, Allan R},
	journal={The American Mathematical Monthly},
	volume={109},
	number={4},
	pages={338--361},
	year={2002},
	publisher={Taylor \& Francis}
}

@article{martin2024geometric,
	title={A geometric study of circle packings and ideal class groups},
	author={Martin, Daniel E},
	journal={Discrete \& Computational Geometry},
	volume={72},
	number={1},
	pages={181--208},
	year={2024},
	publisher={Springer}
}

@misc{Mathematica,
	author = {,Wolfram Research Inc.},
	title = {Mathematica, {V}ersion 13.1},
	url = {https://www.wolfram.com/mathematica},
	note = {Champaign, IL, 2022}
}

@article{mauldon1962sets,
	title={Sets of equally inclined spheres},
	author={Mauldon, JG},
	journal={Canadian Journal of Mathematics},
	volume={14},
	pages={509--516},
	year={1962},
	publisher={Cambridge University Press}
}

@article{maxwell,
title = {Sphere packings and hyperbolic reflection groups},
journal = {Journal of Algebra},
volume = {79},
number = {1},
pages = {78-97},
year = {1982},
issn = {0021-8693},
doi = {https://doi.org/10.1016/0021-8693(82)90318-0},
url = {https://www.sciencedirect.com/science/article/pii/0021869382903180},
author = {Maxwell, G.}
}

@misc{nakamura2014localglobal,
      title={The local-global principle for integral bends in orthoplicial {A}pollonian sphere packings}, 
      author={K. Nakamura},
      year={2014},
      eprint={1401.2980},
      archivePrefix={arXiv},
      primaryClass={math.NT}
}

@article{RR20,
title = {Ball packings for links},
journal = {European Journal of Combinatorics},
volume = {96},
pages = {103351},
year = {2021},
issn = {0195-6698},
doi = {https://doi.org/10.1016/j.ejc.2021.103351},
url = {https://www.sciencedirect.com/science/article/pii/S0195669821000433},
author = {J. L. {Ramírez Alfonsín} and I. Rasskin}
}

@article{RR21_1,
      title={A polytopal generalization of {A}pollonian packings and {D}escartes' theorem}, 
      author={J. L. {Ramírez Alfonsín} and I. Rasskin},
      year={2021},
      eprint={2107.09432},
      archivePrefix={arXiv},
      primaryClass={math.CO}
}

@article{RR2024links,
	title={Links in orthoplicial Apollonian packings},
	author={J. L. {Ramírez Alfonsín} and I. Rasskin},
	journal={European Journal of Combinatorics},
	volume={122},
	pages={104017},
	year={2024},
	publisher={Elsevier}
}

@phdthesis{rasskin_thesis,
  TITLE = {{A polytopal approach to Apollonian packings and discrete knotted structures}},
  AUTHOR = {Rasskin, I.},
  URL = {https://hal.archives-ouvertes.fr/tel-03480927},
  SCHOOL = {{Universit{\'e} de Montpellier}},
  YEAR = {2021},
  MONTH = Dec,
  TYPE = {Theses},
  HAL_ID = {tel-03480927},
  HAL_VERSION = {v1},
}

@inproceedings{schulte04,
  title={Symmetry of Polytopes and Polyhedra},
author = {Schulte, E.},
  booktitle={Handbook of Discrete and Computational Geometry, 2nd Ed.},
  year={2004}
}

@article{SHEYDVASSER201941,
title = {Quaternion orders and sphere packings},
journal = {Journal of Number Theory},
volume = {204},
pages = {41-98},
year = {2019},
issn = {0022-314X},
doi = {https://doi.org/10.1016/j.jnt.2019.03.014},
url = {https://www.sciencedirect.com/science/article/pii/S0022314X19301192},
author = {Sheydvasser, A.}
}

@Article{soddy1936,
author={Soddy, F.},
title={{The Kiss Precise}},
journal={Nature},
year={1936},
month={Jun},
day={01},
volume={137},
number={3477},
pages={1021-1021},
issn={1476-4687},
doi={10.1038/1371021a0},
url={https://doi.org/10.1038/1371021a0}
}

@article{springborn2005unique,
  title={{A unique representation of polyhedral types. Centering via M{\"o}bius transformations}},
  author={Springborn, B. A.},
  journal={Mathematische Zeitschrift},
  volume={249},
  number={3},
  pages={513--517},
  year={2005},
  publisher={Springer}
}

@article{stange2015bianchi,
author = {Stange, K. E.},
year = {2015},
month = {05},
pages = {},
title = {The {A}pollonian structure of {B}ianchi groups},
volume = {370},
journal = {Transactions of the American Mathematical Society},
doi = {10.1090/tran/7111}
}

@article{Zhang+2018+71+110,
	author = {Zhang, X. },
	doi = {doi:10.1515/crelle-2015-0042},
	url = {https://doi.org/10.1515/crelle-2015-0042},
	title = {On the local-global principle for integral {A}pollonian 3-circle packings: },
	journal = {Journal für die Reine und Angewandte Mathematik (Crelles Journal)},
	number = {737},
	volume = {2018},
	year = {2018},
	pages = {71--110}
}
	
	\appendix

\section{Apollonian arrangements and sections of the regular 4-polytopes}\label{sec:appendix}

\subsection{Simplex $\{3,3,3\}$}
\begin{figure}[H]
	\includegraphics[align=c,width=.22\textwidth]{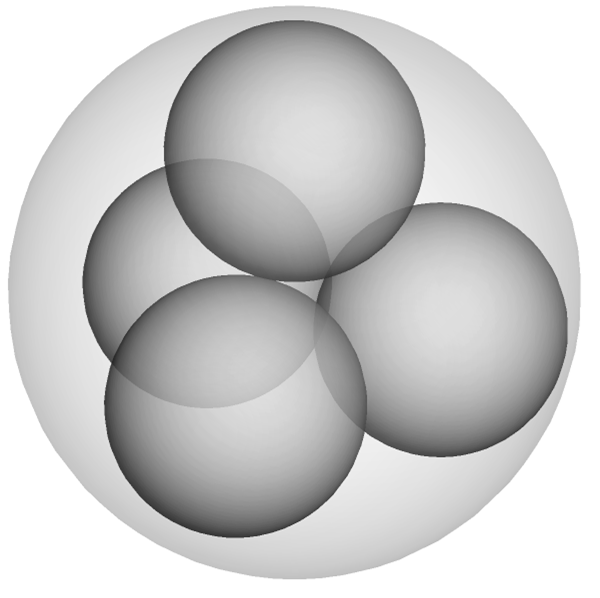} \hspace{.1cm}  
	\includegraphics[align=c,width=.22\textwidth]{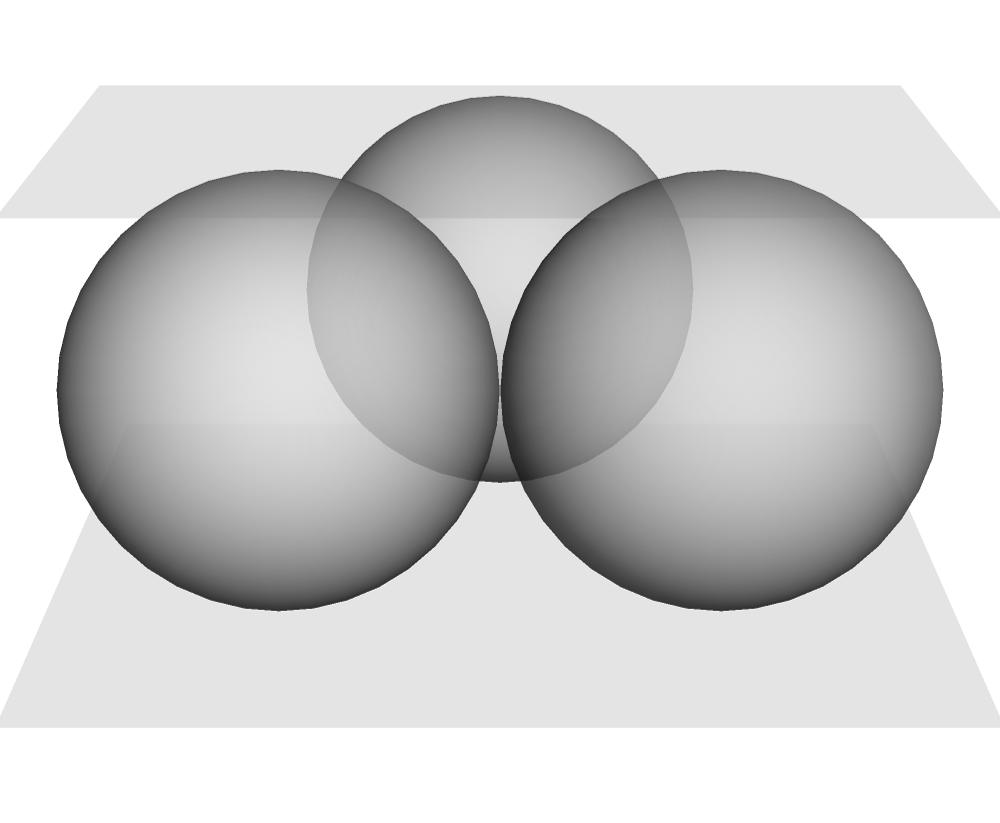} \hspace{.1cm} 
	\includegraphics[align=c,width=.22\textwidth]{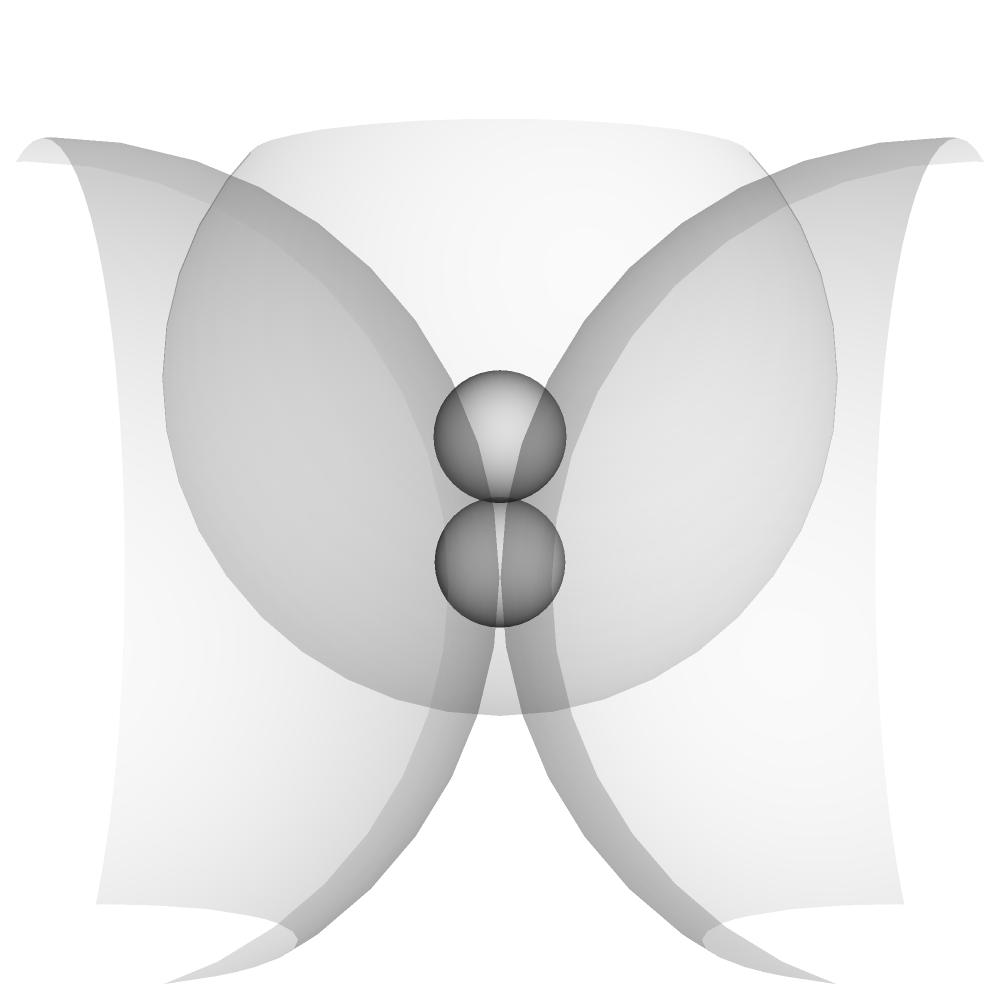} \hspace{.1cm} 
	\includegraphics[align=c,width=.22\textwidth]{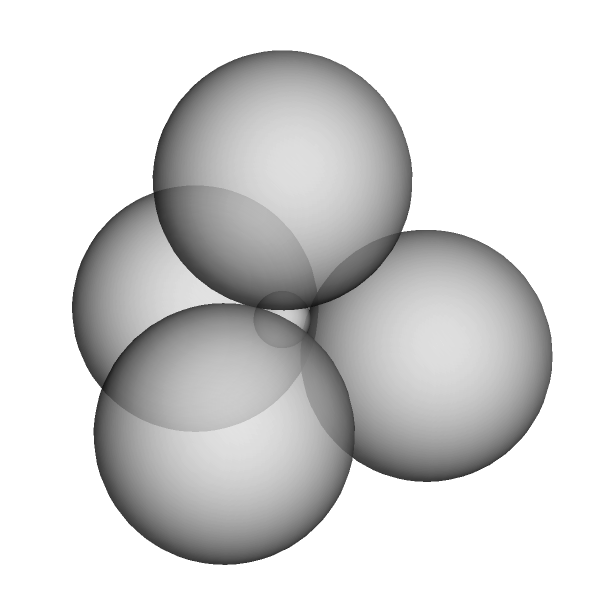}
	\caption{(From left to right) Arrangement projections of a vertex-centered, edge-centered, ridge-centered and facet-centered canonical $4$-simplex.}
	\label{fig:centeretetr}
\end{figure}

In Figure \ref{fig:centeretetr}, we show four simplicial sphere packings obtained by the arrangement projections of face-centered canonical $4$-simplices.

The full symmetry group $\Gamma_{\{3,3,3\}}$ is isomorphic to $\langle \mathbf{R}_1,\mathbf R_2,\mathbf{R}_3,\mathbf R_4,\mathbf S_f\rangle <\mathrm{SL}_5(\mathbb Z)$ where
\begin{equation}
\begin{gathered}
	\mathbf{R}_1=\left(\begin{array}{ccccc}
		& 1 &   &   &   \\
		1 &   &   &   &   \\
		&   & 1 &   &   \\
		&   &   & 1 &   \\
		&   &   &   & 1 \\
	\end{array}\right)\quad
	\mathbf{R}_2=\left(
\begin{array}{ccccc}
	1 &   &   &   &   \\
	&   & 1 &   &   \\
	& 1 &   &   &   \\
	&   &   & 1 &   \\
	&   &   &   & 1 \\
\end{array}
\right)\quad
\mathbf{R}_3=\left(
\begin{array}{ccccc}
	1 &   &   &   &   \\
	& 1 &   &   &   \\
	&   &   & 1 &   \\
	&   & 1 &   &   \\
	&   &   &   & 1 \\
\end{array}
\right)\\
\mathbf R_4=
\left(
\begin{array}{ccccc}
	1 &   &   &   &   \\
	& 1 &   &   &   \\
	&   & 1 &   &   \\
	&   &   &   & 1 \\
	&   &   & 1 &   \\
\end{array}
\right)\quad
\mathbf S_f=\left(
\begin{array}{ccccc}
	1 &   &   &   &   \\
	& 1 &   &   &   \\
	&   & 1 &   &   \\
	&   &   & 1 &   \\
	1 & 1 & 1 & 1 & -1 \\
\end{array}
\right)
\end{gathered}
\end{equation}	

\begin{figure}[H]
	\centering
	\begin{tikzpicture}[scale=2.5] 
		\node at (-3.2,0) {	\includegraphics[clip,trim=0 60 0 60,  align=c,width=.49\textwidth]{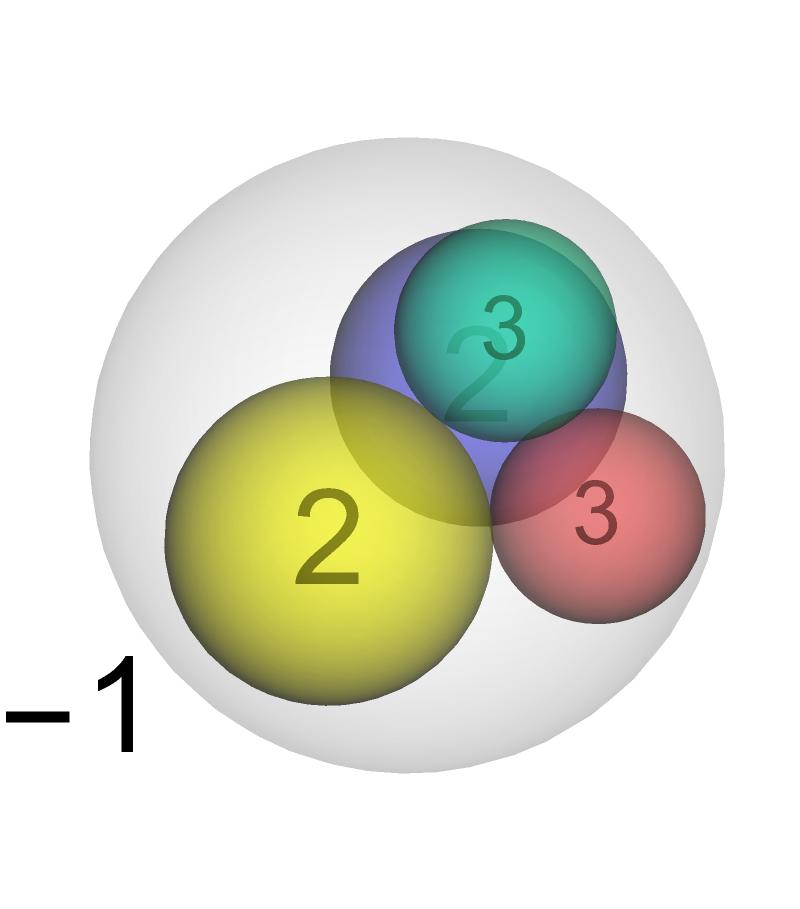}};
		\node at (0,0) {	\includegraphics[ align=c,width=.4\textwidth]{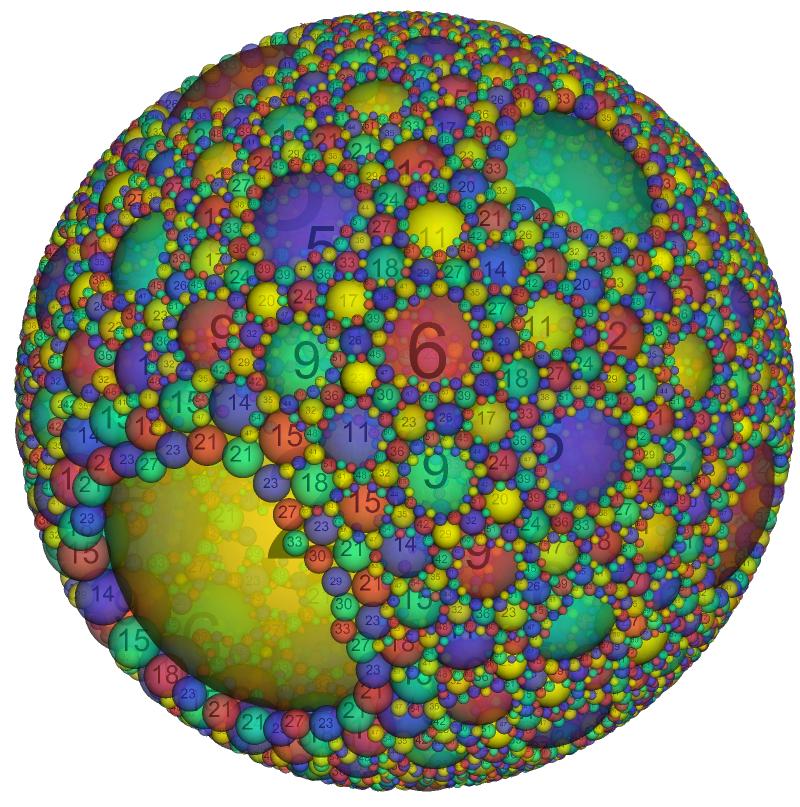}} ;
	\end{tikzpicture}
	
	\caption{
		A simplicial sphere packing satisfying the integrality condition (left) and the corresponding integral crystallographic packing  $\mathscr P_{\{3,3,3\}}(-1,2,2,3)$ (right).
	}
	\label{fig:apo333}
\end{figure}

Any fundamental bend vector $\mathbf b=(b_1,b_2,b_3,b_4,b_5)^T$ of a simplicial sphere packing 
satisfies the quadratic equation $\mathbf{b}^T\mathbf Q_{\{3,3,3\}}\mathbf{b}=0$ for the bisymmetric matrix
\begin{align}\label{eq:descartes333}
	\mathbf Q_{\{3,3,3\}}=
	\left(
	\begin{array}{ccccc}
		2 & -1 & -1 & -1 & -1 \\
		\ast & 2 & -1 & -1 & \ast \\
		\ast &\ast & 2 & \ast & \ast \\
		\ast & \ast & \ast & \ast & \ast \\
		\ast & \ast & \ast & \ast & \ast \\
	\end{array}
	\right)
\end{align} 
The latter is equivalent to the classic Soddy's quadratic equation \cite{soddy1936}
\begin{align}
	(b_1+b_2+b_3+b_4+b_5)^2= 3(b_1^2+b_2^2+b_3^2+b_4^2+b_5^2).
\end{align}
The integrality condition of Corollary \ref{cor:integrality} states that if  $	b_1,b_2,b_3,b_4,\sqrt{\Delta_{\{3,3,3\}}}\in\mathbb Z$ where
\begin{align}\label{eq:inttetra}
\Delta_{\{3,3,3\}}
=3\left((b_1+b_2+b_3+b_4)^2-2(b_1^2+b_2^2+b_3^2+b_4^2)\right)
\end{align}	
then the simplicial crystallographic packing $\mathscr P_{\{3,3,3\}}(b_1,b_2,b_3,b_4)$ is integral (see Figure \ref{fig:apo333}).

\begin{figure}[H]
	\centering
	\begin{tabular}{cc}
		&\begin{tikzpicture}[scale=2] 
			\node at (-3,0) {		\begin{tikzpicture}[scale=1.4] 		
		
		\newcommand\xmin{1}
\newcommand\xmax{2.414}

\newcommand\ymin{1}
\newcommand\ymax{2}
		
		\draw[very thin] (-\xmin,1) -- (\xmax,1);
		\draw[very thin] (0,0.5) circle (.5cm);
		\draw[very thin] (1.,0.5) circle (.5cm) ;
		
		\draw[very thick,red] (-\xmin,0) -- (\xmax,0) node[pos=.95, below] {$S_v$};
		
\draw[blue,thick] (-\xmin,0.5) -- (\xmax,.5)	node[pos=.95, above] {$R_1$};
\draw[thick,blue] (0,0) circle (1cm) node [xshift=-1.3cm,yshift=-1.3cm] {$R_2$};
\draw[blue,thick] (0.5,-\ymin) -- (0.5,\ymax) node[pos=.95,right] {$R_3$};
\draw[blue,thick] (0,-\ymin) -- (0,\ymax) node[pos=.95,left] {$S_f$};

	\end{tikzpicture}};
			
			\node {	\includegraphics[align=c,width=0.4\textwidth]{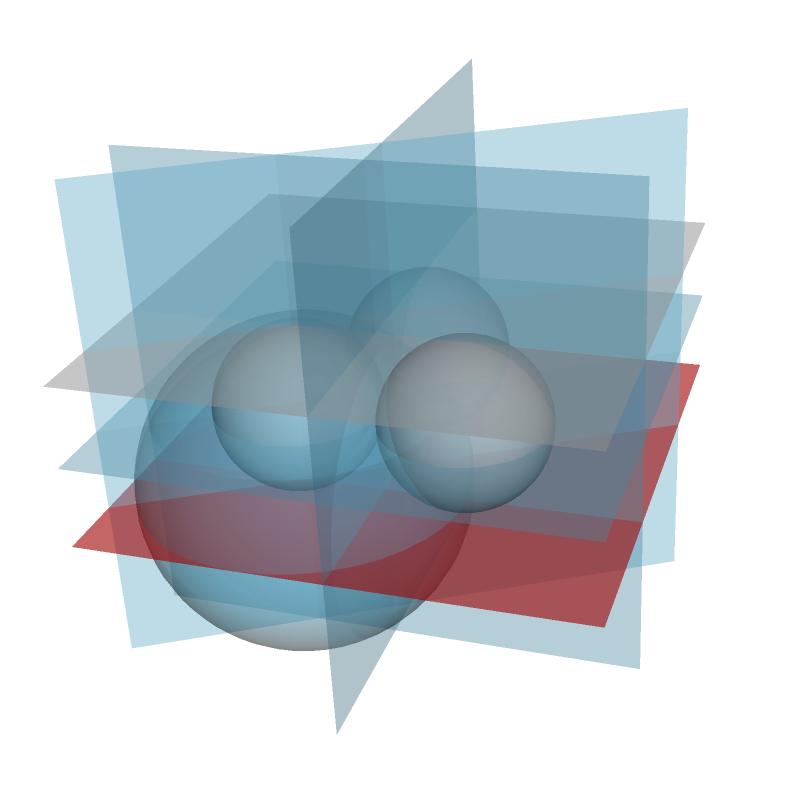} };
			\node at (1.28,0.15) [red] {$S_v$};
			\node at (-1.15,1.09) [blue] {$S_f$};
			\node at (1.3,.4) [blue] {$R_1$};
			\node at (-.55,-1.1) [blue] {$R_2$};
			\node at (.35,1.4) [blue] {$R_3$};
			\node at (1.2,1.2) [blue] {$R_4$};
		\end{tikzpicture}
		
	\end{tabular}
	\vspace{-.5cm}
	\caption{
		The strip packing with the fundamental symmetries of the tetrahedron (left) and the $4$-simplex (right).
	}
	\label{fig:symaposec333}
\end{figure}

\begin{figure}[H]
	\centering
	
	\begin{tikzpicture}
		\begin{scope}
			\begin{scope}[xshift=-5.5cm]
				\node at (0,0) {\includegraphics[trim=0 0 40 60,clip,align=c,width=0.35\textwidth]{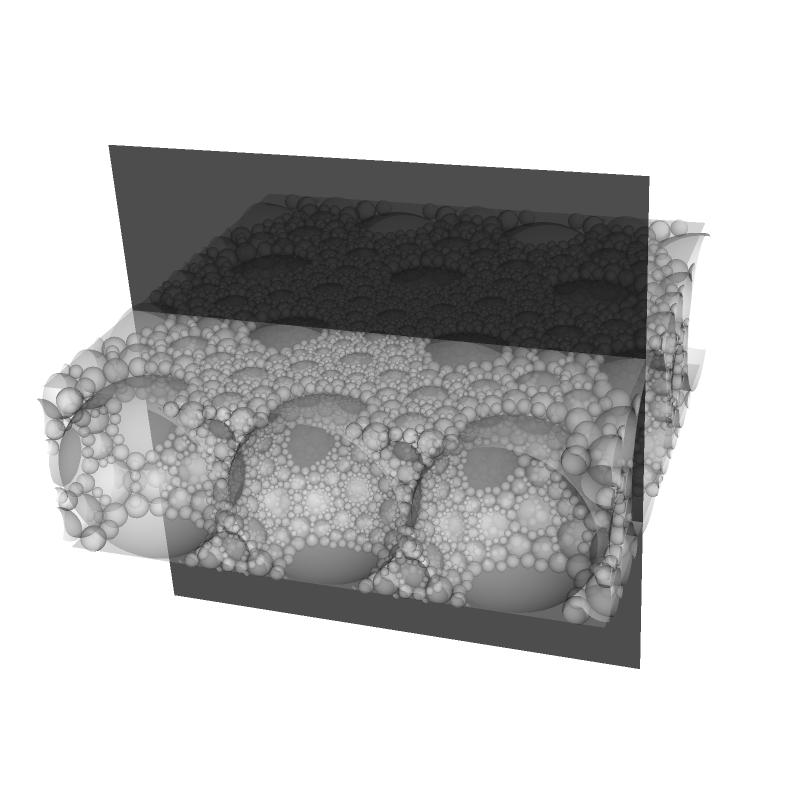}};
				
			\end{scope}
			\begin{scope}[xshift=-0cm]
				\node at (0,0) {\includegraphics[trim=0 0 40 60,clip,align=c,width=0.35\textwidth]{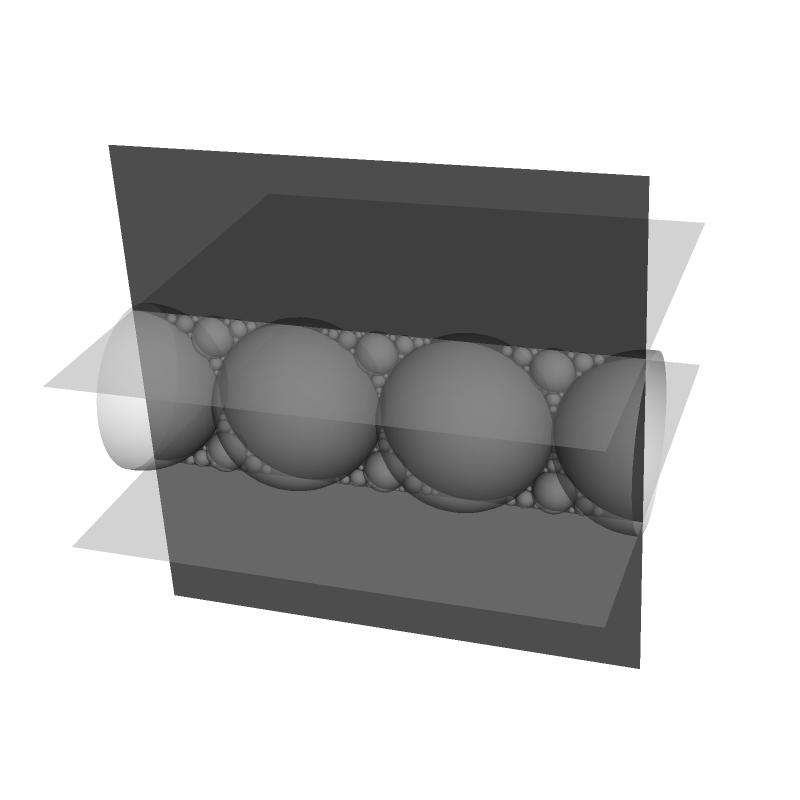}};				
			\end{scope}
			\begin{scope}[xshift=5.cm]
				\node at (0,0) {\includegraphics[align=c,width=0.25\textwidth]{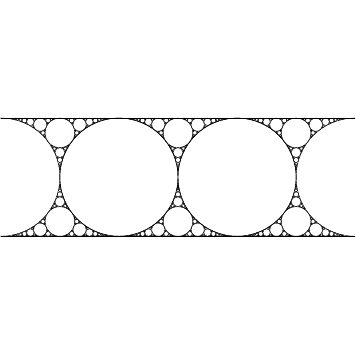}};
				
			\end{scope}
		\end{scope}

	\end{tikzpicture}
	
	\vspace{-1cm}
	
	\caption{
(From left to right) The simplicial crystallographic packing $\mathscr{P}_{\{3,3,3\}}$ with a cutting sphere $\Sigma_{\{3,3,3\}}^{\{3,3\}}$,
the tetrahedral Apollonian section $\mathscr S_{\{3,3,3\}}^{\{3,3\}}$ with $\Sigma_{\{3,3,3\}}^{\{3,3\}}$, and the tetrahedral crystallographic packing $\mathscr{P}_{\{3,3\}}$. 		
	}
	\label{fig:sectionsgeo333}
\end{figure}

 Every simplicial crystallographic packing contains a tetrahedral Apollonian section $\mathscr{S}^{\{3,3\}}_{\{3,3,3\}}$. In the strip packing, the cutting sphere is $\Sigma_{\{3,3,3\}}^{\{3,3\}}=S_f$ (see Figures \ref{fig:symaposec333}, \ref{fig:sectionsgeo333}). The homomorphism between the full symmetry groups is given by
\begin{align}\label{eq:hom333}
		\begin{tabular}{cccccc}
		&$\Gamma_{\{3,3\}}$&$\xrightarrow{\phi_{\{3,3,3\}}^{\{3,3\}}}$&$\Gamma_{\{3,3,3\}}$\\
		&$r_1$&$\longmapsto$& $r_1$\\
		&$r_2$&$\longmapsto$& $r_2$\\
		&$r_3$&$\longmapsto$& $r_3$\\
		&$s_f$&$\longmapsto$& $(r_4s_f)^3$\\
	\end{tabular}
\end{align}


\subsection{Orthoplex $\{3,3,4\}$}\label{sec:334}
\begin{figure}[H]
	\includegraphics[align=c,width=.22\textwidth]{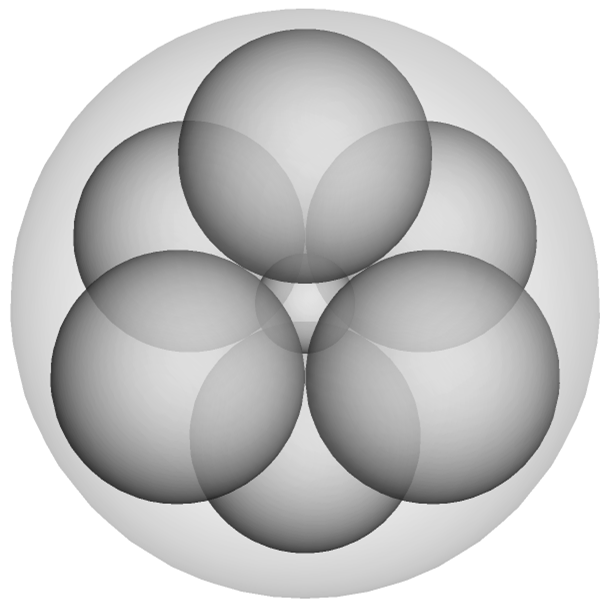} \hspace{.1cm}  
	\includegraphics[align=c,width=.22\textwidth]{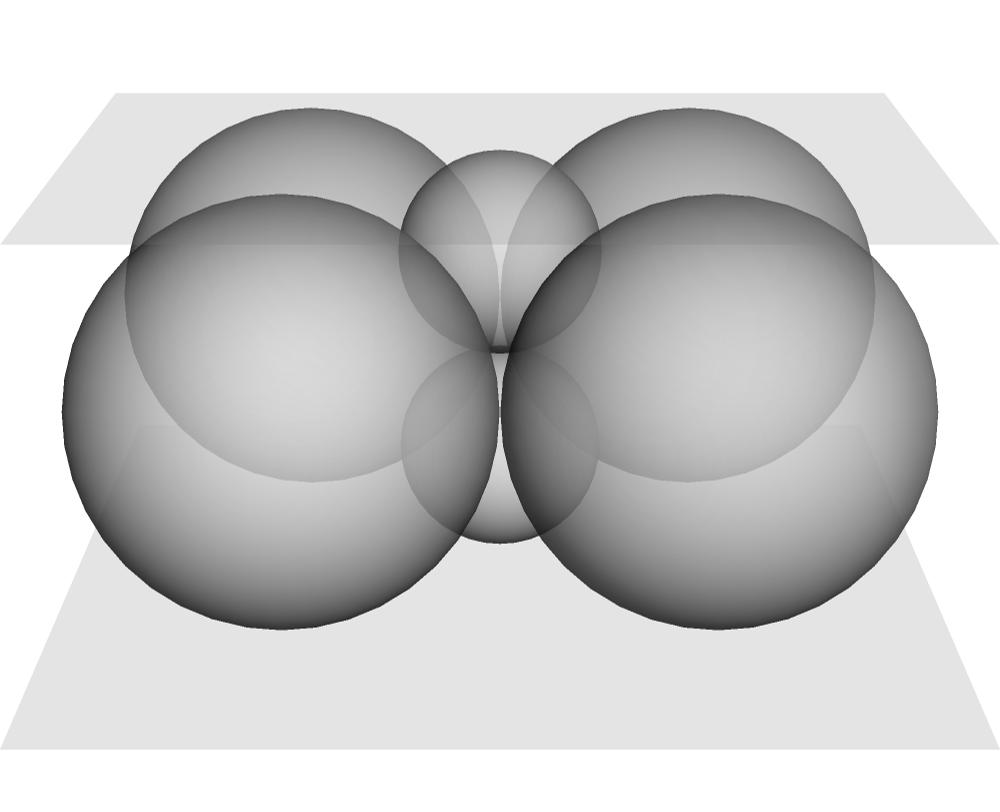} \hspace{.1cm} 
	\includegraphics[align=c,width=.22\textwidth]{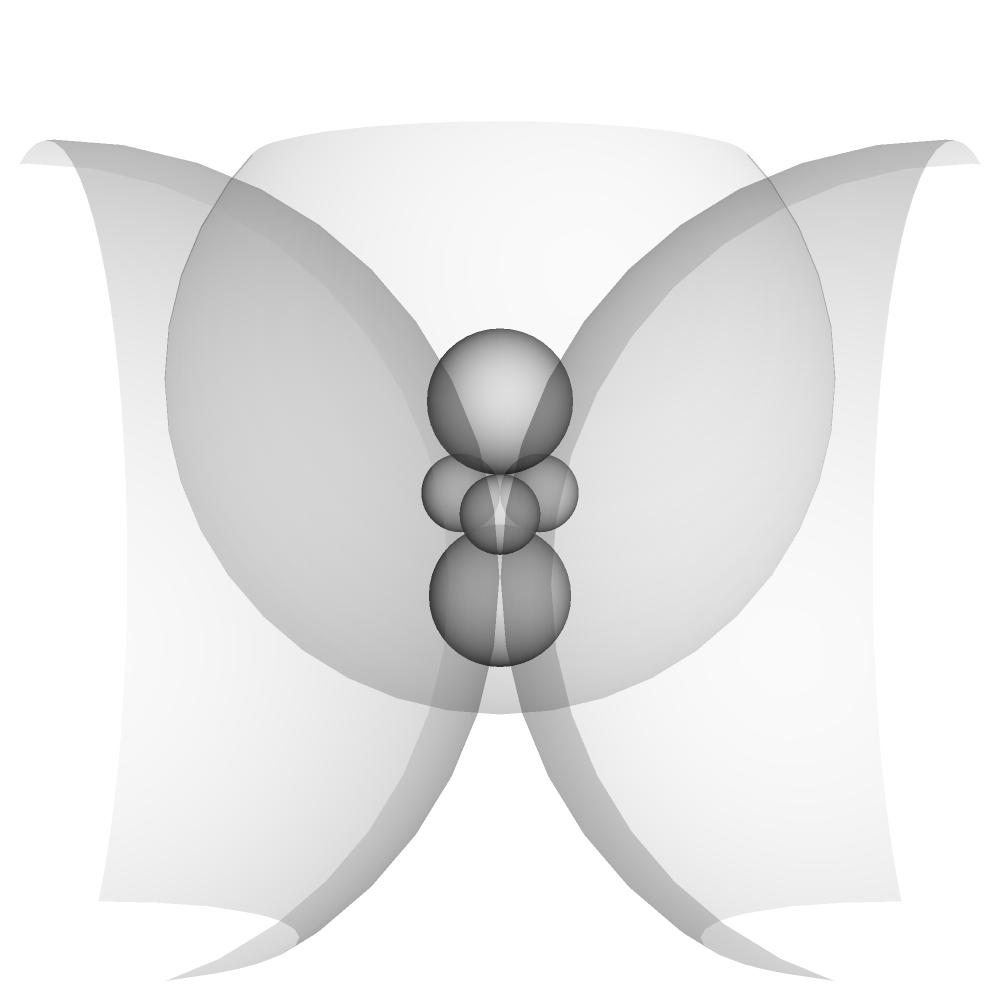} \hspace{.1cm} 
	\includegraphics[align=c,width=.22\textwidth]{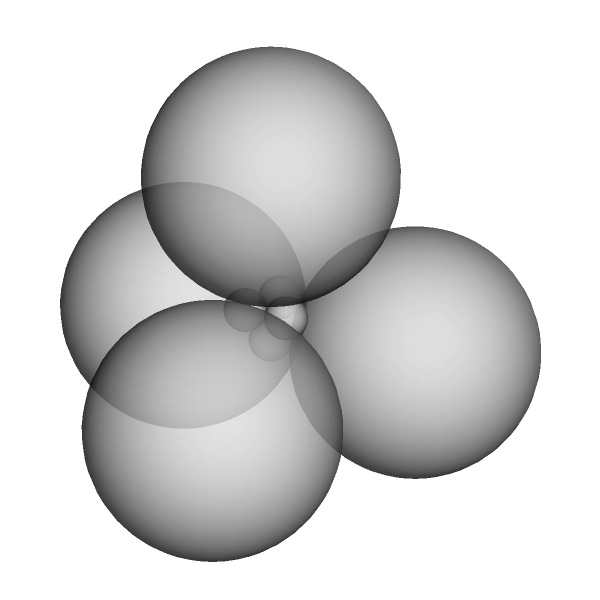}\\
	\caption{(From left to right) Arrangement projections of a vertex-centered, edge-centered, ridge-centered and facet-centered canonical orthoplex.}
	\label{fig:centeredhoct}
\end{figure}

In Figure \ref{fig:centeredhoct}, we show four orthoplicial sphere packings obtained by the arrangement projections of face-centered canonical orthoplices.

The full symmetry group $\Gamma_{\{3,3,4\}}$ is isomorphic to $\langle \mathbf{R}_1,\mathbf R_2,\mathbf{R}_3,\mathbf R_4,\mathbf S_f\rangle <\mathrm{SL}_5(\mathbb Z)$ where
\begin{equation}
\begin{gathered}
	\mathbf{R}_1= \left(
	\begin{array}{ccccc}
		& 1 &   &   &   \\
		1 &   &   &   &   \\
		&   & 1 &   &   \\
		&   &   & 1 &   \\
		1 & -1 &   &   & 1 \\
	\end{array}
	\right)\quad
	\mathbf{R}_2= \left(
	\begin{array}{ccccc}
		1 &   &   &   &   \\
		&   & 1 &   &   \\
		& 1 &   &   &   \\
		&   &   & 1 &   \\
		&   &   &   & 1 \\
	\end{array}
	\right)\quad
	\mathbf{R}_3= \left(
	\begin{array}{ccccc}
		1 &   &   &   &   \\
		& 1 &   &   &   \\
		&   &   & 1 &   \\
		&   & 1 &   &   \\
		&   &   &   & 1 \\
	\end{array}
	\right)\\
	\mathbf{R}_4=	 \left(
	\begin{array}{ccccc}
		1 &   &   &   &   \\
		& 1 &   &   &   \\
		&   & 1 &   &   \\
		1 &   &   & -1 & 1 \\
		&   &   &   & 1 \\
	\end{array}
	\right)\quad
	\mathbf S_f= \left(
	\begin{array}{ccccc}
		1 &   &   &   &   \\
		& 1 &   &   &   \\
		&   & 1 &   &   \\
		&   &   & 1 &   \\
		& 2 & 2 & 2 & -1 \\
	\end{array}
	\right)
\end{gathered}
\end{equation}	

Any fundamental bend vector $\mathbf B=(b_1,b_2,b_3,b_4,b_5)^T$ of an orthoplicial sphere packing 
satisfies the quadratic equation $\mathbf{b}^T\mathbf Q_{\{3,3,4\}}\mathbf{b}=0$ for the bisymmetric matrix
\begin{align}\label{eq:descartes334}
	\mathbf Q_{\{3,3,4\}}=
	\left(
	\begin{array}{ccccc}
		2 & -2 & -2 & -2 & 0 \\
		\ast & 4 & 0 & 0 & \ast \\
		\ast &\ast & 4 & \ast & \ast \\
		\ast & \ast & \ast & \ast & \ast \\
		\ast & \ast & \ast & \ast & \ast \\
	\end{array}
	\right)
\end{align} 

The latter  is equivalent to the following quadratic equation\footnote{Also equivalent to the quadratic equation given in \cite{Dias2014TheLP,nakamura2014localglobal}. }.
\begin{align}
	(b_1-b_5)^2+(b_1-2b_2+b_5)^2+(b_1-2b_3+b_5)^2+(b_1-2b_4+b_5)^2=2(b_1+b_5)^2.
\end{align}

\begin{figure}[H]
	\centering
	\begin{tikzpicture}[scale=2.5] 
		\node at (-3.2,0) {	\includegraphics[clip,trim=0 60 0 60, align=c,width=.49\textwidth]{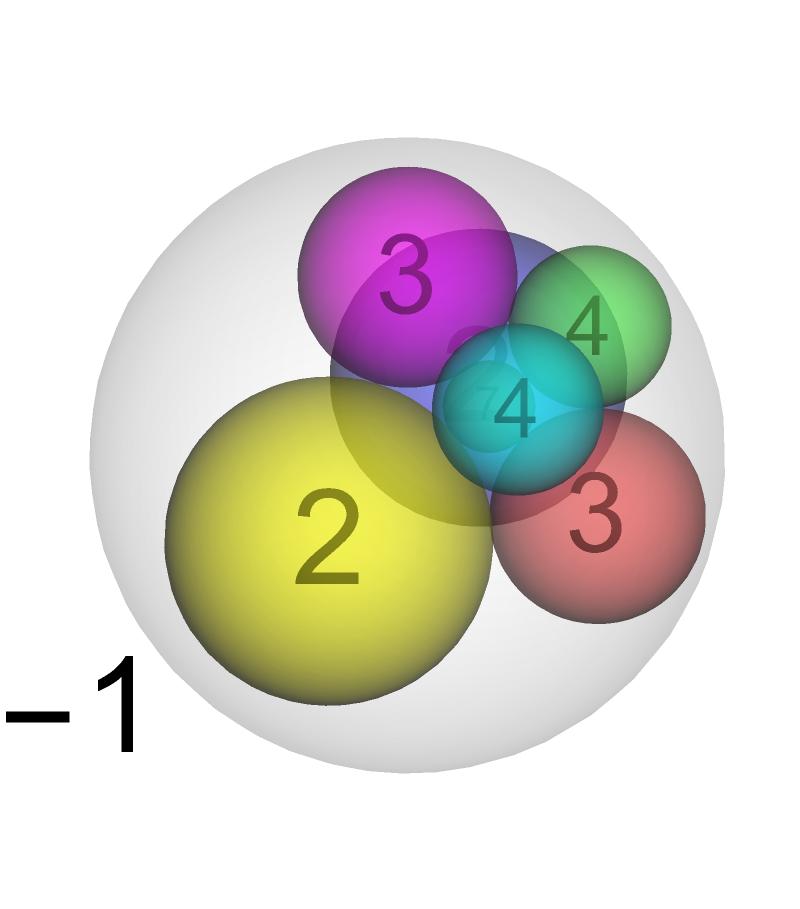}};
		
		\node at (0,0) {	\includegraphics[align=c,width=.4\textwidth]{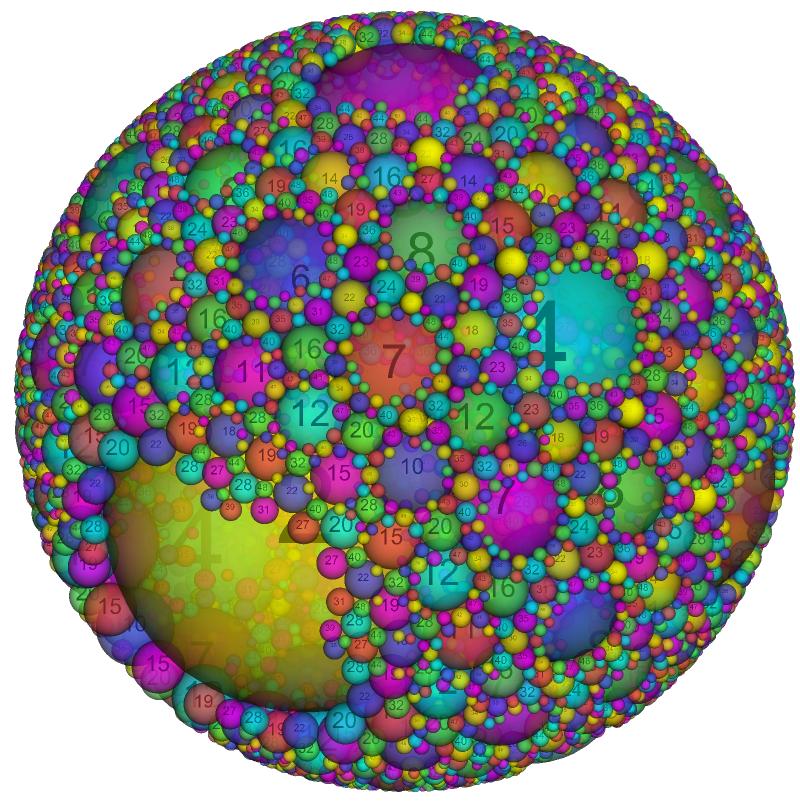}} ;
	\end{tikzpicture}
	\caption{
		An orthoplicial sphere packing satisfying the integrality condition (left) and the corresponding integral crystallographic packing  $\mathscr P_{\{3,3,4\}}(-1,2,2,3)$ (right).
	}
	\label{fig:apo334}
\end{figure}

The integrality condition of Corollary \ref{cor:integrality} states that if  $b_1,b_2,b_3,b_4,\sqrt{\Delta_{\{3,3,4\}}}\in\mathbb Z$ where
\begin{align}\label{eq:intortho}
	\Delta_{\{3,3,4\}}=	(b_1+b_2+b_3+b_4)^2-2(b_1^2+b_2^2+b_3^2+b_4^2)
\end{align}	
then the orthoplicial crystallographic packing $\mathscr P_{\{3,3,4\}}(b_1,b_2,b_3,b_4)$ is integral (see Figure \ref{fig:apo334}). 

\begin{figure}[H]\label{fig:strips334}
	\centering
	
	\begin{tikzpicture}
		
		\begin{scope}
			\begin{scope}[xshift=-5.5cm]
				\node at (0,0) {	\begin{tikzpicture}[scale=1.4] 		
		
		\newcommand\xmin{1}
\newcommand\xmax{2.414}

\newcommand\ymin{1}
\newcommand\ymax{2}
		
		\draw[very thin] (-\xmin,1) -- (\xmax,1);
		\draw[very thin] (0,0.5) circle (.5cm);
		\draw[very thin] (1.,0.5) circle (.5cm) ;
		
		\draw[very thick,red] (-\xmin,0) -- (\xmax,0) node[pos=.95, below] {$S_v$};
		
\draw[blue,thick] (-\xmin,0.5) -- (\xmax,.5)	node[pos=.95, above] {$R_1$};
\draw[thick,blue] (0,0) circle (1cm) node [xshift=-1.3cm,yshift=-1.3cm] {$R_2$};
\draw[blue,thick] (0.5,-\ymin) -- (0.5,\ymax) node[pos=.95,right] {$R_3$};
\draw[blue,thick] (0,-\ymin) -- (0,\ymax) node[pos=.95,left] {$S_f$};

	\end{tikzpicture}};
				
			\end{scope}
			\begin{scope}[xshift=-0cm]
				\node at (0,0) {	\begin{tikzpicture}[scale=1.4] 

		\newcommand\xmin{1}
		\newcommand\xmax{2.414}

		\newcommand\ymin{1}
		\newcommand\ymax{2}

	\draw[very thin] (-\xmin,1) -- (\xmax,1);
	\draw[very thin] (0,0.5) circle (.5cm);
	\draw[very thin] (1.414,0.5) circle (.5cm) ;
	\draw[very thin] (0.707,0.25) circle (.25cm);
	\draw[very thin] (0.707,0.75) circle (.25cm);
	
	\draw[very thick,red] (-\xmin,0) -- (\xmax,0) node[pos=.95, below] {$S_v$};
	
\draw[blue,thick] (-\xmin,0.5) -- (\xmax,.5)	node[pos=.95, above] {$R_1$};
\draw[thick,blue] (0,0) circle (1cm) node [xshift=-1.3cm,yshift=-1.3cm] {$R_2$};
\draw[blue,thick] (0.707,-\ymin) -- (0.707,\ymax) node[pos=.95,right] {$R_3$};
\draw[blue,thick] (0,-\ymin) -- (0,\ymax) node[pos=.95,left] {$S_f$};

	\end{tikzpicture}};
				
			\end{scope}
			\begin{scope}[xshift=5.5cm]
				\node at (0,0) {	\begin{tikzpicture}[scale=1.4]

		\newcommand\xmin{1}
		\newcommand\xmax{2}

		\newcommand\ymin{1}
		\newcommand\ymax{2}

		\clip (-1.1,-1.1) rectangle (2.1,2.1);

	\draw[very thin] (-\xmin,1.414) -- (\xmax,1.414);
\draw[very thin, draw=black] (1.000,1.061) circle (0.3536cm);
\draw[very thin, draw=black] (1.000,0.3536) circle (0.3536cm);
\draw[very thin, draw=black] (0,1.061) circle (0.3536cm);
\draw[very thin, draw=black] (0,0.3536) circle (0.3536cm);
\draw[very thin, draw=black] (0.5000,0.8839) circle (0.1768cm);
\draw[very thin, draw=black] (0.5000,0.5303) circle (0.1768cm);

\draw[very thick,red] (-\xmin,0) -- (\xmax,0)
node[pos=.95, below] {$S_v$};

\draw[thick,blue] (-\xmin,0.707) -- (\xmax,.707)	node[pos=.95, above] {$R_1$};
\draw[thick,blue] (0,0) circle (1cm) node [xshift=-1.3cm,yshift=-1.3cm] {$R_2$};
\draw[thick,blue] (0.5,-\ymin) -- (0.5,\ymax) node[pos=.95,right] {$R_3$};
\draw[thick,blue] (0,-\ymin) -- (0,\ymax) node[pos=.95,left] {$S_f$};

%
		
	\end{tikzpicture}};
				
			\end{scope}
		\end{scope}	
		
		\begin{scope}[scale=2,yshift=-2.7cm] 
			\node {	\includegraphics[align=c,width=0.4\textwidth]{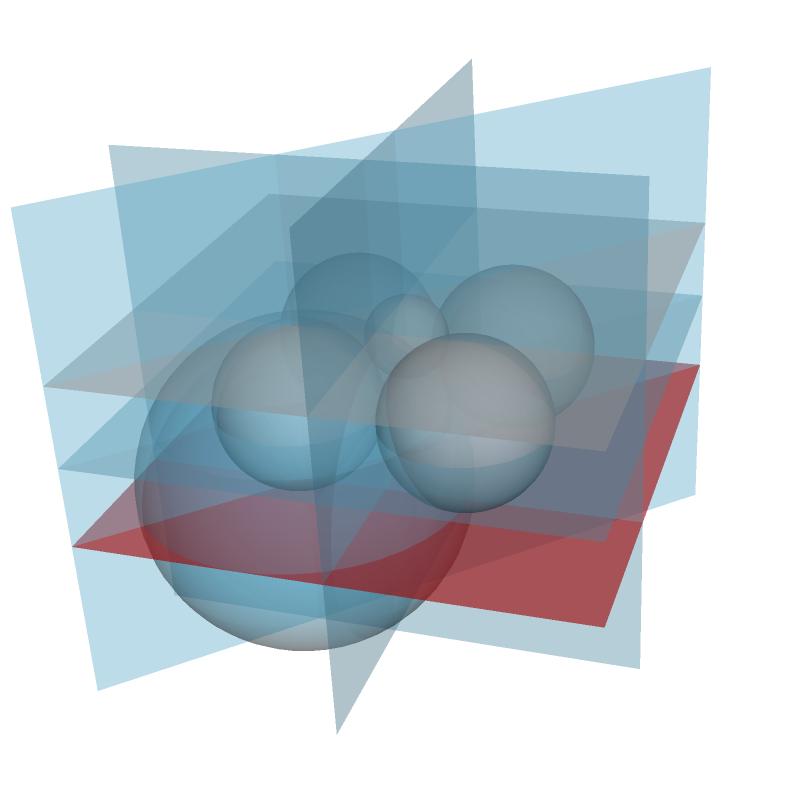} };
			\node at (1.28,0.15) [red] {$S_v$};
			\node at (-1.15,1.09) [blue] {$S_f$};
			\node at (1.3,.4) [blue] {$R_1$};
			\node at (-.55,-1.1) [blue] {$R_2$};
			\node at (.35,1.4) [blue] {$R_3$};
			\node at (1.32,1.4) [blue] {$R_4$};
		\end{scope}	
	\end{tikzpicture}
	
	\vspace{-1cm}
	
	\caption{
		The strip packing with the fundamental symmetries of the tetrahedron (top left), the octahedron (top center), the cube (top right) and the  orthoplex (bottom).
	}
	\label{fig:symaposec334}
\end{figure}

\begin{figure}[H]
	\centering
	
	\begin{tikzpicture}
		\begin{scope}
			\begin{scope}[xshift=-5.5cm]
				\node at (0,0) {\includegraphics[trim=0 0 40 60,clip,align=c,width=0.35\textwidth]{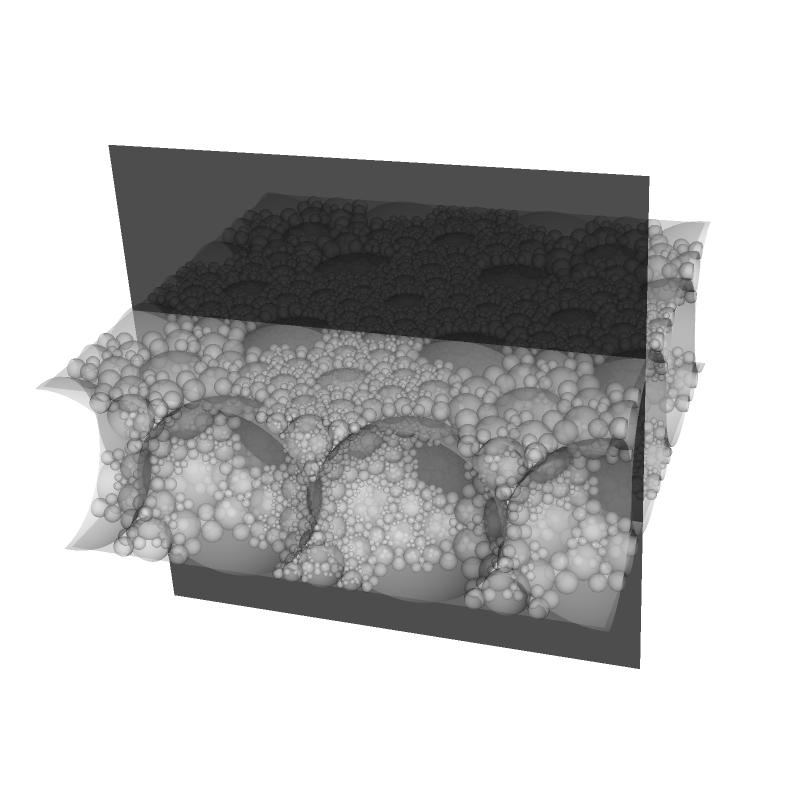}};
				
			\end{scope}
			\begin{scope}[xshift=-0cm]
				\node at (0,0) {\includegraphics[trim=0 0 40 60,clip,align=c,width=0.35\textwidth]{img/sections/aposection33_334b}};
				
			\end{scope}
			\begin{scope}[xshift=5.cm]
				\node at (0,0) {\includegraphics[align=c,width=0.25\textwidth]{img/strips/33packing.pdf}};
				
			\end{scope}
		\end{scope}

	\end{tikzpicture}
	
	\vspace{-1cm}
	
	\caption{
(From left to right) The orthoplicial crystallographic packing $\mathscr{P}_{\{3,3,4\}}$ with a cutting sphere $\Sigma_{\{3,3,4\}}^{\{3,3\}}$, the tetrahedral Apollonian section $\mathscr S_{\{3,3,4\}}^{\{3,3\}}$ with $\Sigma_{\{3,3,4\}}^{\{3,3\}}$, and the tetrahedral crystallographic packing $\mathscr{P}_{\{3,3\}}$. 		
	}
	\label{fig:sectionsgeo334_33}
\end{figure}

\begin{figure}[H]
	\centering
	
	\begin{tikzpicture}
		\begin{scope}
			\begin{scope}[xshift=-5.5cm]
				\node at (0,0) {\includegraphics[trim=0 0 40 20,clip,align=c,width=0.35\textwidth]{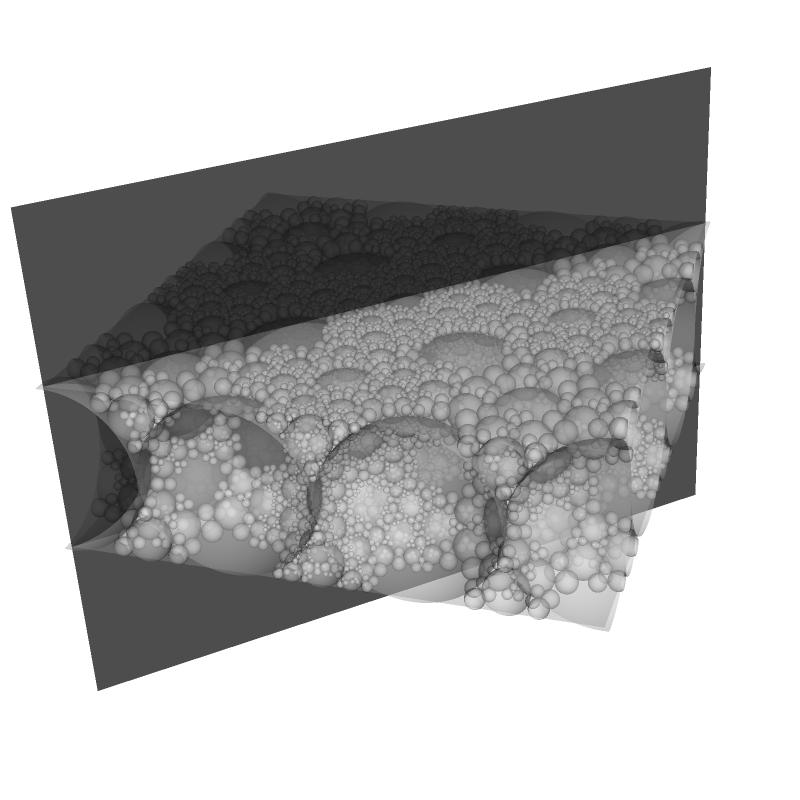}};
				
			\end{scope}
			\begin{scope}[xshift=-0cm]
				\node at (0,0) {\includegraphics[trim=0 0 40 20,clip,align=c,width=0.35\textwidth]{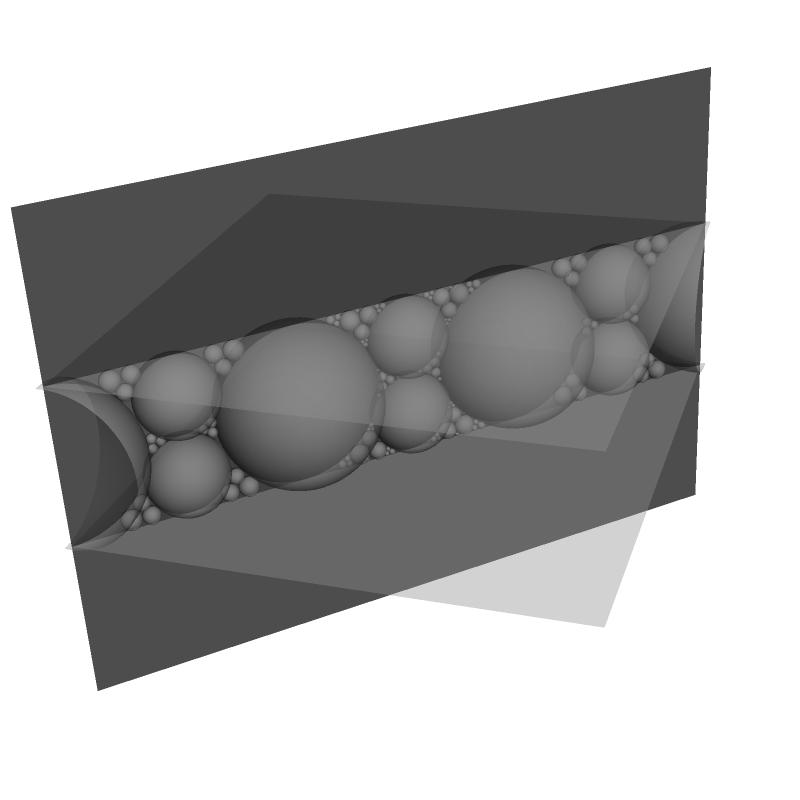}};
				
			\end{scope}
			\begin{scope}[xshift=5.cm]
				\node at (0,0) {\includegraphics[align=c,width=0.25\textwidth]{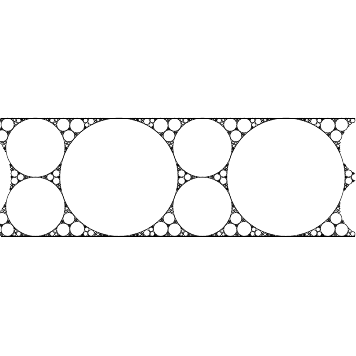}};
				
			\end{scope}
		\end{scope}

	\end{tikzpicture}
	
	\vspace{-1cm}
	
	\caption{
(From left to right) The orthoplicial crystallographic packing $\mathscr{P}_{\{3,3,4\}}$ with a cutting sphere $\Sigma_{\{3,3,4\}}^{\{3,4\}}$,
the octahedral Apollonian section $\mathscr S_{\{3,3,4\}}^{\{3,4\}}$ with $\Sigma_{\{3,3,4\}}^{\{3,4\}}$, and the octahedral crystallographic packing $\mathscr{P}_{\{3,4\}}$. 		
	}
	\label{fig:sectionsgeo334_34}
\end{figure}

\begin{figure}[H]
	\centering
	
	\begin{tikzpicture}
		\begin{scope}
			\begin{scope}[xshift=-5.5cm]
				\node at (0,0) {\includegraphics[trim=0 0 30 0,clip,align=c,width=0.35\textwidth]{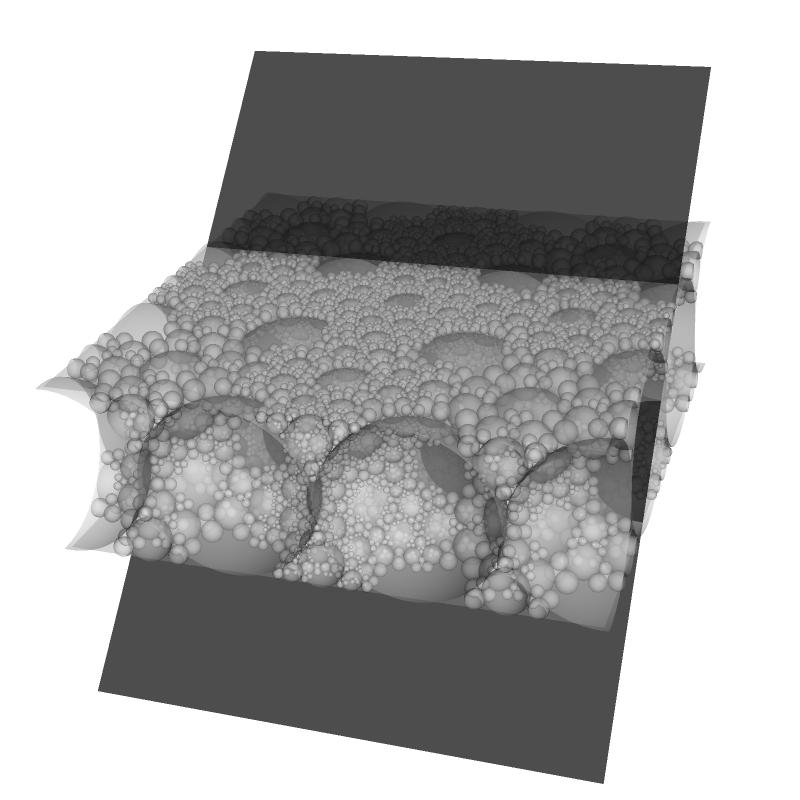}};
			\end{scope}
			\begin{scope}[xshift=-0cm]
				\node at (0,0) {\includegraphics[trim=0 0 40 0,clip,align=c,width=0.35\textwidth]{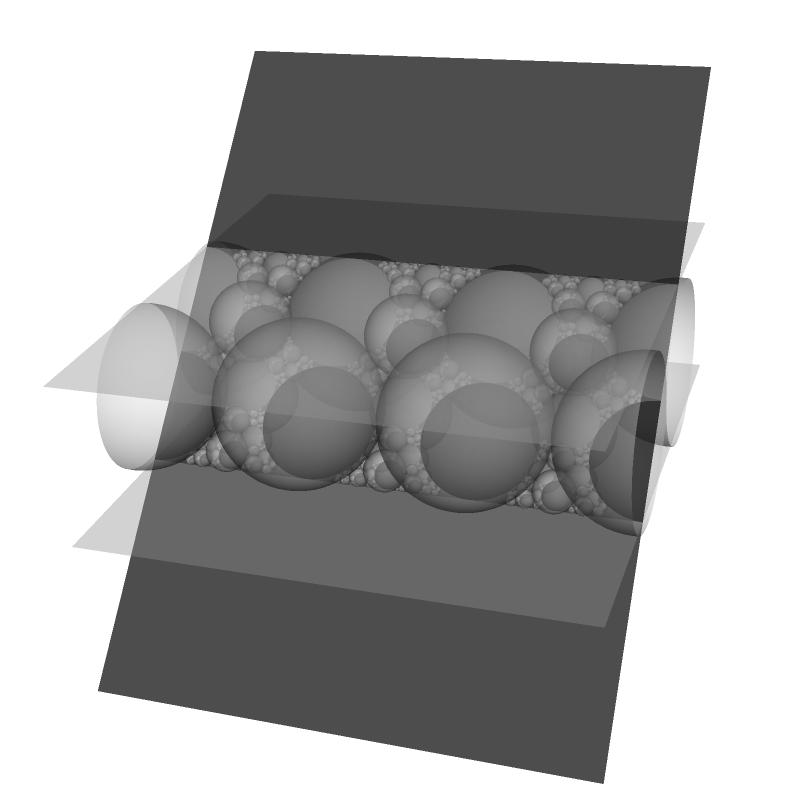}};
				
			\end{scope}
			\begin{scope}[xshift=5.cm]
				\node at (0,0) {\includegraphics[align=c,width=0.25\textwidth]{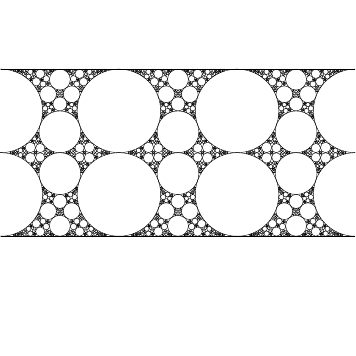}};
			\end{scope}
		\end{scope}

	\end{tikzpicture}

	\caption{
		(From left to right) The orthoplicial crystallographic packing $\mathscr{P}_{\{3,3,4\}}$ with a cutting sphere $\Sigma_{\{3,3,4\}}^{\{4,3\}}$, the cubic Apollonian section $\mathscr S_{\{3,3,4\}}^{\{4,3\}}$ with $\Sigma_{\{3,3,4\}}^{\{4,3\}}$, and the cubic crystallographic packing $\mathscr{P}_{\{4,3\}}$. 		
	}
	\label{fig:sectionsgeo334_43}
\end{figure}

Every orthoplicial crystallographic packing contains a tetrahedral Apollonian section $\mathscr{S}^{\{3,3\}}_{\{3,3,4\}}$, an octahedral Apollonian section $\mathscr{S}^{\{3,4\}}_{\{3,3,4\}}$ and a cubic Apollonian section $\mathscr{S}^{\{4,3\}}_{\{3,3,4\}}$.  In the strip packing, the corresponding cutting spheres are  $\Sigma_{\{3,3,4\}}^{\{3,3\}}=S_f$, $\Sigma_{\{3,3,4\}}^{\{3,4\}}=R_4$ and $\Sigma_{\{3,3,4\}}^{\{4,3\}}$ is the plane passing through the origin and normal vector $(0,-1,1)$ (see Figures \ref{fig:symaposec334}, \ref{fig:sectionsgeo334_33}, \ref{fig:sectionsgeo334_34}, \ref{fig:sectionsgeo334_43}). The homomorphisms between the full symmetry groups are
\begin{align}\label{eq:hom334}
	\begin{tabular}{cccccc}
		&$\Gamma_{\{3,3\}}$&$\xrightarrow{\phi_{\{3,3,4\}}^{\{3,3\}}}$&$\Gamma_{\{3,3,4\}}$\\
		&$r_1$&$\longmapsto$& $r_1$\\
		&$r_2$&$\longmapsto$& $r_2$\\
		&$r_3$&$\longmapsto$& $r_3$\\
		&$s_f$&$\longmapsto$& $(r_4s_f)^2$\\
	\end{tabular}&&
	\begin{tabular}{cccccc}
		&$\Gamma_{\{3,4\}}$&$\xrightarrow{\phi_{\{3,3,4\}}^{\{3,4\}}}$&$\Gamma_{\{3,3,4\}}$\\
		&$r_1$&$\longmapsto$& $r_1$\\
		&$r_2$&$\longmapsto$& $r_2$\\
		&$r_3$&$\longmapsto$& $r_3r_4r_3$\\
		&$s_f$&$\longmapsto$& $(r_4s_f)^2$\\
	\end{tabular}&&
	\begin{tabular}{cccccc}
		&$\Gamma_{\{4,3\}}$&$\xrightarrow{\phi_{\{3,3,4\}}^{\{4,3\}}}$&$\Gamma_{\{3,3,4\}}$\\
		&$r_1$&$\longmapsto$& $r_1r_4r_3r_4$\\
		&$r_2$&$\longmapsto$& $r_2$\\
		&$r_3$&$\longmapsto$& $r_3$\\
		&$s_f$&$\longmapsto$& $r_4s_fr_4$\\
	\end{tabular}
\end{align}


		\subsection{Hypercube $\{4,3,3\}$}\label{sec:433}
\begin{figure}[H]
			\includegraphics[align=c,width=.22\textwidth]{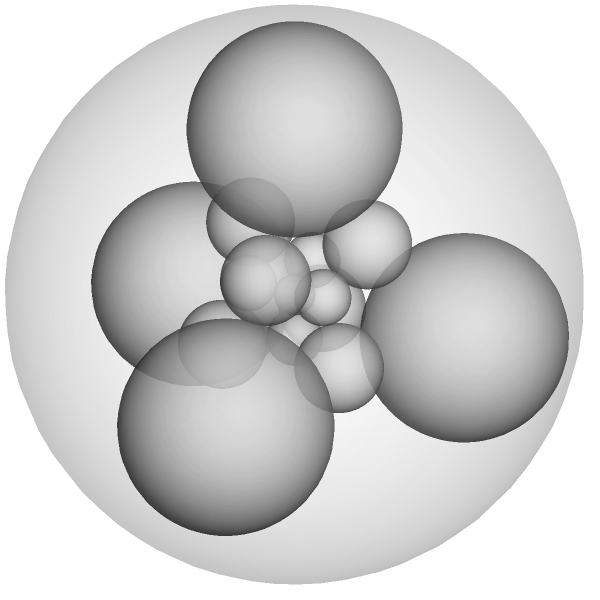} \hspace{.1cm}  
			\includegraphics[align=c,width=.22\textwidth]{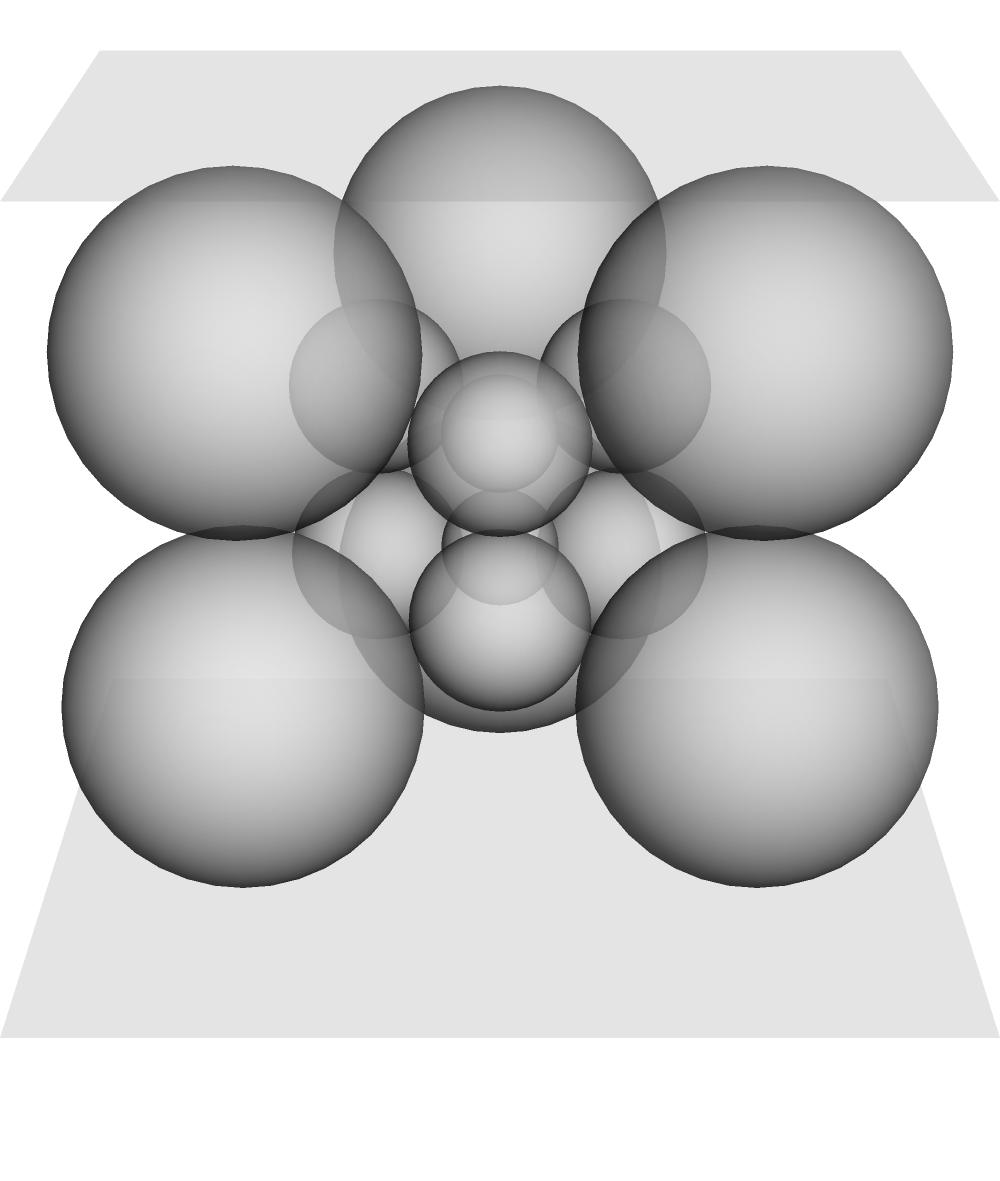} \hspace{.1cm} 
			\includegraphics[align=c,width=.22\textwidth]{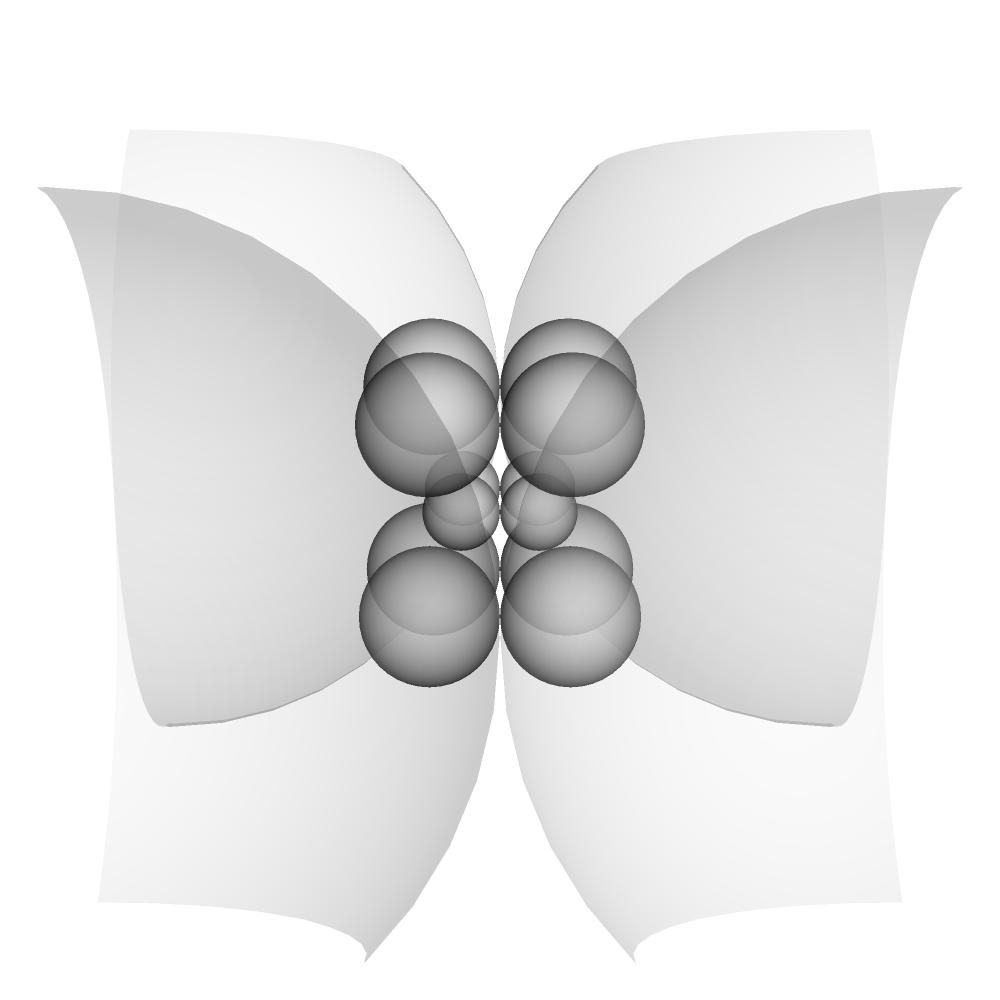} \hspace{.1cm} 
			\includegraphics[align=c,width=.22\textwidth]{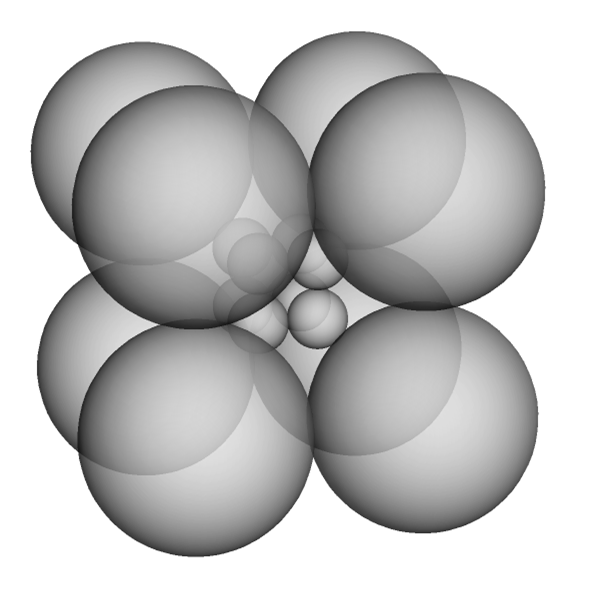}
			\caption{(From left to right) Arrangement projections of a vertex-centered, edge-centered, ridge-centered and facet-centered canonical hypercube.}
			\label{fig:centeredhcub}
		\end{figure}

		In Figure \ref{fig:centeredhcub}, we show four polytopal sphere packings obtained by the arrangement projections of face-centered canonical hypercubes.
	The full symmetry group $\Gamma_{\{4,3,3\}}$ is isomorphic to $\langle \mathbf{R}_1,\mathbf R_2,\mathbf{R}_3,\mathbf R_4,\mathbf S_f\rangle <\mathrm{SL}_5(\mathbb Z)$ where
\begin{equation}
\begin{gathered}
	\mathbf{R}_1=\left(
	\begin{array}{ccccc}
		& 1 &  &  &  \\
		1 &  &  &  &  \\
		1 & -1 & 1 &  &  \\
		1 & -1 &  & 1 &  \\
		1 & -1 &  &  & 1 \\
	\end{array}
	\right)\quad
	\mathbf{R}_2= \left(
	\begin{array}{ccccc}
		1 &  &  &  &  \\
		1 & -1 & 1 &  &  \\
		&  & 1 &  &  \\
		&  &  & 1 &  \\
		&  &  &  & 1 \\
	\end{array}
	\right)\quad
	\mathbf{R}_3= \left(
	\begin{array}{ccccc}
		1 &  &  &  &  \\
		& 1 &  &  &  \\
		& 1 & -1 & 1 &  \\
		&  &  & 1 &  \\
		&  &  &  & 1 \\
	\end{array}
	\right)\\[.2cm]
	\mathbf{R}_4= \left(
	\begin{array}{ccccc}
		1 &  &  &  &  \\
		& 1 &  &  &  \\
		&  & 1 &  &  \\
		&  & 1 & -1 & 1 \\
		&  &  &  & 1 \\
	\end{array}
	\right)\quad
	\mathbf S_f=\left(
	\begin{array}{ccccc}
		1 &  &  &  &  \\
		& 1 &  &  &  \\
		&  & 1 &  &  \\
		&  &  & 1 &  \\
		1 &  &  & 3 & -1 \\
	\end{array}
	\right)
\end{gathered}
\end{equation}	

Any fundamental bend vector $\mathbf b=(b_1,b_2,b_3,b_4,b_5)^T$ of a hypercubic sphere packing 
satisfies the quadratic equation $\mathbf{b}^T\mathbf Q_{\{4,3,3\}}\mathbf{b}=0$ for the bisymmetric matrix
\begin{align}\label{eq:descartes334}
	\mathbf Q_{\{4,3,3\}}=
 	\left(
	\begin{array}{ccccc}
		2 & -3 & 0 & 0 & -1 \\
		\ast & 6 & -3 & 0 & \ast \\
		\ast & \ast & 6 & \ast & \ast \\
		\ast & \ast & \ast & \ast & \ast \\
		\ast & \ast & \ast & \ast & \ast \\
	\end{array}
	\right)
\end{align} 

The latter  is equivalent to the following quadratic equation
\begin{align}
	(b_1-b_2)^2+(b_2-b_3)^2+(b_3-b_4)^2+(b_4-b_5)^2=\frac13(b_1+b_5)^2.
\end{align}

\begin{figure}[H]
	\centering
	\begin{tikzpicture}[scale=2.5] 
		\node at (-3.2,0) {	\includegraphics[clip,trim=0 60 0 60,align=c,width=.49\textwidth]{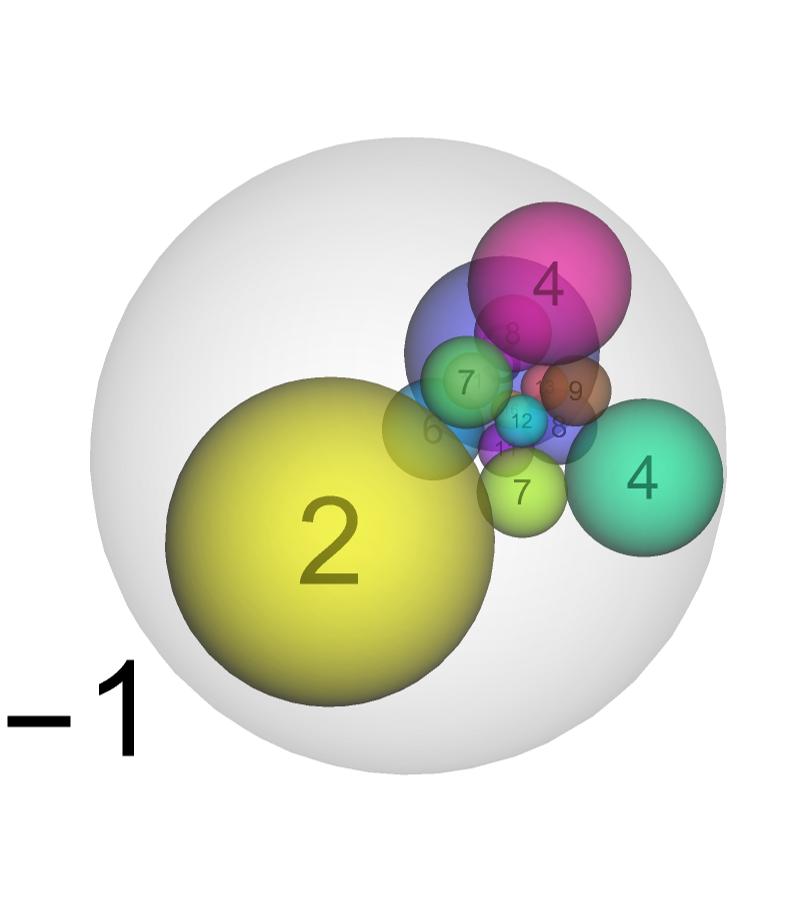}};
		
		\node at (0,0) {	\includegraphics[align=c,width=.4\textwidth]{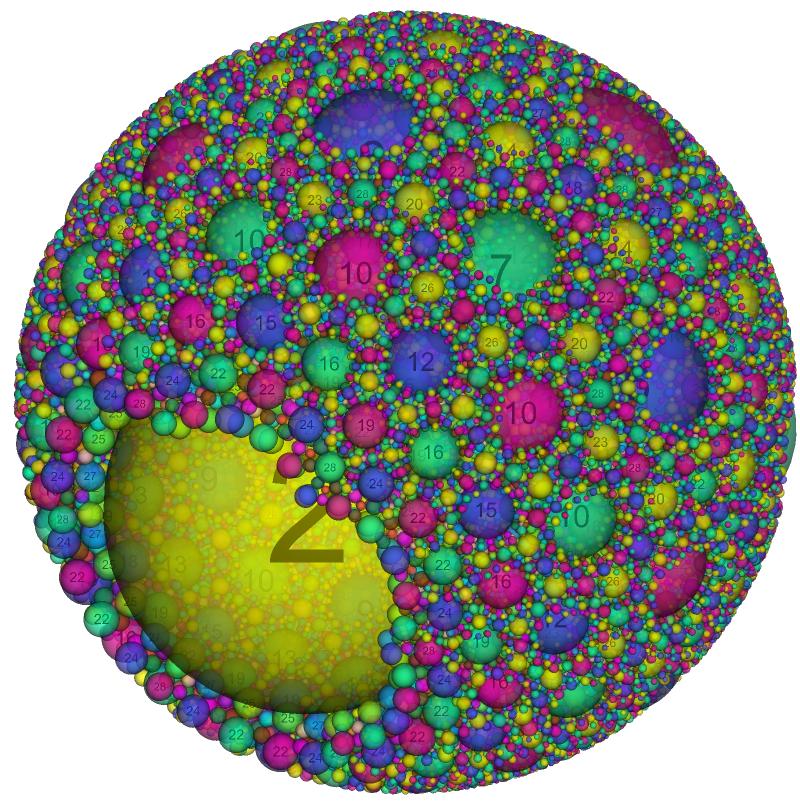}} ;
	\end{tikzpicture}
	\vspace{-.3cm}
	\caption{
		An hypercubic sphere packing satisfying the integrality condition (left) and the corresponding integral crystallographic packing  $\mathscr P_{\{4,3,3\}}(-1,2,6,11)$ (right).
	}
	\label{fig:apo433}
\end{figure}

The integrality condition of Corollary \ref{cor:integrality} states that if  $b_1,b_2,b_3,b_4,\sqrt{\Delta_{\{4,3,3\}}}\in\mathbb Z$ where
\begin{align}\label{eq:intcube}
\Delta_{\{4,3,3\}}=3\left((b_1+b_4)^2-2(b_1-b_2)^2-2(b_2-b_3)^2-2(b_3-b_4)^2)\right)
\end{align}
then the hypercubic crystallographic packing $\mathscr P_{\{4,3,3\}}(b_1,b_2,b_3,b_4)$ is integral (see Figure \ref{fig:apo433}). 

\begin{figure}[H]
			\centering
			\begin{tabular}{cc}
				&\begin{tikzpicture}[scale=2.] 
					\node at (-3,0) {		\begin{tikzpicture}[scale=1.4]

		\newcommand\xmin{1}
		\newcommand\xmax{2}

		\newcommand\ymin{1}
		\newcommand\ymax{2}

		\clip (-1.1,-1.1) rectangle (2.1,2.1);

	\draw[very thin] (-\xmin,1.414) -- (\xmax,1.414);
\draw[very thin, draw=black] (1.000,1.061) circle (0.3536cm);
\draw[very thin, draw=black] (1.000,0.3536) circle (0.3536cm);
\draw[very thin, draw=black] (0,1.061) circle (0.3536cm);
\draw[very thin, draw=black] (0,0.3536) circle (0.3536cm);
\draw[very thin, draw=black] (0.5000,0.8839) circle (0.1768cm);
\draw[very thin, draw=black] (0.5000,0.5303) circle (0.1768cm);

\draw[very thick,red] (-\xmin,0) -- (\xmax,0)
node[pos=.95, below] {$S_v$};

\draw[thick,blue] (-\xmin,0.707) -- (\xmax,.707)	node[pos=.95, above] {$R_1$};
\draw[thick,blue] (0,0) circle (1cm) node [xshift=-1.3cm,yshift=-1.3cm] {$R_2$};
\draw[thick,blue] (0.5,-\ymin) -- (0.5,\ymax) node[pos=.95,right] {$R_3$};
\draw[thick,blue] (0,-\ymin) -- (0,\ymax) node[pos=.95,left] {$S_f$};

%
		
	\end{tikzpicture}};
					
					\node {	\includegraphics[clip,trim=0 0 0 20,align=c,width=0.4\textwidth]{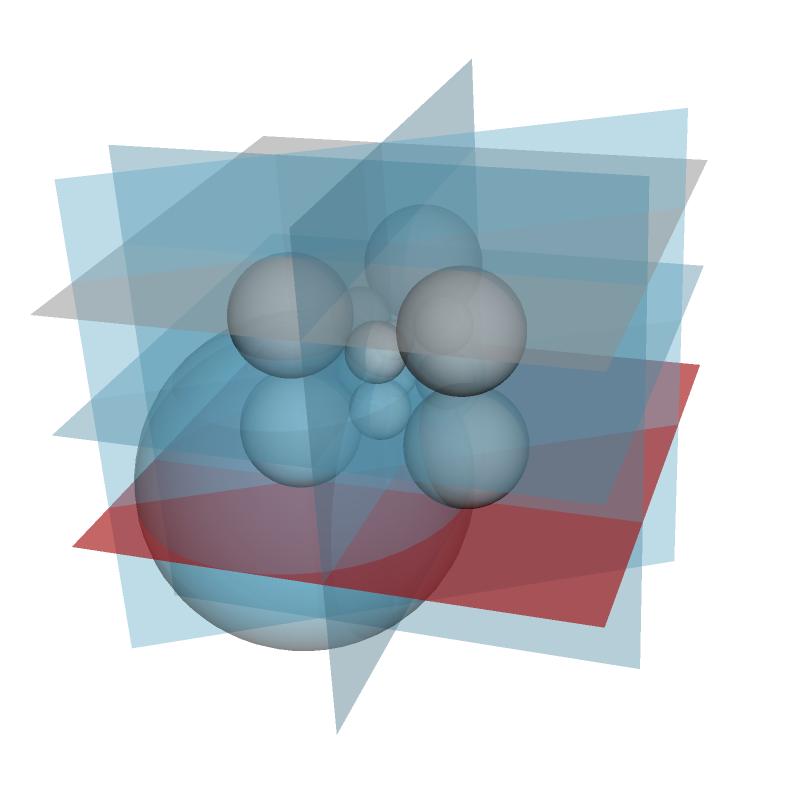} };
					\node at (1.25,0.15) [red] {$S_v$};
					\node at (-1.15,1.09) [blue] {$S_f$};
					\node at (1.24,.5) [blue] {$R_1$};
			 		\node at (-.55,-1.1) [blue] {$R_2$};					
					\node at (.35,1.4) [blue] {$R_3$};
					\node at (1.2,1.2) [blue] {$R_4$};
					
				\end{tikzpicture}
				
			\end{tabular}
			
			\vspace{-1.2cm}
			\caption{
				The strip packing with the fundamental symmetries of the cube (left) and the hypercube (right).
			}
			\label{fig:symaposec433}
		\end{figure}

\begin{figure}[H]
	\centering
	
	\begin{tikzpicture}
		\begin{scope}
			\begin{scope}[xshift=-5.5cm]
				\node at (0,0) {\includegraphics[trim=0 0 40 60,clip,align=c,width=0.35\textwidth]{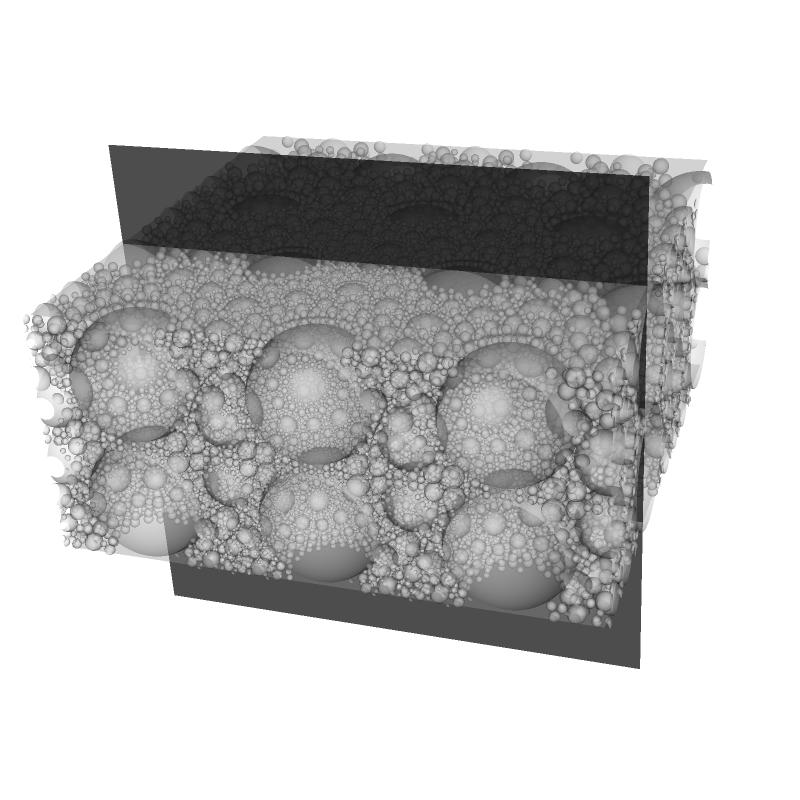}};
				
			\end{scope}
			\begin{scope}[xshift=-0cm]
				\node at (0,0) {\includegraphics[trim=0 0 40 60,clip,align=c,width=0.35\textwidth]{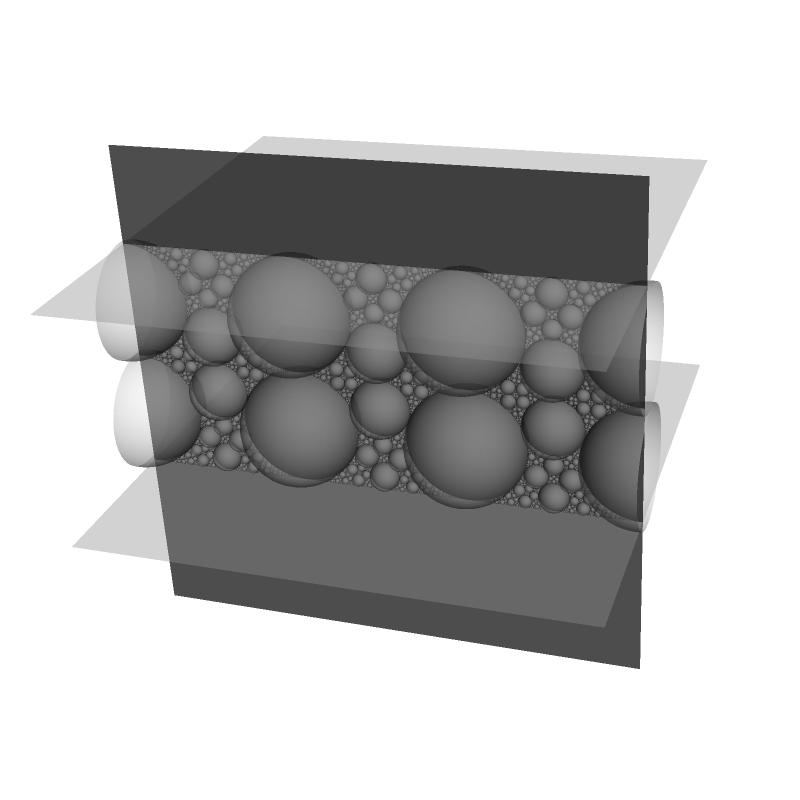}};
				
			\end{scope}
			\begin{scope}[xshift=5.cm]
				\node at (0,0) {\includegraphics[align=c,width=0.25\textwidth]{img/strips/43packing.pdf}};
			\end{scope}
		\end{scope}

	\end{tikzpicture}
	
	\vspace{-1.cm}
	
	\caption{ 
(From left to right) The hypercubic crystallographic packing $\mathscr{P}_{\{4,3,3\}}$ with a cutting sphere $\Sigma_{\{4,3,3\}}^{\{4,3\}}$, the cubic Apollonian section $\mathscr S_{\{4,3,3\}}^{\{4,3\}}$ with $\Sigma_{\{4,3,3\}}^{\{4,3\}}$, and the cubic crystallographic packing $\mathscr{P}_{\{4,3\}}$. 		
	}
	\label{fig:sectionsgeo433}
\end{figure}

Every hypercubic crystallographic packing contains a geometric cubic Apollonian section  $\mathscr{S}^{\{4,3\}}_{\{4,3,3\}}$. In the strip packing, the cutting sphere is  $\Sigma_{\{4,3,3\}}^{\{4,3\}}=S_f$	(see Figures \ref{fig:symaposec433}, \ref{fig:sectionsgeo433}). The homomorphism between the full symmetry groups is described below.
		\begin{align}\label{eq:hom433}
			\begin{tabular}{cccccc}
				&$\Gamma_{\{4,3\}}$&$\xrightarrow{\phi_{\{4,3,3\}}^{\{4,3\}}}$&$\Gamma_{\{4,3,3\}}$\\
				&$r_1$&$\longmapsto$& $r_1$\\
				&$r_2$&$\longmapsto$& $r_2$\\
				&$r_3$&$\longmapsto$& $r_3$\\
				&$s_f$&$\longmapsto$& $(r_4s_f)^3$\\
			\end{tabular}
		\end{align}
			\vspace{-.5cm}
			\vspace{-.5cm}

\subsection{24-cell $\{3,4,3\}$}\label{sec:343}

\begin{figure}[H]
			\includegraphics[align=c,width=.22\textwidth]{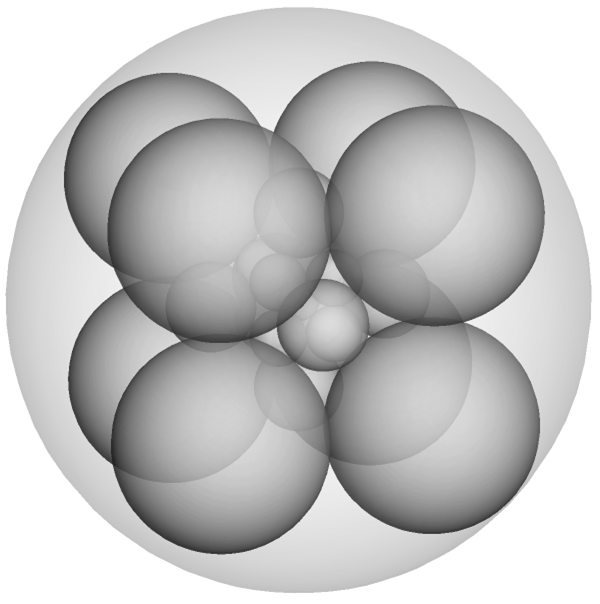} \hspace{.1cm}  
			\includegraphics[align=c,width=.22\textwidth]{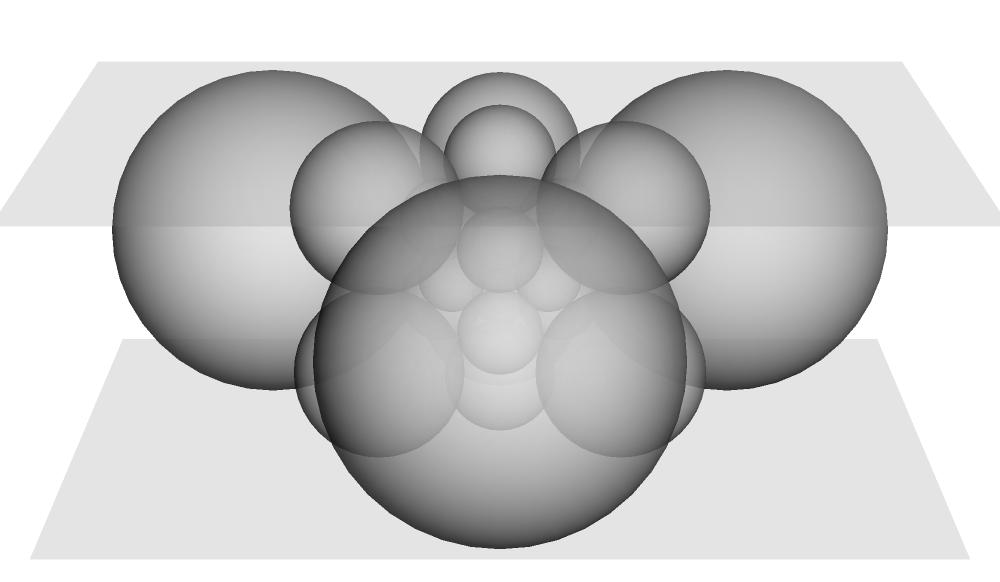} \hspace{.1cm} 
			\includegraphics[align=c,width=.22\textwidth]{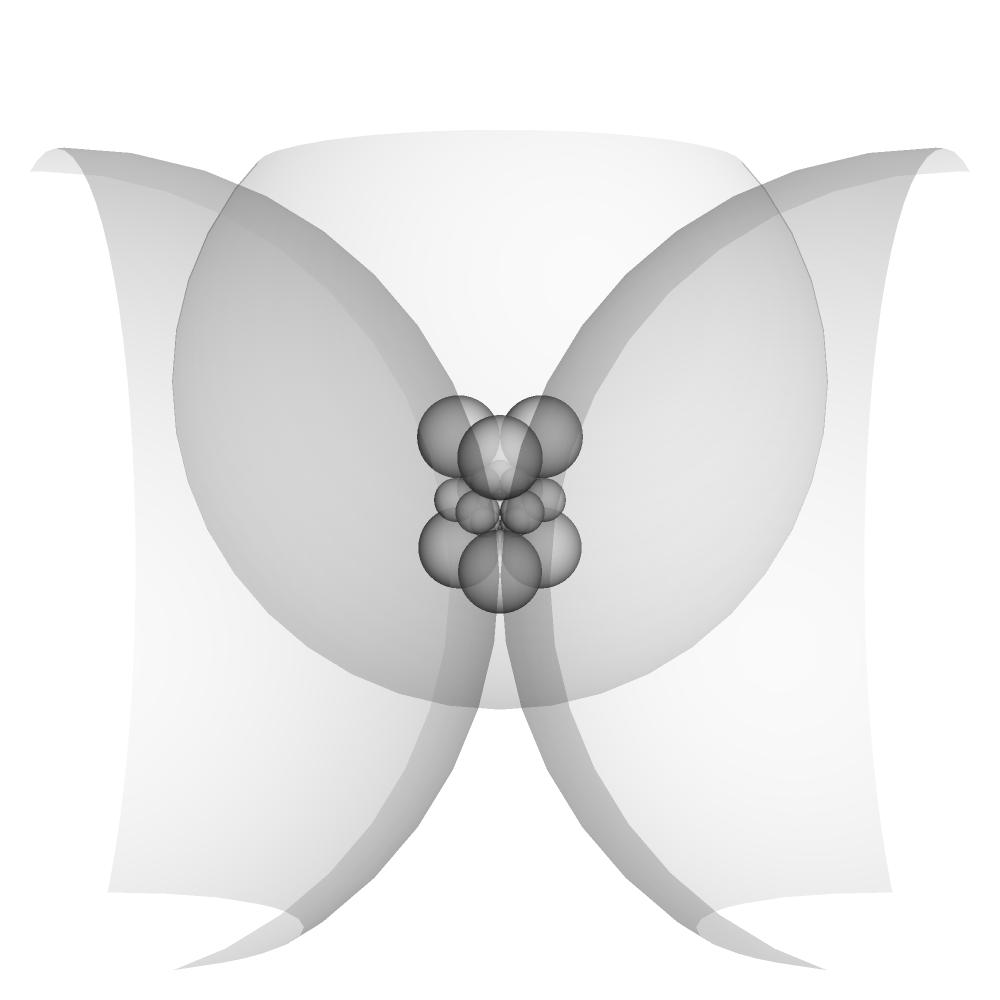} \hspace{.1cm} 
			\includegraphics[align=c,width=.22\textwidth]{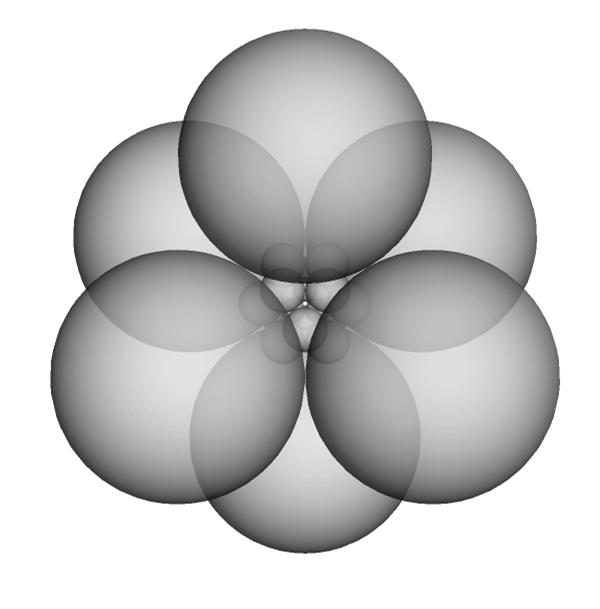}
			\vspace{-.2cm}
			\caption{(From left to right) Arrangement projections of a vertex-centered, edge-centered, ridge-centered and facet-centered canonical 24-cell.}
			\label{fig:centered24}
		\end{figure}

		In Figure \ref{fig:centered24}, we show four polytopal sphere packings obtained by the arrangement projections of face-centered canonical 24-cells.
		\vspace{-.2cm}
			The full symmetry group $\Gamma_{\{3,4,3\}}$ is isomorphic to $\langle \mathbf{R}_1,\mathbf R_2,\mathbf{R}_3,\mathbf R_4,\mathbf S_f\rangle <\mathrm{SL}_5(\mathbb Z)$ where
\begin{equation}
	\begin{gathered}
	\mathbf{R}_1= \left(
	\begin{array}{ccccc}
		& 1 &  &  &  \\
		1 &  &  &  &  \\
		&  & 1 &  &  \\
		1 & -1 &  & 1 &  \\
		1 & -1 &  &  & 1 \\
	\end{array}
	\right)\quad
	\mathbf{R}_2= \left(
	\begin{array}{ccccc}
		1 &  &  &  &  \\
		&  & 1 &  &  \\
		& 1 &  &  &  \\
		&  &  & 1 &  \\
		& 1 & -1 &  & 1 \\
	\end{array}
	\right)\quad
	\mathbf{R}_3= \left(
	\begin{array}{ccccc}
		1 &  &  &  &  \\
		& 1 &  &  &  \\
		1 &  & -1 & 1 &  \\
		&  &  & 1 &  \\
		1 &  & -2 & 1 & 1 \\
	\end{array}
	\right)\\[.2cm]
	\mathbf{R}_4= \left(
	\begin{array}{ccccc}
		1 &  &  &  &  \\
		& 1 &  &  &  \\
		&  & 1 &  &  \\
		& 1 &  & -1 & 1 \\
		&  &  &  & 1 \\
	\end{array}
	\right)\quad
	\mathbf S_f=
	\left(
	\begin{array}{ccccc}
		1 &  &  &  &  \\
		& 1 &  &  &  \\
		&  & 1 &  &  \\
		&  &  & 1 &  \\
		1 &  & 2 & 3 & -1 \\
	\end{array}
	\right)
\end{gathered}
\end{equation}			

Any fundamental bend vector $\mathbf b=(b_1,b_2,b_3,b_4,b_5)^T$ of a 24-cell sphere packing 
satisfies the quadratic equation $\mathbf{b}^T\mathbf Q_{\{3,4,3\}}\mathbf{b}=0$ for the bisymmetric matrix
\begin{align}\label{eq:descartes334}
	\mathbf Q_{\{3,4,3\}}=
	\left(
	\begin{array}{ccccc}
		2 & -3 & -2 & 0 & -1 \\
		\ast & 6 & 0 & -3 & \ast \\
		\ast & \ast & 8 & \ast & \ast \\
		\ast & \ast & \ast & \ast & \ast \\
		\ast & \ast & \ast & \ast & \ast \\
	\end{array}
	\right)
\end{align}

\begin{figure}[H]
	\centering
	\begin{tikzpicture}[scale=2.5] 
		\node at (-3.2,0) {	\includegraphics[clip,trim=0 60 0 60,align=c,width=.49\textwidth]{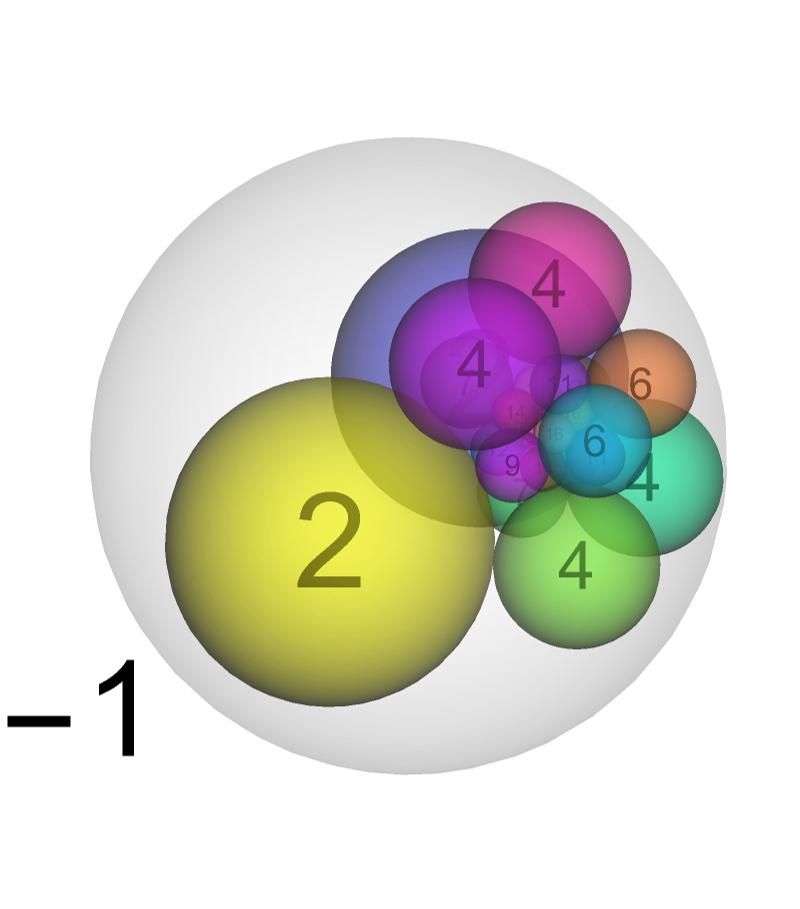}};
		
		\node at (0,0) {	\includegraphics[align=c,width=.4\textwidth]{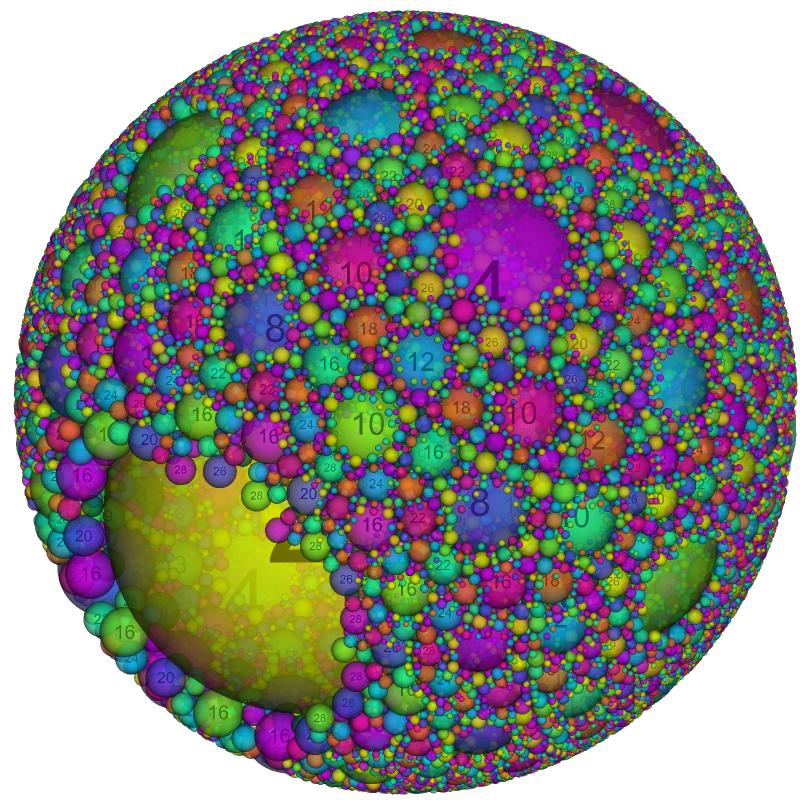}} ;
		
	\end{tikzpicture}
	\caption{
		A 24-cell sphere packing satisfying the integrality condition (left) and the corresponding integral crystallographic packing  $\mathscr P_{\{3,4,3\}}(-1,2,2,7)$ (right).
	}
	
	\label{fig:apo343}
\end{figure}

The integrality condition of Corollary \ref{cor:integrality} states that if  $b_1,b_2,b_3,b_4,\sqrt{\Delta_{\{3,4,3\}}}\in\mathbb Z$ where
\begin{align}\label{eq:int24cell}
	\Delta_{\{3,4,3\}}=3\left(2(b_1+b_4)^2-(b_1-b_4)^2-(b_1-2b_2+b_4)^2-(b_1-2b_3+b_4)^2)\right)
\end{align}
then the 24-cell crystallographic packing $\mathscr P_{\{3,4,3\}}(b_1,b_2,b_3,b_4)$ is integral (see Figure \ref{fig:apo343}). 	

\begin{figure}[H]
		\centering
		\begin{tikzpicture}

				\begin{scope}[xshift=-5.cm]
					\node at (0,0) {	\begin{tikzpicture}[scale=1.4] 

		\newcommand\xmin{1}
		\newcommand\xmax{2.414}

		\newcommand\ymin{1}
		\newcommand\ymax{2}

	\draw[very thin] (-\xmin,1) -- (\xmax,1);
	\draw[very thin] (0,0.5) circle (.5cm);
	\draw[very thin] (1.414,0.5) circle (.5cm) ;
	\draw[very thin] (0.707,0.25) circle (.25cm);
	\draw[very thin] (0.707,0.75) circle (.25cm);
	
	\draw[very thick,red] (-\xmin,0) -- (\xmax,0) node[pos=.95, below] {$S_v$};
	
\draw[blue,thick] (-\xmin,0.5) -- (\xmax,.5)	node[pos=.95, above] {$R_1$};
\draw[thick,blue] (0,0) circle (1cm) node [xshift=-1.3cm,yshift=-1.3cm] {$R_2$};
\draw[blue,thick] (0.707,-\ymin) -- (0.707,\ymax) node[pos=.95,right] {$R_3$};
\draw[blue,thick] (0,-\ymin) -- (0,\ymax) node[pos=.95,left] {$S_f$};

	\end{tikzpicture}};
				\end{scope}
				\begin{scope}[xshift=0cm]
					\node at (0,0) {	\begin{tikzpicture}[scale=1.4]

		\newcommand\xmin{1}
		\newcommand\xmax{2}

		\newcommand\ymin{1}
		\newcommand\ymax{2}

		\clip (-1.1,-1.1) rectangle (2.1,2.1);

	\draw[very thin] (-\xmin,1.414) -- (\xmax,1.414);
\draw[very thin, draw=black] (1.000,1.061) circle (0.3536cm);
\draw[very thin, draw=black] (1.000,0.3536) circle (0.3536cm);
\draw[very thin, draw=black] (0,1.061) circle (0.3536cm);
\draw[very thin, draw=black] (0,0.3536) circle (0.3536cm);
\draw[very thin, draw=black] (0.5000,0.8839) circle (0.1768cm);
\draw[very thin, draw=black] (0.5000,0.5303) circle (0.1768cm);

\draw[very thick,red] (-\xmin,0) -- (\xmax,0)
node[pos=.95, below] {$S_v$};

\draw[thick,blue] (-\xmin,0.707) -- (\xmax,.707)	node[pos=.95, above] {$R_1$};
\draw[thick,blue] (0,0) circle (1cm) node [xshift=-1.3cm,yshift=-1.3cm] {$R_2$};
\draw[thick,blue] (0.5,-\ymin) -- (0.5,\ymax) node[pos=.95,right] {$R_3$};
\draw[thick,blue] (0,-\ymin) -- (0,\ymax) node[pos=.95,left] {$S_f$};

%
		
	\end{tikzpicture}};
				\end{scope}
			
			\begin{scope}[scale=2,xshift=2.7cm] 
				\node {	\includegraphics[clip,trim=0 20 0 20,align=c,width=0.4\textwidth]{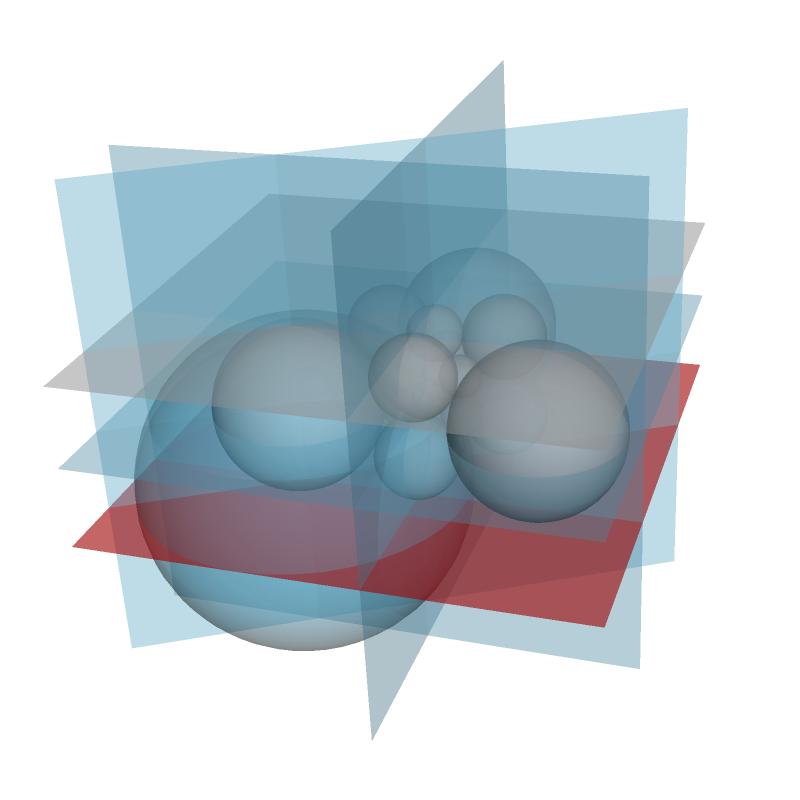} };
				\node at (1.3,0) [red] {$S_v$};
				\node at (-1.15,1.09) [blue] {$S_f$};
				\node at (1.3,.4) [blue] {$R_1$};
				\node at (-.55,-1.1) [blue] {$R_2$};
				\node at (.35,1.4) [blue] {$R_3$};
				\node at (1.2,1.2) [blue] {$R_4$};
			\end{scope}
		\end{tikzpicture}
		\vspace{-1cm}
		\caption{
			The strip packing with the walls of the fundamental symmetries of the octahedron (left), the cube (center) and the 24-cell (right).
		}
		\label{fig:symaposec343}
	\end{figure}

\begin{figure}[H]
		\centering
		
		\begin{tikzpicture}
			\begin{scope}
				\begin{scope}[xshift=-5.5cm]
					\node at (0,0) {\includegraphics[trim=0 0 40 60,clip,align=c,width=0.35\textwidth]{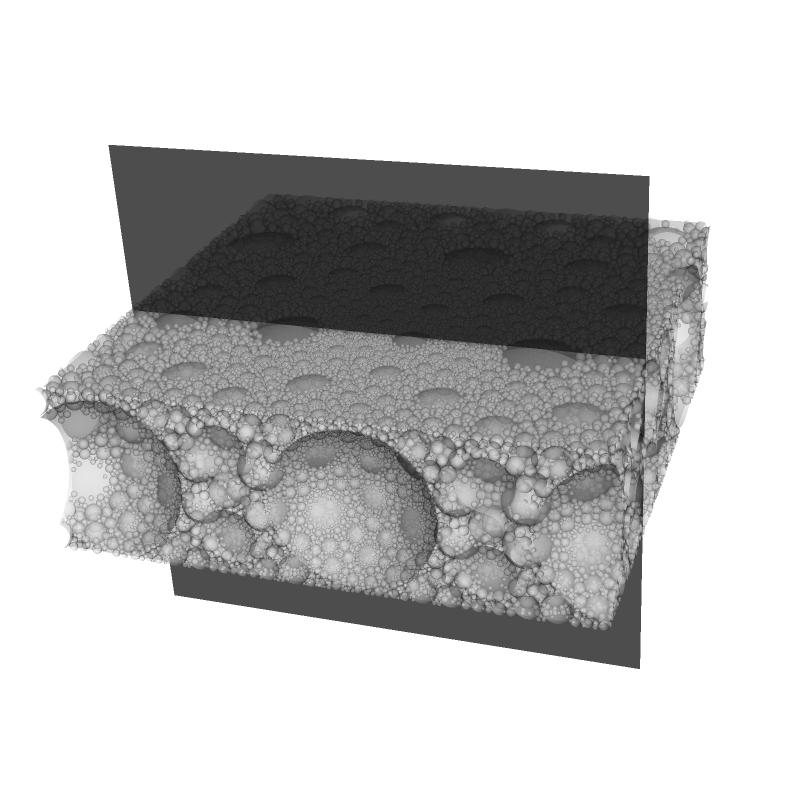}};
					
				\end{scope}
				\begin{scope}[xshift=-0cm]
					\node at (0,0) {\includegraphics[trim=0 0 40 60,clip,align=c,width=0.35\textwidth]{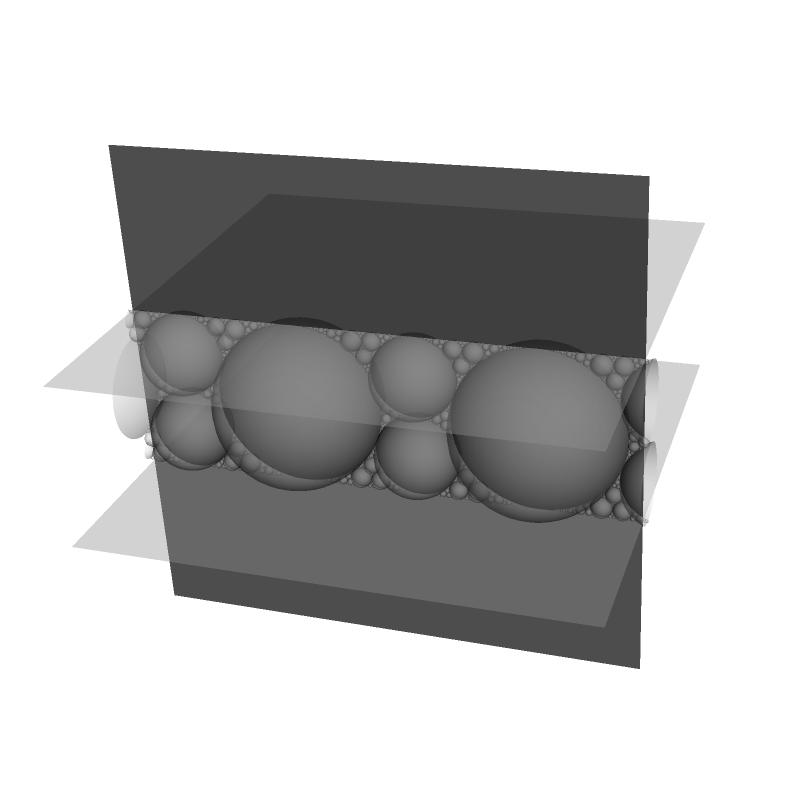}};

				\end{scope}
				\begin{scope}[xshift=5.cm]
					\node at (0,0) {\includegraphics[align=c,width=0.25\textwidth]{img/strips/34packing.pdf}};
					
				\end{scope}
			\end{scope}

		\end{tikzpicture}
		
		
		\caption{
The octahedral Apollonian section $\mathscr S_{\{3,4,3\}}^{\{3,4\}}$. 
(From left to right) The 24-cell crystallographic packing $\mathscr{P}_{\{3,4,3\}}$ with a cutting sphere $\Sigma_{\{3,4,3\}}^{\{3,4\}}$, the octahedral Apollonian section $\mathscr S_{\{3,4,3\}}^{\{3,4\}}$ with $\Sigma_{\{3,4,3\}}^{\{3,4\}}$, and the octahedral crystallographic packing $\mathscr{P}_{\{3,4\}}$. 		
		}
		\label{fig:sectionsgeo343_34}
	\end{figure}

\begin{figure}[H]
	\centering
	
	\begin{tikzpicture}
		\begin{scope}
			\begin{scope}[xshift=-5.5cm]
				\node at (0,0) {\includegraphics[trim=0 0 40 0,clip,align=c,width=0.35\textwidth]{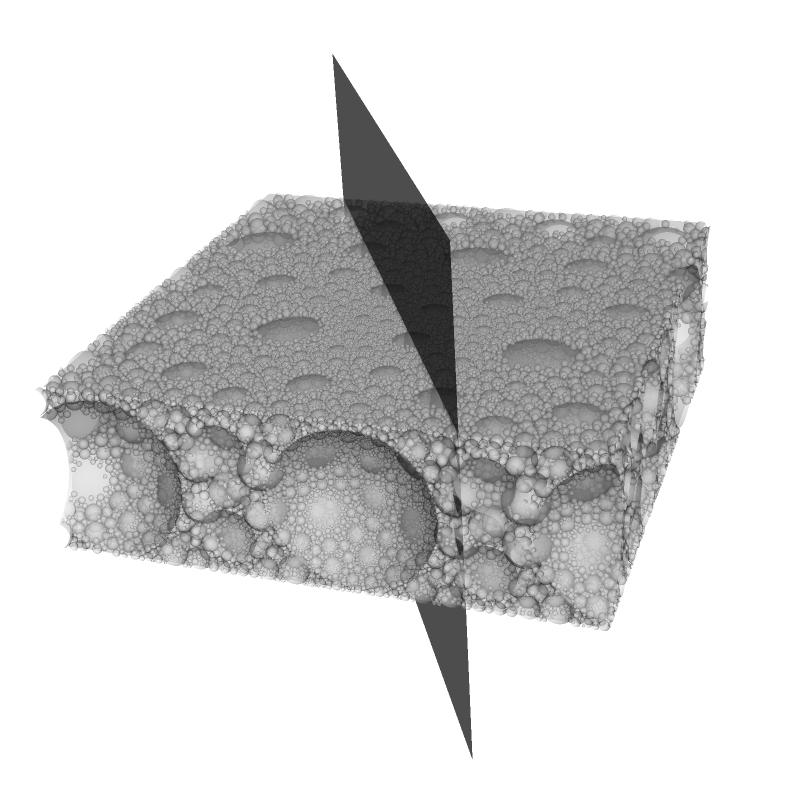}};
				
			\end{scope}
			\begin{scope}[xshift=-0cm]
				\node at (0,0) {\includegraphics[trim=0 0 40 0,clip,align=c,width=0.35\textwidth]{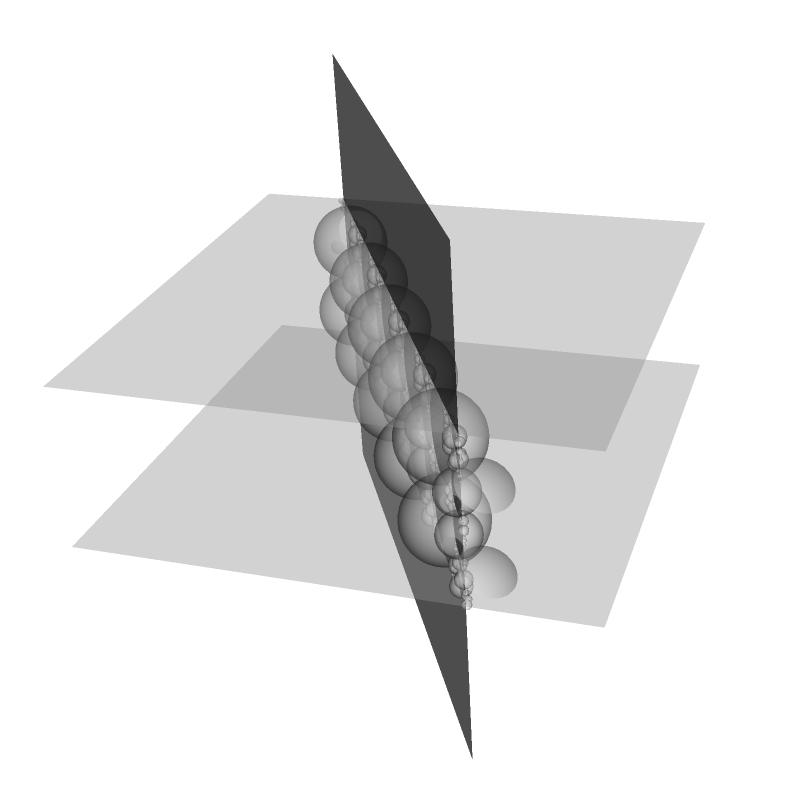}};

			\end{scope}
			\begin{scope}[xshift=5.cm]
				\node at (0,0) {\includegraphics[align=c,width=0.25\textwidth]{img/strips/43packing.pdf}};

			\end{scope}
		\end{scope}

	\end{tikzpicture}
	
	\vspace{-.5cm}
	
	\caption{
(From left to right) The 24-cell crystallographic packing $\mathscr{P}_{\{3,4,3\}}$ with a cutting sphere $\Sigma_{\{3,4,3\}}^{\{4,3\}}$, the cubic Apollonian section $\mathscr S_{\{3,4,3\}}^{\{4,3\}}$ with $\Sigma_{\{3,4,3\}}^{\{4,3\}}$, and the cubic crystallographic packing $\mathscr{P}_{\{4,3\}}$. 		
	}
	\label{fig:sectionsgeo343_43}
\end{figure}

Every 24-cell crystallographic packing contains an octahedral Apollonian section $\mathscr{S}_{\{3,4,3\}}^{\{3,4\}}$ and a cubic Apollonian section $\mathscr{S}_{\{3,4,3\}}^{\{4,3\}}$. In the strip packing, the cutting spheres are  $\Sigma_{\{3,4,3\}}^{\{3,4\}}=S_f$ and the plane $\Sigma_{\{3,4,3\}}^{\{4,3\}}$ passing through $(\frac1{\sqrt2},0,0)$ with normal vector $(\sqrt{3},1,0)$ (see Figures \ref{fig:symaposec343}, \ref{fig:sectionsgeo343_34}, \ref{fig:sectionsgeo343_43}). The homomorphisms between the full symmetry groups are described below.
		\begin{align}\label{eq:hom343}
			\begin{tabular}{cccccc}
				&$\Gamma_{\{3,4\}}$&$\xrightarrow{\phi_{\{3,4,3\}}^{\{3,4\}}}$&$\Gamma_{\{3,4,3\}}$\\
				&$r_1$&$\longmapsto$& $r_1$\\
				&$r_2$&$\longmapsto$& $r_2$\\
				&$r_3$&$\longmapsto$& $r_3$\\
				&$s_f$&$\longmapsto$& $(r_4s_f)^3$\\
			\end{tabular}&&
			\begin{tabular}{cccccc}
				&$\Gamma_{\{4,3\}}$&$\xrightarrow{\phi_{\{3,4,3\}}^{\{4,3\}}}$&$\Gamma_{\{3,4,3\}}$\\
				&$r_1$&$\longmapsto$& $r_1$\\
				&$r_2$&$\longmapsto$& $r_2r_3r_2$\\
				&$r_3$&$\longmapsto$& $r_4$\\
				&$s_f$&$\longmapsto$& $r_3s_f$\\
			\end{tabular}
		\end{align}

		\vspace{-1cm}

\subsection{600-cell $\{3,3,5\}$}

\begin{figure}[H]

				\includegraphics[align=c,width=.22\textwidth]{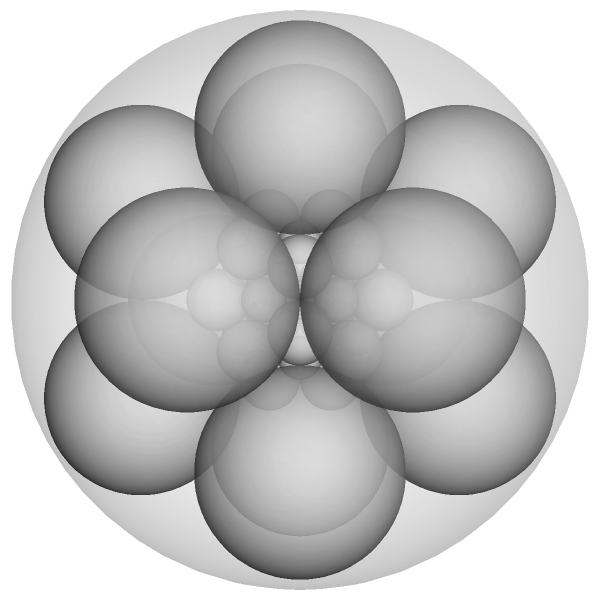} \hspace{.1cm}  
				\includegraphics[align=c,width=.22\textwidth]{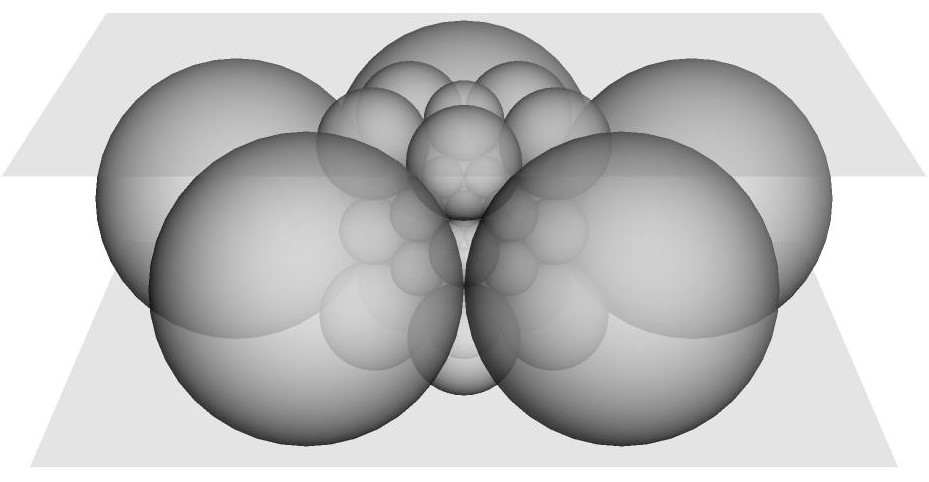} \hspace{.1cm} 
				\includegraphics[align=c,width=.22\textwidth]{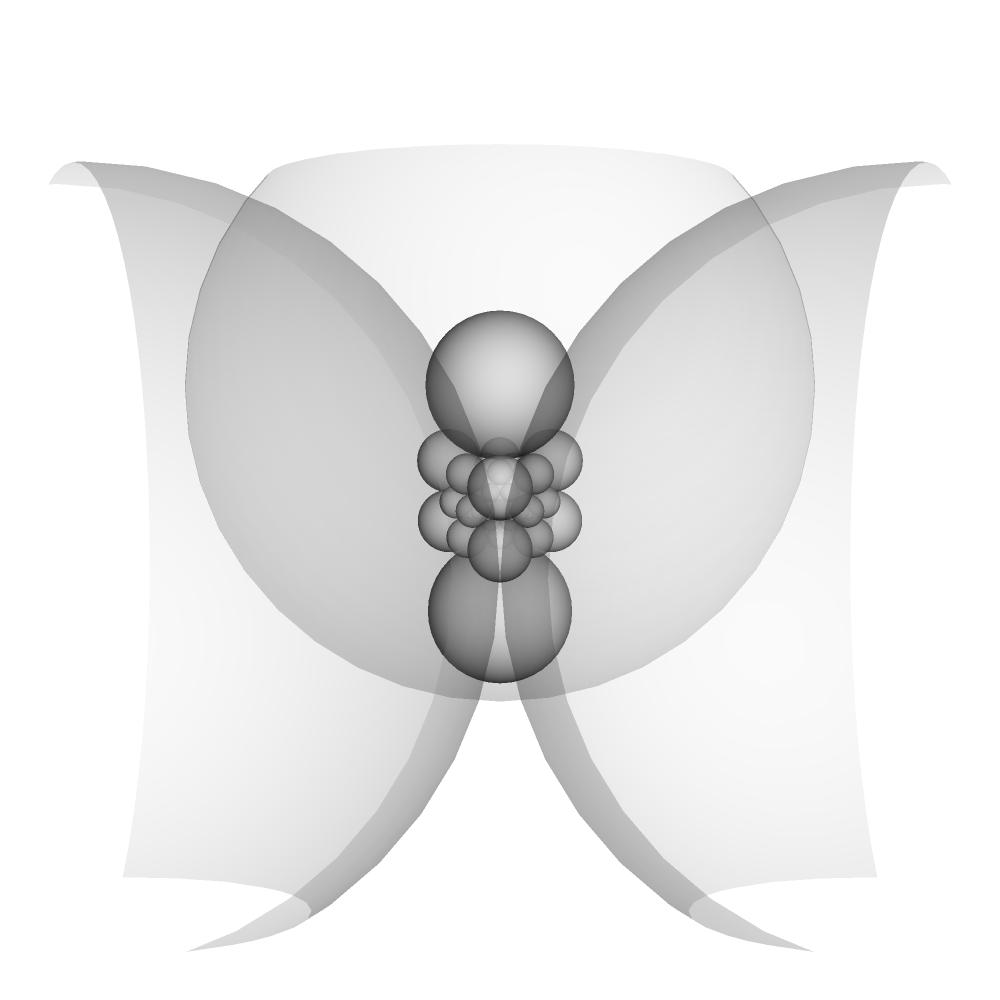} \hspace{.1cm} 
				\includegraphics[align=c,width=.22\textwidth]{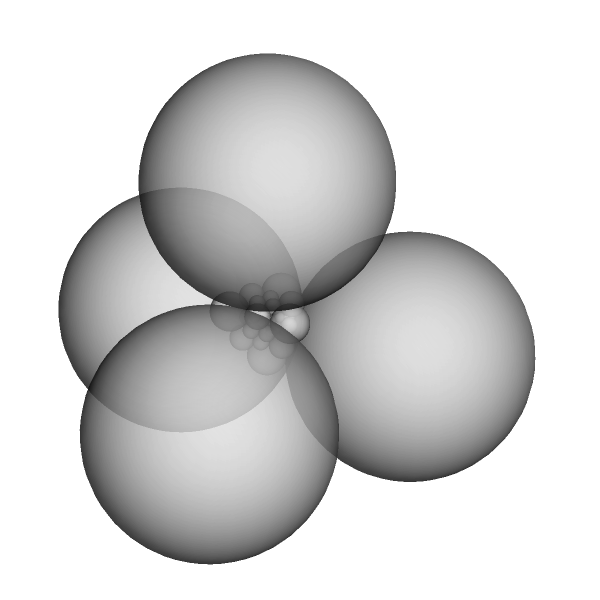}
				\vspace{-.4cm}
				\caption{(From left to right) Arrangement projections of a vertex-centered, edge-centered, ridge-centered and facet-centered canonical 600-cell.}
				\label{fig:centered600}
			\end{figure}

			In Figure \ref{fig:centered600}, we show four polytopal sphere packings obtained by the arrangement projections of face-centered canonical 600-cells.
			\vspace{-.4cm}

\vspace{-.3cm}
			The full symmetry group $\Gamma_{\{3,3,5\}}$ is isomorphic to $\langle \mathbf{R}_1,\mathbf R_2,\mathbf{R}_3,\mathbf R_4,\mathbf S_f\rangle <\mathrm{SL}_5(\mathbb Z[\varphi])$ where
		\begin{equation}
			\begin{gathered}
	\mathbf{R}_1= \left(
	\begin{array}{ccccc}
		& 1 &  &  &  \\
		1 &  &  &  &  \\
		&  & 1 &  &  \\
		&  &  & 1 &  \\
		\varphi  & -\varphi  &  &  & 1 \\
	\end{array}
	\right)\quad
	\mathbf{R}_2= \left(
	\begin{array}{ccccc}
		1 &  &  &  &  \\
		&  & 1 &  &  \\
		& 1 &  &  &  \\
		&  &  & 1 &  \\
		&  &  &  & 1 \\
	\end{array}
	\right)\quad
	\mathbf{R}_3= \left(
	\begin{array}{ccccc}
		1 &  &  &  &  \\
		& 1 &  &  &  \\
		&  &  & 1 &  \\
		&  & 1 &  &  \\
		&  &  &  & 1 \\
	\end{array}
	\right)\\
	\mathbf{R}_4=
	\left(
	\begin{array}{ccccc}
		1 &  &  &  &  \\
		& 1 &  &  &  \\
		&  & 1 &  &  \\
		\varphi  &  &  & -\varphi  & 1 \\
		1 &  &  & -\varphi & \varphi  \\
	\end{array}
	\right)\quad
	\mathbf S_f=
	\left(
	\begin{array}{ccccc}
		1 &  &  &  &  \\
		& 1 &  &  &  \\
		&  & 1 &  &  \\
		&  &  & 1 &  \\
		-\varphi^{-1} & \varphi ^2 & \varphi ^2 & \varphi ^2 & -1 \\
	\end{array}
	\right)
\end{gathered}
\end{equation}

Any fundamental bend vector $\mathbf b=(b_1,b_2,b_3,b_4,b_5)^T$ of a 600-cell sphere packing 
satisfies the quadratic equation $\mathbf{b}^T\mathbf Q_{\{3,3,5\}}\mathbf{b}=0$ for the bisymmetric matrix
\begin{align}\label{eq:descartes334}
	\mathbf Q_{\{3,3,5\}}=
	\left(
	\begin{array}{ccccc}
		2 & -\varphi ^2 & -\varphi ^2 & -\varphi ^2 & \varphi^{-1} \\
		\ast & 2 \varphi ^2 & \varphi  & \varphi  & \ast \\
		\ast & \ast & 2 \varphi ^2 & \ast  & \ast \\
		\ast & \ast  & \ast  & \ast & \ast \\
		\ast  & \ast & \ast& \ast & \ast\\
	\end{array}
	\right)
\end{align} 
				
\begin{figure}[H]
				\centering
				\begin{tikzpicture}[scale=2.5] 
					\node at (-3.2,0) {	\includegraphics[clip,trim=0 60 0 60,align=c,width=.49\textwidth]{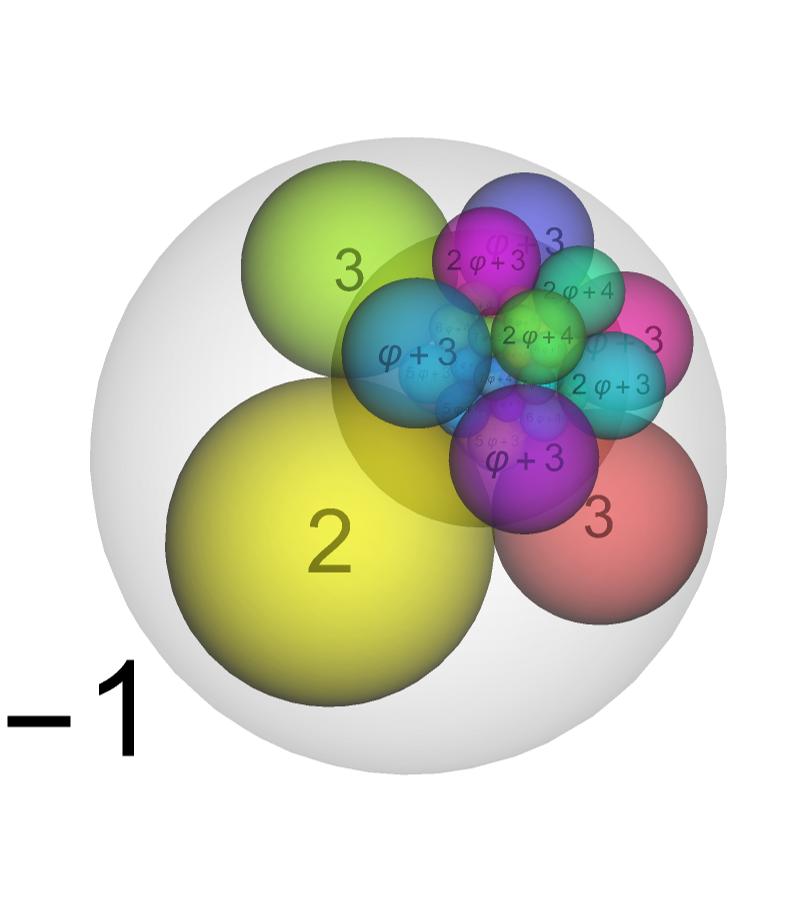}};
					
					\node at (0,0) {	\includegraphics[align=c,width=.4\textwidth]{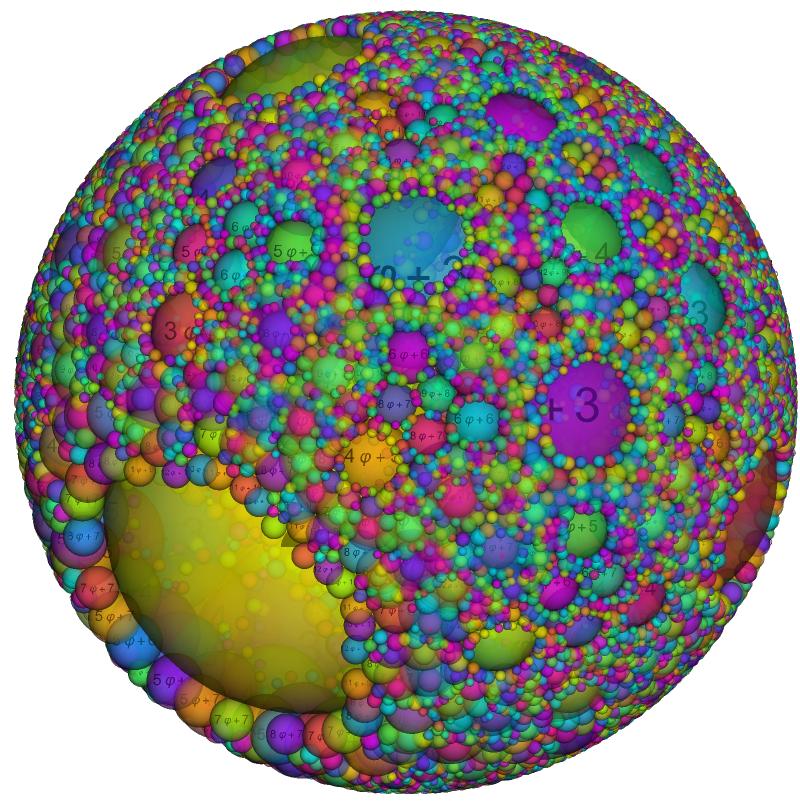}} ;
				\end{tikzpicture}
				\caption{
					A 600-cell sphere packing satisfying the integrality condition (left) and the corresponding $\mathbb Z[\varphi]$-integral Apollonian arrangement $\mathscr P_{\{3,3,5\}}(-1,2,2,3)$ (right).					
				}
				\label{fig:apo335}
			\end{figure}

\begin{figure}[H]
			\centering
			
			\begin{tikzpicture}
				
				\begin{scope}
					
					\begin{scope}[xshift=-5.5cm]
						\node at (0,0) {	\begin{tikzpicture}[scale=1.4] 		
		
		\newcommand\xmin{1}
\newcommand\xmax{2.414}

\newcommand\ymin{1}
\newcommand\ymax{2}
		
		\draw[very thin] (-\xmin,1) -- (\xmax,1);
		\draw[very thin] (0,0.5) circle (.5cm);
		\draw[very thin] (1.,0.5) circle (.5cm) ;
		
		\draw[very thick,red] (-\xmin,0) -- (\xmax,0) node[pos=.95, below] {$S_v$};
		
\draw[blue,thick] (-\xmin,0.5) -- (\xmax,.5)	node[pos=.95, above] {$R_1$};
\draw[thick,blue] (0,0) circle (1cm) node [xshift=-1.3cm,yshift=-1.3cm] {$R_2$};
\draw[blue,thick] (0.5,-\ymin) -- (0.5,\ymax) node[pos=.95,right] {$R_3$};
\draw[blue,thick] (0,-\ymin) -- (0,\ymax) node[pos=.95,left] {$S_f$};

	\end{tikzpicture}};

					\end{scope}
					\begin{scope}[xshift=0cm]
						\node at (0,0) {	\begin{tikzpicture}[scale=1.4] 		
		
		\newcommand\xmin{1}
\newcommand\xmax{2.414}

\newcommand\ymin{1}
\newcommand\ymax{2}
		
\draw[very thin] (-\xmin,1.) -- (\xmax,1);
\draw[very thin, draw=black] (0,0.5000) circle (0.5000cm);
\draw[very thin, draw=black] (1.618,0.5000) circle (0.5000cm);
\draw[very thin, draw=black] (1.000,0.8090) circle (0.1910cm);
\draw[very thin, draw=black] (0.6180,0.8090) circle (0.1910cm);
\draw[very thin, draw=black] (1.000,0.1910) circle (0.1910cm);
\draw[very thin, draw=black] (0.6180,0.1910) circle (0.1910cm);
\draw[very thin, draw=black] (0.6180,0.5000) circle (0.1180cm);
\draw[very thin, draw=black] (1.000,0.5000) circle (0.1180cm);
\draw[very thin, draw=black] (0.8090,0.4045) circle (0.09549cm);
\draw[very thin, draw=black] (0.8090,0.5955) circle (0.09549cm);

		\draw[very thick,red] (-\xmin,0) -- (\xmax,0) node[pos=.95, below] {$S_v$};
		
\draw[blue,thick] (-\xmin,0.5) -- (\xmax,.5)	node[pos=.95, above] {$R_1$};
\draw[thick,blue] (0,0) circle (1cm) node [xshift=-1.3cm,yshift=-1.3cm] {$R_2$};
\draw[blue,thick] (0.809,-\ymin) -- (0.809,\ymax) node[pos=.95,right] {$R_3$};
\draw[blue,thick] (0,-\ymin) -- (0,\ymax) node[pos=.95,left] {$S_f$};

	\end{tikzpicture}};
						
					\end{scope}
					
					\begin{scope}[scale=2.,xshift=3cm] 
						\node {	\includegraphics[clip,trim=0 0 0 20,align=c,width=0.4\textwidth]{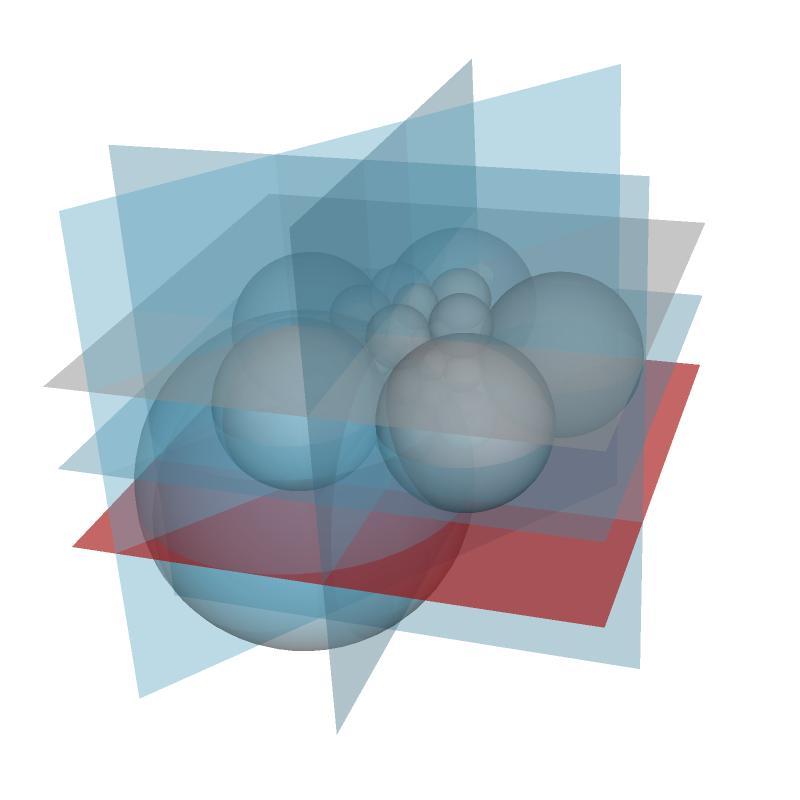} };
						
						\node at (1.25,0.15) [red] {$S_v$};
						\node at (-1.15,1.09) [blue] {$S_f$};
						\node at (1.3,.4) [blue] {$R_1$};
						\node at (-.55,-1.1) [blue] {$R_2$};
						\node at (.35,1.4) [blue] {$R_3$};
						\node at (.96,1.36) [blue] {$R_4$};
						
					\end{scope}
				\end{scope}
							
			\end{tikzpicture}
			
		\vspace{-1cm}
		\caption{
			The strip packing with the walls of the fundamental symmetries of the tetrahedron (left), the icosahedron (center) and the 600-cell (right).
		}
		\label{fig:symaposec335}
	\end{figure}

			The integrality condition of Corollary \ref{cor:integrality} states that if  $b_1,b_2,b_3,b_4,\sqrt{\Delta_{\{3,3,5\}}}\in\mathbb Z[\varphi]$ where
			\begin{align}\label{eq:inthico}
				\Delta_{\{3,3,5\}}=(2+\varphi)\left((b_1+b_2+b_3+b_4)^2-2(b_1^2+b_2^2+b_3^2+b_4^2)\right)
			\end{align}
			then the 600-cell Apollonian arrangement $\mathscr P_{\{3,3,5\}}(b_1,b_2,b_3,b_4)$  is $\mathbb Z[\varphi]$-integral (see Figure \ref{fig:apo335}).
			\vspace{-.3cm}			
			\vspace{-.2cm}
			Every Apollonian arrangement of the 600-cell contains a tetrahedral Apollonian section $\mathscr{S}^{\{3,3\}}_{\{3,3,5\}}$ and an icosahedral Apollonian section $\mathscr{P}^{\{3,5\}}_{\{3,3,5\}}$. In the strip packing, the cutting spheres are the plane $\Sigma_{\{3,3,5\}}^{\{3,3\}}=S_f$ and the plane $\Sigma_{\{3,3,5\}}^{\{3,5\}}$ passing through the origin with normal vector $(-\sqrt{7-4\varphi},1,0)$ (see Fig. \ref{fig:symaposec335}). 	The homomorphisms between the full symmetry groups are given by

			\begin{align}\label{eq:hom335}
				\begin{tabular}{cccccc}
					&$\Gamma_{\{3,3\}}$&
					$\xrightarrow{\phi_{\{3,3,5\}}^{\{3,3\}}}$
					&$\Gamma_{\{3,3,5\}}$\\
					&$r_1$&$\longmapsto$& $r_1$\\
					&$r_2$&$\longmapsto$& $r_2$\\
					&$r_3$&$\longmapsto$& $r_3$\\
					&$s_f$&$\longmapsto$& $(r_4s_f)^5$\\
				\end{tabular}&&
				\begin{tabular}{cccccc}
					&$\Gamma_{\{3,5\}}$&
					$\xrightarrow{\phi_{\{3,3,5\}}^{\{3,5\}}}$
					&$\Gamma_{\{3,3,5\}}$\\
					&$r_1$&$\longmapsto$& $r_1$\\
					&$r_2$&$\longmapsto$& $r_2$\\
					&$r_3$&$\longmapsto$& $r_3r_4r_3$\\
					&$s_f$&$\longmapsto$& $(r_4s_f)^5$\\
				\end{tabular}
			\end{align}
			\vspace{-.5cm}
					\begin{figure}[H]
	\centering
	
	\begin{tikzpicture}
		\begin{scope}
			\begin{scope}[xshift=-5.5cm]
				\node at (0,0) {\includegraphics[trim=0 0 40 20,clip,align=c,width=0.35\textwidth]{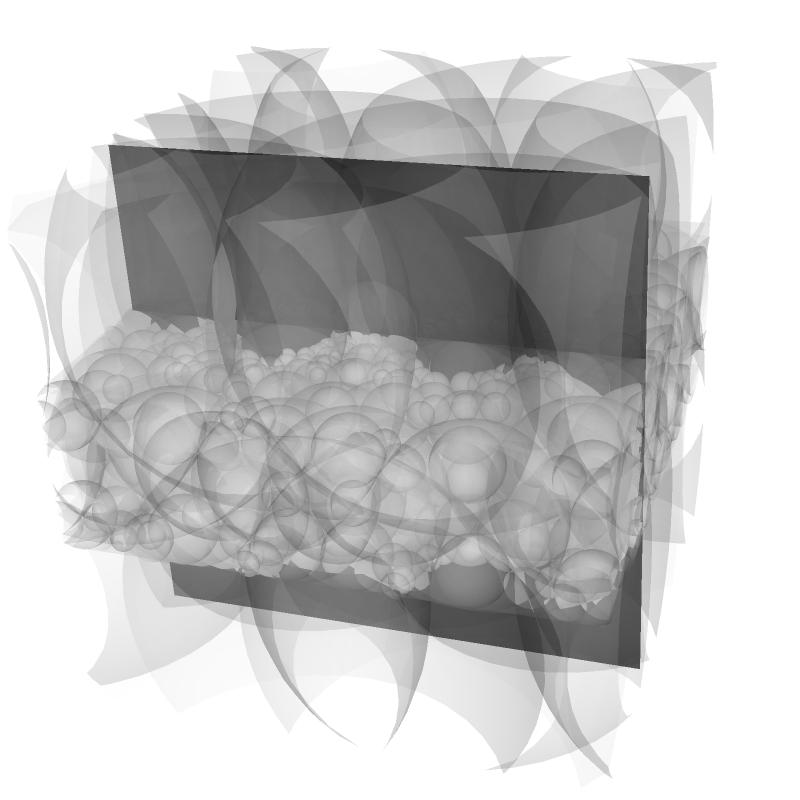}};
				
			\end{scope}
			\begin{scope}[xshift=-0cm]
				\node at (0,0) {\includegraphics[trim=0 0 40 20,clip,align=c,width=0.35\textwidth]{img/sections/aposection33_334b}};
				
			\end{scope}
			\begin{scope}[xshift=5.cm]
				\node at (0,0) {\includegraphics[align=c,width=0.25\textwidth]{img/strips/33packing.pdf}};
				
			\end{scope}
		\end{scope}

	\end{tikzpicture}
	
			\vspace{-.5cm}
\caption{
(From left to right) The 600-cell Apollonian arrangement $\mathscr{P}_{\{3,3,5\}}$ with a cutting sphere $\Sigma_{\{3,3,5\}}^{\{3,3\}}$, the tetrahedral Apollonian section $\mathscr S_{\{3,3,5\}}^{\{3,3\}}$ with $\Sigma_{\{3,3,5\}}^{\{3,3\}}$, and the tetrahedral crystallographic packing $\mathscr{P}_{\{3,3\}}$. 		
}

	\label{fig:sectionsgeo335_33}
\end{figure}
		\vspace{-.5cm}
					\begin{figure}[H]
	\centering
	
	\begin{tikzpicture}
		\begin{scope}
			\begin{scope}[xshift=-5.5cm]
				\node at (0,0) {\includegraphics[trim=0 0 40 20,clip,align=c,width=0.35\textwidth]{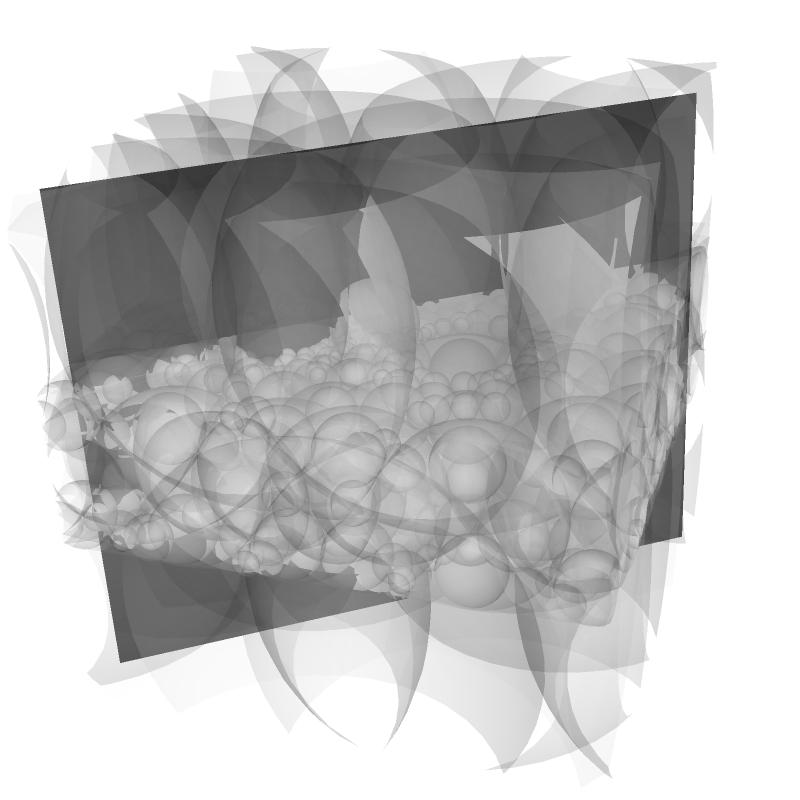}};
				
			\end{scope}
			\begin{scope}[xshift=-0cm]
				\node at (0,0) {\includegraphics[trim=0 0 40 20,clip,align=c,width=0.35\textwidth]{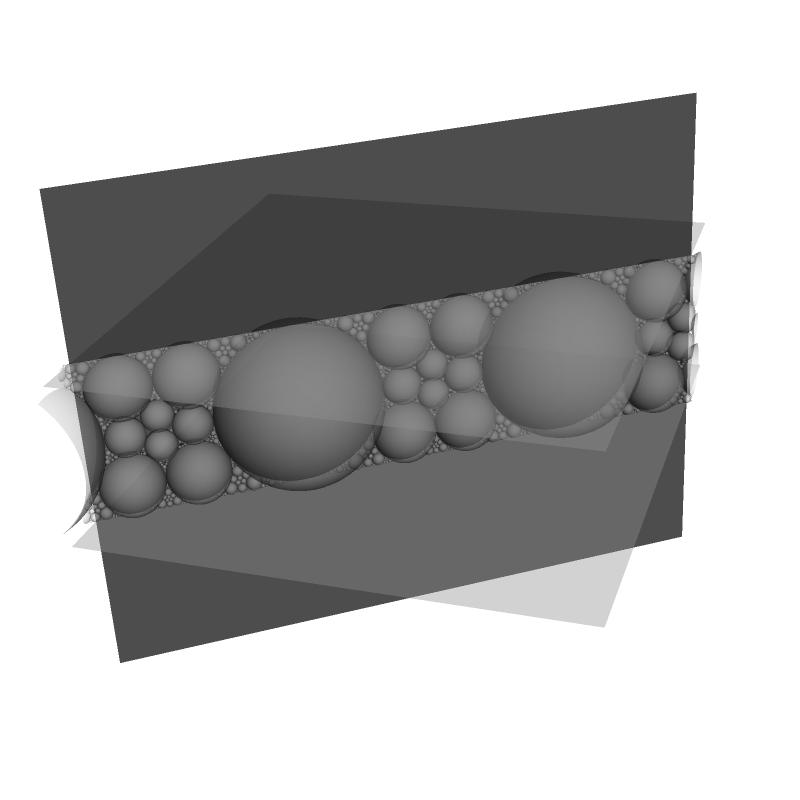}};
				
			\end{scope}
			\begin{scope}[xshift=5.cm]
				\node at (0,0) {\includegraphics[align=c,width=0.25\textwidth]{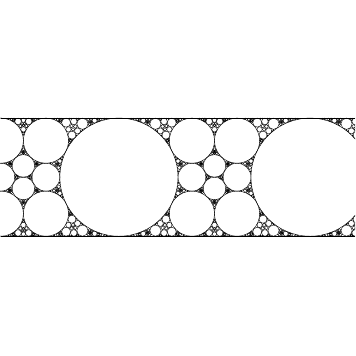}};
				
			\end{scope}
		\end{scope}

	\end{tikzpicture}
				\vspace{-.5cm}
	
	\caption{
(From left to right) The 600-cell Apollonian arrangement $\mathscr{P}_{\{3,3,5\}}$ with a cutting sphere $\Sigma_{\{3,3,5\}}^{\{3,5\}}$, the icosahedral Apollonian section $\mathscr S_{\{3,3,5\}}^{\{3,5\}}$ with $\Sigma_{\{3,3,5\}}^{\{3,5\}}$, and the icosahedral crystallographic packing $\mathscr{P}_{\{3,5\}}$. 		
}
	\label{fig:sectionsgeo335_35}
\end{figure}

\subsection{120-cell $\{5,3,3\}$}\label{sec:533}

\begin{figure}[H]
				\includegraphics[align=c,width=.22\textwidth]{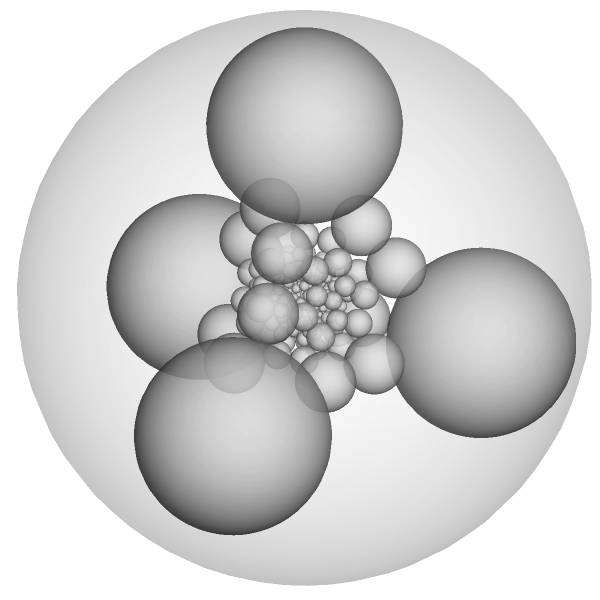} \hspace{.1cm}  
				\includegraphics[align=c,width=.22\textwidth]{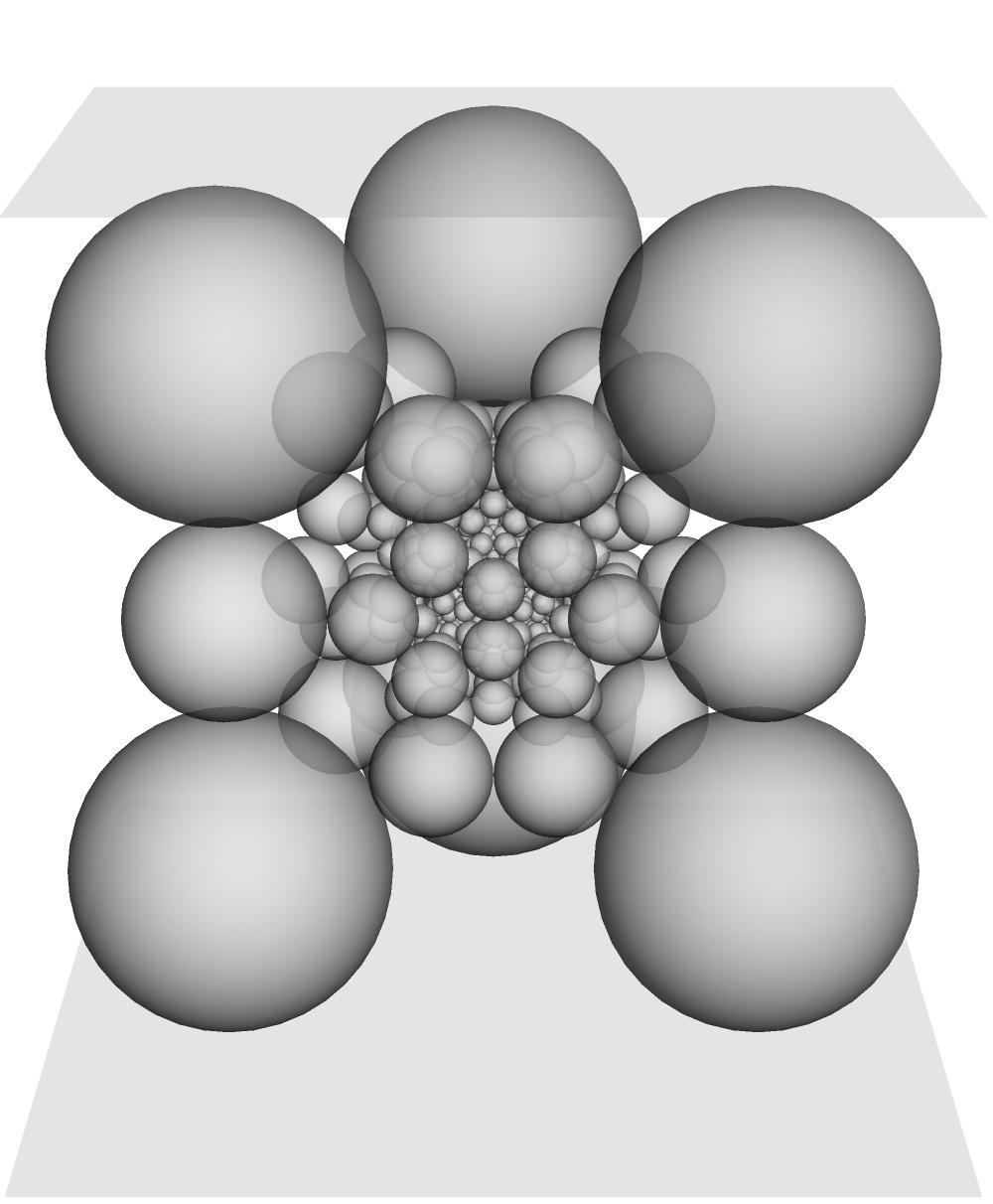} \hspace{.1cm} 
				\includegraphics[align=c,width=.22\textwidth]{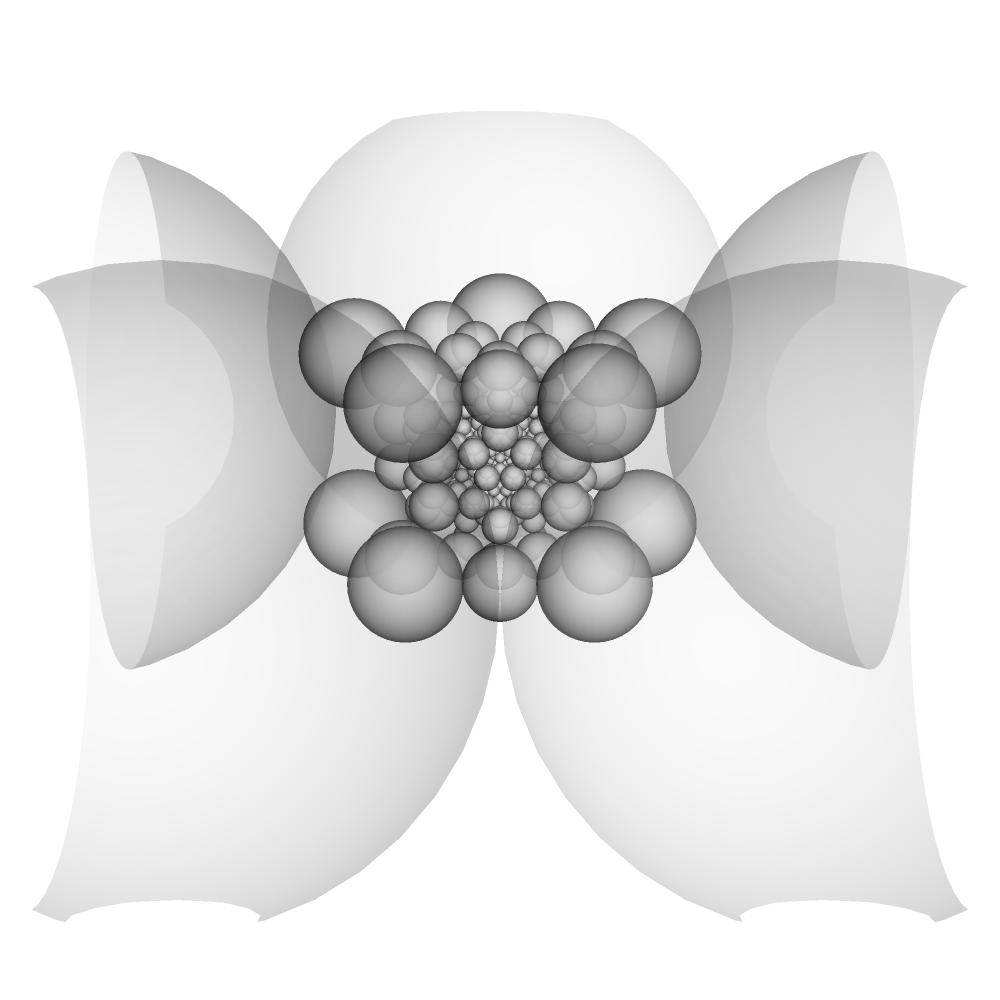} \hspace{.1cm} 
				\includegraphics[align=c,width=.22\textwidth]{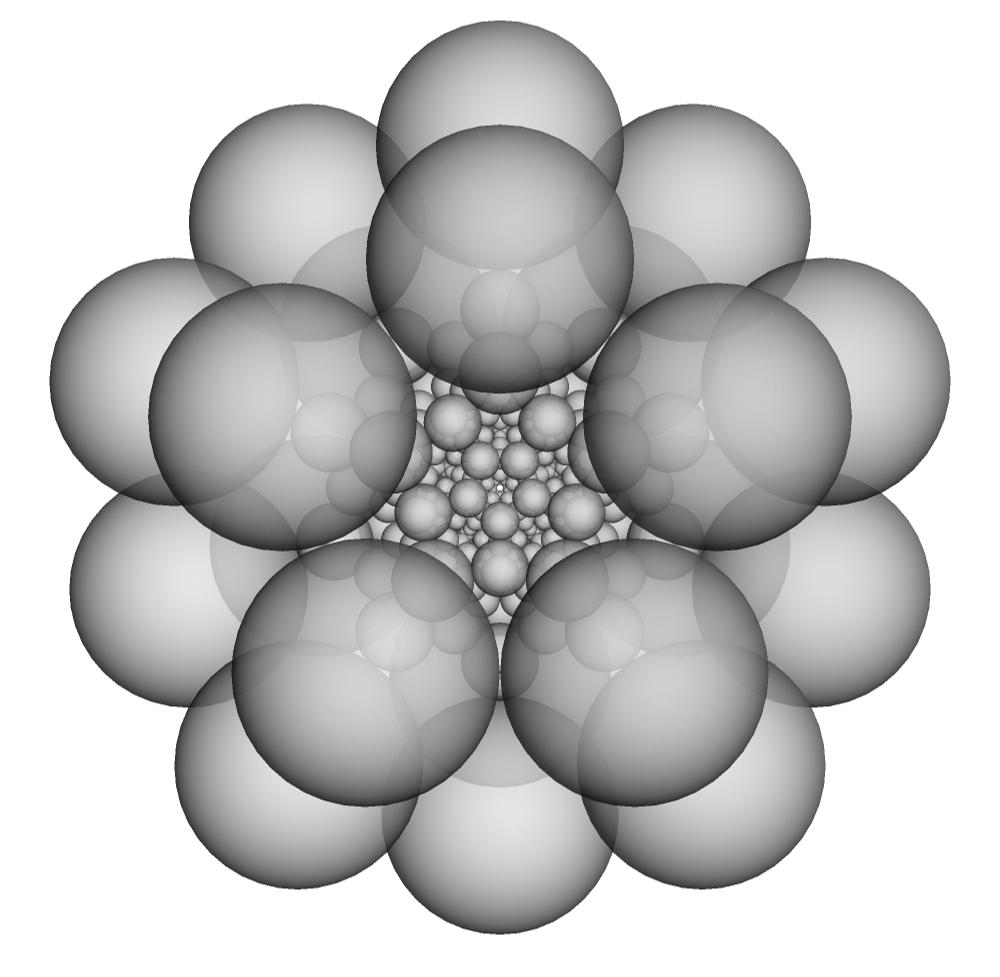}\\
				\caption{(From left to right) Arrangement projections of a vertex-centered, edge-centered, ridge-centered and facet-centered canonical 120-cell.}
				\label{fig:centered120}
	\end{figure}

In Figure \ref{fig:centered120}, we show four polytopal sphere packings obtained by the arrangement projections of face-centered canonical 120-cells.

				The full symmetry group $\Gamma_{\{5,3,3\}}$ is isomorphic to $\langle \mathbf{R}_1,\mathbf R_2,\mathbf{R}_3,\mathbf R_4,\mathbf S_f\rangle <\mathrm{SL}_5(\mathbb Z[\varphi])$ where
				\begin{equation}
				\begin{gathered}
	\mathbf{R}_1=
	\left(
	\begin{array}{ccccc}
		& 1 &  &  &  \\
		1 &  &  &  &  \\
		\varphi  & -\varphi  & 1 &  &  \\
		\varphi ^2 & -\varphi ^2 &  & 1 &  \\
		\varphi ^3 & -\varphi ^3 &  &  & 1 \\
	\end{array}
	\right)\quad
	\mathbf{R}_2= \left(
	\begin{array}{ccccc}
		1 &  &  &  &  \\
		\varphi  & -\varphi  & 1 &  &  \\
		1 & -\varphi  & \varphi  &  &  \\
		\varphi  & -\varphi ^2 & 1 & 1 &  \\
		\varphi ^2 & -\varphi ^3 & \varphi  &  & 1 \\
	\end{array}
	\right)\quad
	\mathbf{R}_3= \left(
	\begin{array}{ccccc}
		1 &  &  &  &  \\
		& 1 &  &  &  \\
		& \varphi  & -\varphi  & 1 &  \\
		& 1 & -\varphi  & \varphi  &  \\
		& \varphi  & -\varphi ^2 & 1 & 1 \\
	\end{array}
	\right)\\
	\mathbf{R}_4=
	\left(
	\begin{array}{ccccc}
		1 &  &  &  &  \\
		& 1 &  &  &  \\
		&  & 1 &  &  \\
		&  & \varphi  & -\varphi  & 1 \\
		&  & 1 & -\varphi & \varphi  \\
	\end{array}
	\right)\quad
	\mathbf S_f=
	\left(
	\begin{array}{ccccc}
		1 &  &  &  &  \\
		& 1 &  &  &  \\
		&  & 1 &  &  \\
		&  &  & 1 &  \\
		1 & -\varphi^{-1}& -\varphi^{-1} & \varphi ^3 & -1 \\
	\end{array}
	\right)
\end{gathered}
\end{equation}	
				
Any fundamental bend vector $\mathbf b=(b_1,b_2,b_3,b_4,b_5)^T$ of a 120-cell sphere packing 
satisfies the quadratic equation $\mathbf{b}^T\mathbf Q_{\{5,3,3\}}\mathbf{b}=0$ for the bisymmetric matrix
\begin{align}\label{eq:descartes334}
	\mathbf Q_{\{5,3,3\}}=
	\left(
	\begin{array}{ccccc}
		2 & -\varphi ^3 & \varphi^{-1}  & \varphi^{-1} & -1 \\
		\ast & 4 \varphi ^2 & -2 \varphi ^2 & -\varphi^{-2} &\ast \\
		\ast &\ast  & 4 \varphi ^2 &\ast  & \ast \\
		\ast& \ast &\ast  &\ast  & \ast \\
		\ast	&\ast  & \ast & \ast &\ast  \\
	\end{array}
	\right)
\end{align}

\begin{figure}[H]
					\centering
					\begin{tikzpicture}[scale=2.5] 
						\node at (-3.2,0) {	\includegraphics[clip,trim=0 60 0 60,align=c,width=.49\textwidth]{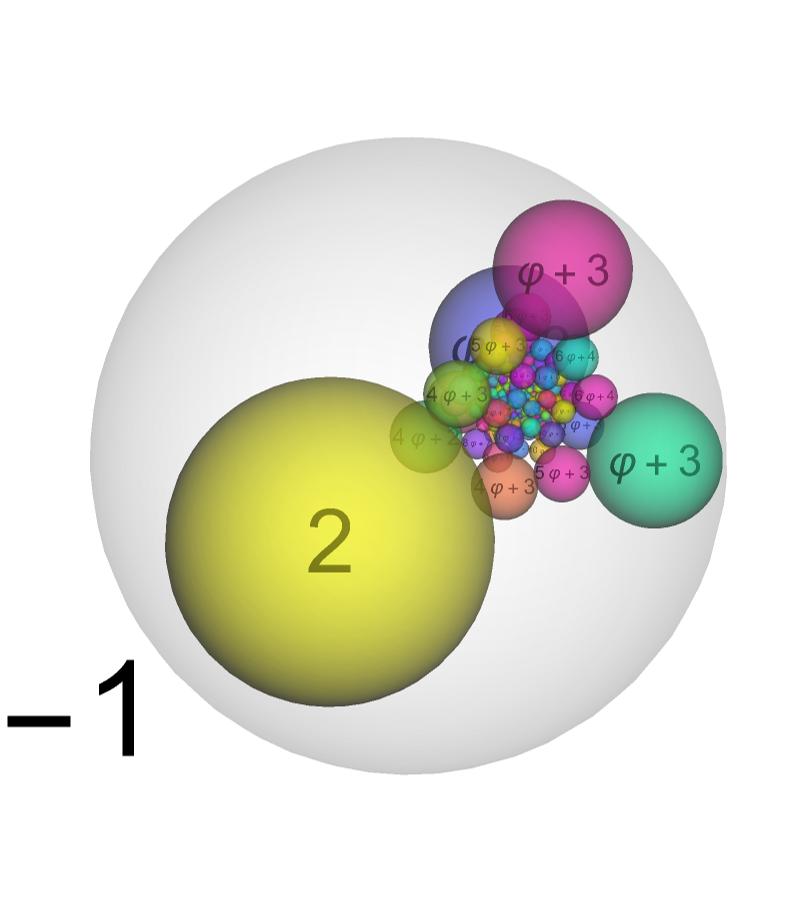}};
						
						\node at (0,0) {	\includegraphics[align=c,width=.4\textwidth]{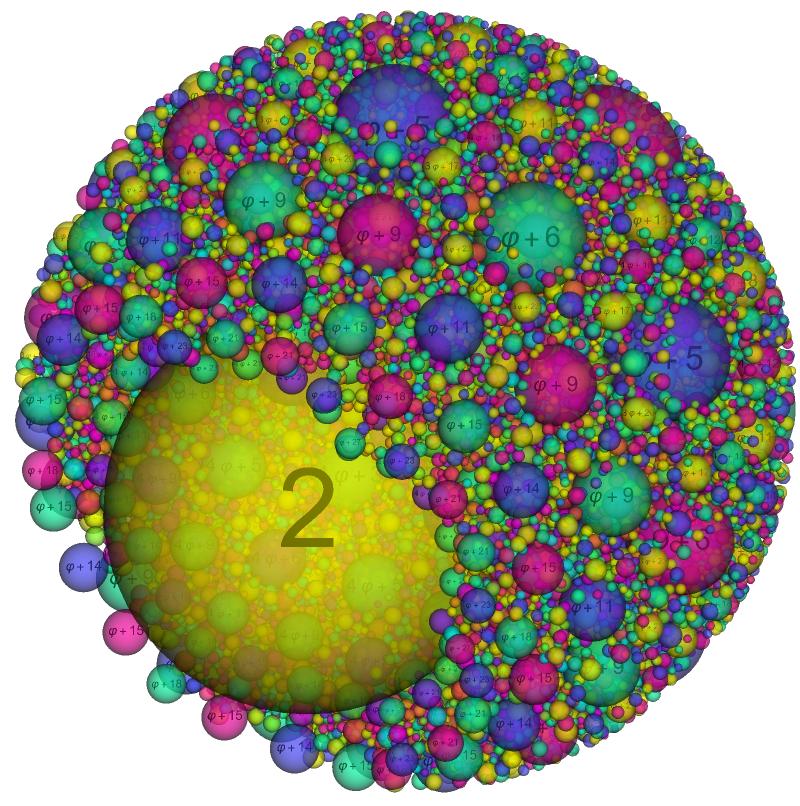}} ;

					\end{tikzpicture}
					
					\caption{
						A 120-cell packing satisfying the integrality condition (left) and the corresponding $\mathbb Z[\varphi]$-integral crystallographic packing $\mathscr P_{\{5,3,3\}}(-1,2,2+4\varphi,7+8\varphi)$ (right).
					}
					\label{fig:apo533}
				\end{figure}

				The integrality condition of Corollary \ref{cor:integrality} states that if  $b_1,b_2,b_3,b_4,\sqrt{\Delta_{\{5,3,3\}}}\in\mathbb Z[\varphi]$ where
				\begin{align}\label{eq:int120}
					\Delta_{\{5,3,3\}}=-3 Q_{\{5,3\}}(b_1,b_2,b_3,b_4)
				\end{align}
				and $ Q_{\{5,3\}}$ is one of the quadratic forms described in \eqref{eq:descartes3d},	then the 120-cell crystallographic packing $\mathscr P_{\{5,3,3\}}(b_1,b_2,b_3,b_4)$ is $\mathbb Z[\varphi]$-integral (see Figure \ref{fig:apo533}). 
				
\begin{figure}[H]
					\centering
					\begin{tabular}{cc}
						&\begin{tikzpicture}[scale=2] 
							\node at (-3,0) {		\begin{tikzpicture}[scale=1.4] 		
		
		\newcommand\xmin{1}
\newcommand\xmax{2.414}

\newcommand\ymin{1}
\newcommand\ymax{2}
		
\draw[very thin] (-\xmin,1.618) -- (\xmax,1.618);
\draw[very thin, draw=black] (0,0.3090) circle (0.3090cm);
\draw[very thin, draw=black] (1.000,0.3090) circle (0.3090cm);
\draw[very thin, draw=black] (0,0.8090) circle (0.1910cm);
\draw[very thin, draw=black] (1.00,1.309) circle (0.3090cm);
\draw[very thin, draw=black] (1.000,0.8090) circle (0.1910cm);
\draw[very thin, draw=black] (0,1.309) circle (0.3090cm);
\draw[very thin, draw=black] (0.6180,0.5000) circle (0.1180cm);
\draw[very thin, draw=black] (0.3820,0.5000) circle (0.1180cm);
\draw[very thin, draw=black] (0.6180,1.118) circle (0.1180cm);
\draw[very thin, draw=black] (0.3820,1.118) circle (0.1180cm);
\draw[very thin, draw=black] (0.6180,0.6910) circle (0.07295cm);
\draw[very thin, draw=black] (0.2764,0.8090) circle (0.08541cm);
\draw[very thin, draw=black] (0.3820,0.6910) circle (0.07295cm);
\draw[very thin, draw=black] (0.7236,0.8090) circle (0.08541cm);
\draw[very thin, draw=black] (0.3820,0.9271) circle (0.07295cm);
\draw[very thin, draw=black] (0.6180,0.9271) circle (0.07295cm);
\draw[very thin, draw=black] (0.5000,0.7500) circle (0.05902cm);
\draw[very thin, draw=black] (0.5000,0.8680) circle (0.05902cm);

		\draw[very thick,red] (-\xmin,0) -- (\xmax,0) node[pos=.95, below] {$S_v$};
		
\draw[blue,thick] (-\xmin,0.809) -- (\xmax,.809)	node[pos=.95, above] {$R_1$};
\draw[thick,blue] (0,0) circle (1cm) node [xshift=-1.3cm,yshift=-1.3cm] {$R_2$};
\draw[blue,thick] (0.5,-\ymin) -- (0.5,\ymax) node[pos=.95,right] {$R_3$};
\draw[blue,thick] (0,-\ymin) -- (0,\ymax) node[pos=.95,left] {$S_f$};

	\end{tikzpicture}};
							
							\node {	\includegraphics[align=c,width=0.37\textwidth]{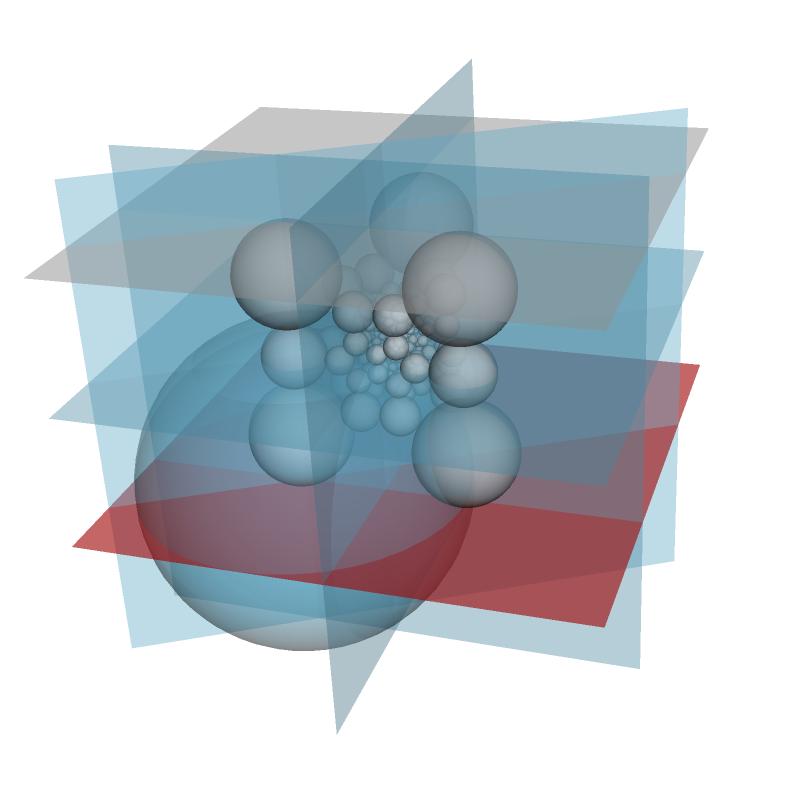} };
							\node at (1.25,0.15) [red] {$S_v$};
							\node at (-1.15,1.09) [blue] {$S_f$};
							\node at (1.25,.6) [blue] {$R_1$};
							\node at (-.6,-1.) [blue] {$R_2$};
							\node at (.35,1.4) [blue] {$R_3$};
							\node at (1.17,1.15) [blue] {$R_4$};

						\end{tikzpicture}
						
					\end{tabular}
						\vspace{-1cm}
					
					\caption{
						The strip packing with the walls of the fundamental symmetries of the dodecahedron (left) and the 120-cell (right).
					}
					\label{fig:symaposec533}
				\end{figure}

\begin{figure}[H]
					\centering
					
					\begin{tikzpicture}
						\begin{scope}
							\begin{scope}[xshift=-5.5cm]
								\node at (0,0) {\includegraphics[trim=0 0 40 40,clip,align=c,width=0.35\textwidth]{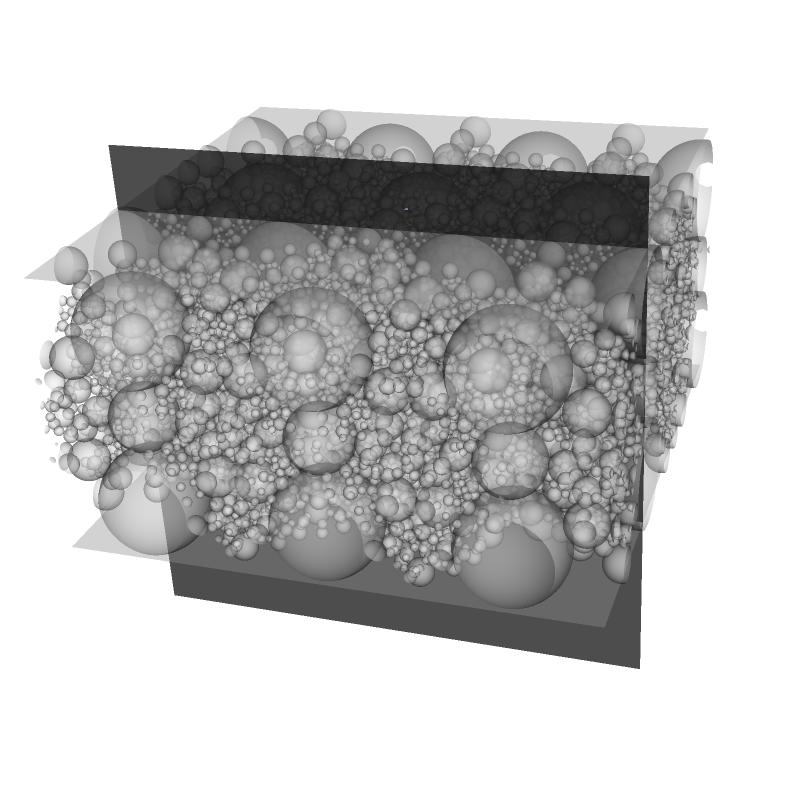}};
								
							\end{scope}
							\begin{scope}[xshift=-0cm]
								\node at (0,0) {\includegraphics[trim=0 0 40 40,clip,align=c,width=0.35\textwidth]{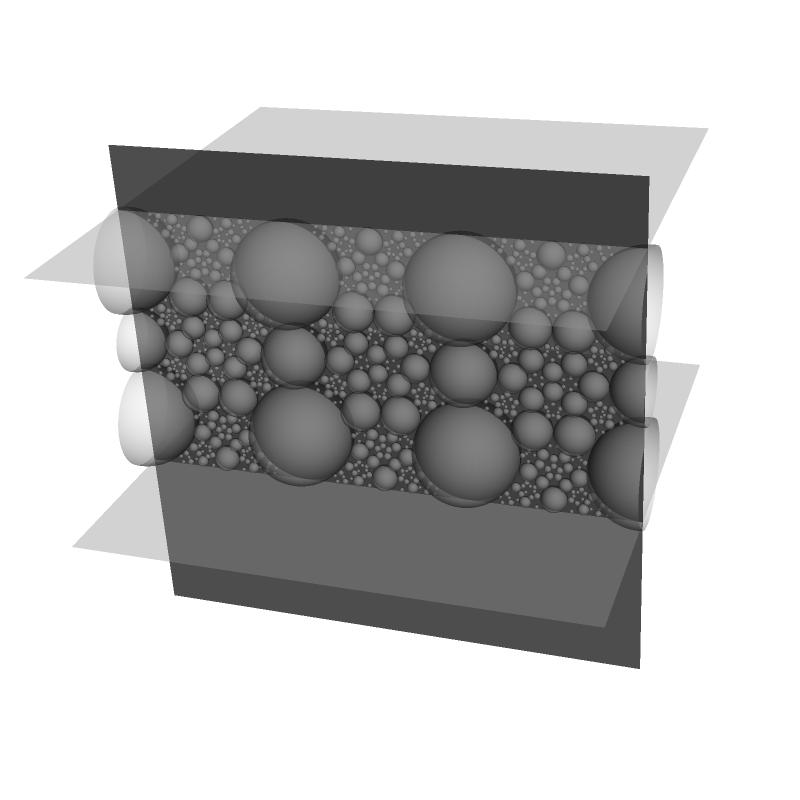}};
								
							\end{scope}
							\begin{scope}[xshift=5.cm]
								\node at (0,0) {\includegraphics[align=c,width=0.25\textwidth]{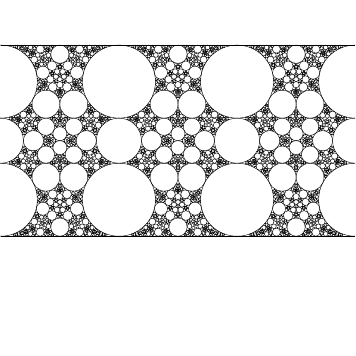}};
								
							\end{scope}
						\end{scope}

					\end{tikzpicture}
					
					\vspace{-1.2cm}
					
					\caption{
(From left to right) The 120-cell crystallographic packing $\mathscr{P}_{\{5,3,3\}}$ with a cutting sphere $\Sigma_{\{5,3,3\}}^{\{5,3\}}$, the dodecahedral Apollonian section $\mathscr S_{\{5,3,3\}}^{\{5,3\}}$ with $\Sigma_{\{5,3,3\}}^{\{5,3\}}$, and the dodecahedral crystallographic packing $\mathscr{P}_{\{5,3\}}$. 		
}	
					
					\label{fig:sectionsgeo533}
				\end{figure}

				Every 120-cell crystallographic packing contains a geometric dodecahedral Apollonian section $\mathscr{S}^{\{5,3\}}_{\{5,3,3\}}$. In the strip packing, the cutting sphere is the plane $\Sigma_{\{5,3,3\}}^{\{5,3\}}=S_f$ (see Figures \ref{fig:symaposec533}, \ref{fig:sectionsgeo533}).  The homomorphism
between the full symmetry groups is described below.

				\begin{align}\label{eq:hom533}
					\begin{tabular}{cccccc}
						&$\Gamma_{\{5,3\}}$&$\xrightarrow{\phi_{\{5,3,3\}}^{\{5,3\}}}$&$\Gamma_{\{5,3,3\}}$\\
						&$r_1$&$\longmapsto$& $r_1$\\
						&$r_2$&$\longmapsto$& $r_2$\\
						&$r_3$&$\longmapsto$& $r_3$\\
						&$s_f$&$\longmapsto$& $(r_4s_f)^3$\\
					\end{tabular}
				\end{align}
				
\end{document}